\documentclass[12pt]{article}

\usepackage{a4wide}   % longer lines
\usepackage{amsmath}   % for more basic mathematical symbols
\usepackage{amssymb}   % for more mathematical symbols
\usepackage{stmaryrd}  % for more mathematical symbols
\usepackage{latexsym}  % for more mathematical symbols,
\usepackage{amsxtra}   % for accents as superscripts
\usepackage{amstext}   % for text in formulas, accents, etc.
\usepackage{bm}        % boldface for non-latin letters
\usepackage{amsthm}    % for theorem-environements
\usepackage{amscd}     % for commutative diagrams
\usepackage[mathscr]{eucal}

\usepackage{epsfig}    % ??????
\usepackage{bbm}       % ??????

\usepackage[latin1]{inputenc}

\usepackage{epic}      %  for graphics
\usepackage{eepic}     % for graphics
\usepackage{calc}      % for computations of parameter settings

\usepackage{xy}        % for small diagrams, e.g. arrows
\xyoption{all}

\cleardoublepage % sorgt für korrekten druck bei doppelseiten

\numberwithin{equation}{section}

\newcounter{counter}
\theoremstyle{plain}
\numberwithin{equation}{subsection}

\newtheorem{thm}[equation]{Theorem}
\newtheorem{lem}[equation]{Lemma}

\newtheorem{prop}[equation]{Proposition}
\newtheorem{cor}[equation]{Corollary}

\theoremstyle{definition}

\newtheorem{defn}[equation]{Definition}
\newtheorem{mumford}[equation]{The Mumford Conjecture}

\newtheorem*{claim}{Claim}
\newtheorem{quests}[equation]{Questions}
\newtheorem{rem}[equation]{Remark}
\newtheorem{example}[equation]{Example}

\newtheorem{conv}[equation]{Convention}

\newtheorem{assumption}[equation]{Assumption}

\newtheorem{quest}[equation]{Question}

\newtheorem{expl}[equation]{Example}
\newtheorem{expls}[equation]{Examples}

\theoremstyle{remark}

\newenvironment{nam}[2]{\medskip {\bf #1 #2.}\it }{\medskip} % ZUM ZITIEREN INTERNER THEOREME

\theoremstyle{plain}
\theoremstyle{definition}

\newcommand{\bC}{\mathbb{C}}
\newcommand{\bD}{\mathbb{D}}

\newcommand{\bG}{\mathbb{G}}
\newcommand{\bH}{\mathbb{H}}

\newcommand{\bN}{\mathbb{N}}

\newcommand{\bP}{\mathbb{P}}
\newcommand{\bQ}{\mathbb{Q}}
\newcommand{\bR}{\mathbb{R}}
\newcommand{\bS}{\mathbb{S}}

\newcommand{\bZ}{\mathbb{Z}}

\newcommand{\bTh}{\mathbb{T}\bold{h}}

\newcommand{\gE}{\bold{E}}
\newcommand{\gF}{\bold{F}}

\newcommand{\gM}{\bold{M}}

\newcommand{\gU}{\bold{U}}

\newcommand{\gZ}{\bold{Z}}

\newcommand{\cA}{\mathcal{A}}
\newcommand{\cB}{\mathcal{B}}
\newcommand{\cC}{\mathcal{C}}

\newcommand{\cF}{\mathcal{F}}

\newcommand{\cO}{\mathcal{O}}

\newcommand{\cS}{\mathcal{S}}
\newcommand{\cT}{\mathcal{T}}

\newcommand{\cV}{\mathcal{V}}

\newcommand{\fM}{\mathfrak{M}}

\newcommand{\actson}{\curvearrowright}

\newcommand{\arf}{\operatorname{arf}}
\newcommand{\Arf}{\operatorname{Arf}}

\newcommand{\cpeins}{\mathbb{C} \mathbb{P}^{1}}
\newcommand{\cpinf}{\mathbb{C} \mathbb{P}^{\infty}}
\newcommand{\eins}{\mathbbm{1}}

\newcommand{\End}{\operatorname{End}}

\newcommand{\Emb}{\operatorname{Emb}}

\newcommand{\Hom}{\operatorname{Hom}}

\newcommand{\num}{\operatorname{num}}
\newcommand{\den}{\operatorname{den}}

\newcommand{\Gr}{\operatorname{Gr}}
\newcommand{\Bun}{\operatorname{Bun}}

\newcommand{\Pic}{\operatorname{Pic}}
\newcommand{\At}{\operatorname{At}}

\newcommand{\Sl}{\operatorname{SL}}

\newcommand{\Sspin}{\mathscr{SPIN}}

\newcommand{\inc}{\operatorname{inc}}

\newcommand{\prt}{\operatorname{prt}}
\newcommand{\Sym}{\operatorname{Sym}}

\newcommand{\vol}{\operatorname{vol}}
\newcommand{\projection}{\operatorname{pr}}

\newcommand{\Spiff}{\operatorname{Sdiff}}
\newcommand{\Sdiff}{\operatorname{Sdiff}}
\newcommand{\Gl}{\operatorname{GL}}
\newcommand{\spingroup}{\widetilde{\Gl_{d}^{+}(\mathbb{R})}}

\newcommand{\Aut}{\operatorname{Aut}}
\newcommand{\Diff}{\operatorname{Diff}}
\newcommand{\Spin}{\operatorname{Spin}}

\newcommand{\mb}[1]{\mathbb{#1}}

\newcommand{\hofib}{\operatorname{hofib}}

\newcommand{\ch}{\operatorname{ch}}
\newcommand{\colim}{\operatorname{colim}}
\newcommand{\ind}{\operatorname{ind}}

\newcommand{\td}{\operatorname{td}}
\newcommand{\tr}{\operatorname{tr}}
\newcommand{\trf}{\operatorname{trf}}
\newcommand{\id}{\operatorname{id}}
\newcommand{\delbar}{\bar{\partial}}
\newcommand{\Th}{\operatorname{Th}}

\newcommand{\bF}{\mathbb{F}}

\newcommand{\Sp}{\operatorname{SP}}
\newcommand{\SO}{\operatorname{SO}}

\newcommand{\Ob}{\mathfrak{O}\mathfrak{b}}
\newcommand{\Mor}{\mathfrak{M}\mathfrak{o}\mathfrak{r}}
\newcommand{\pt}{\ast}
\newcommand{\Mat}{\operatorname{Mat}}

\newcommand{\coker}{\operatorname{coker}}

\begin{document}

\title{Characteristic classes of spin surface bundles: Applications of the Madsen-Weiss theory}

\author{\textbf{Dissertation}\\
\\
zur\\
\\
Erlangung des Doktorgrades (Dr. rer. nat.)\\
\\
der\\
\\
Mathematisch-Naturwissenschaftlichen Fakult\"at\\
\\
der\\
\\
Rheinischen Friedrich-Wilhelms-Universit\"at Bonn\\
\\
\\
\\
\\
vorgelegt von\\
\\
Johannes Felix Ebert\\
\\
aus Bonn\\
\\
e-mail: ebert@math.uni-bonn.de\\
\\
\\
\\
Bonn, Juli 2006\\
}

\date{}

\maketitle

\thispagestyle{empty}

\pagebreak

Angefertigt mit Genehmigung der Mathematisch-Naturwissenschaftlichen Fakult\"at der Rhei\-nischen Friedrich-Wilhelms-Universit\"at Bonn\\
\\
\\
\\
\\
\\
\\
\\
\\
\\
\\
\\
\\
\\
\\
\\
\\
\\
\\
\\
\\
\\
\\
\\
\\
\\
\\
\\
\\
\\
\\
\\
1. Referent: Prof. Dr. Carl-Friedrich B\"odigheimer (Bonn)\\
\\
2. Referent: Prof. Dr. Fritz Grunewald (D\"usseldorf)\\
\\
3. Referentin: Prof. Dr. Ulrike Tillmann (Oxford)\\
\\
\\
Tag der Promotion: 15.11.2006

\thispagestyle{empty}

\pagebreak

\begin{abstract}
In this work, we study topological properties of surface bundles, with an emphasis on surface bundles with a spin structure.\\
We develop a criterion to decide whether a given manifold bundle has a spin structure and specialize it to surface bundles.\\
We study examples of surface bundles, in particular sphere and torus bundles and surface bundles induced by actions of finite groups on Riemann surfaces. The examples are used to show that the obstruction cohomology class against the existence of a spin structure is nonzero.\\
We develop a connection between the Atiyah-Singer index theorem for families of elliptic operators and the modern homotopy theory of moduli spaces of Riemann surfaces due to Tillmann, Madsen and Weiss.\\
This theory and the index theorem is applied to prove that the tautological classes of spin surface bundles satisfy certain divisibility relations. The result is that the divisibility improves, compared with the non-spin-case, by a certain power of $2$. The explicit computations for sphere bundles are used in the proof of the divisibility result.\\
In the last chapter, we use actions of certain finite groups to construct explicit torsion elements in the homotopy groups of the mapping class and compute their order, which relies on methods from algebraic $K$-theory and on the Madsen-Weiss theorem.
\end{abstract}

\thispagestyle{empty}

\pagebreak

\tableofcontents

\pagebreak

\section{Introduction}\label{intro}

In the last decades, the theory of Riemann surfaces and their moduli spaces has seen a rapid development. In particular, there is now an elaborate \emph{homotopy theory} of moduli spaces of Riemann surfaces.\\ 
This homotopy theory of moduli spaces is based on Teichm\"uller theory. Teichm\"uller theory tells us that the moduli space $\fM_g$ of conformal equivalence classes of Riemann surfaces of genus $g$ is the quotient of a space $\cT_g$ homeomorphic to $\bC^{3g-3}$ by the action of the \emph{mapping class group} $\Gamma_g$. This is the group of connected components of $\Diff(M_g)$, the group of orientation-preserving diffeomorphisms of a surface of genus $g$. The action is properly discontinuous and has finite isotropy groups. As a consequence, $\fM_g$ and the classifying space $B \Gamma_g$ are rationally homology equivalent.\\
But the connection between the object $\fM_g$ from complex analysis and the diffeomorphism group $\Diff(M_g)$ from differential topology is much stronger. It comes from a very special property of $2$-dimensional manifolds. Any oriented surface $M$ has a complex structure, but there is a more global statement. The \emph{space} $\cS(M)$ of all complex structures on $M$ is actually \emph{contractible}. As a consequence, Earle and Eells proved in 1969 that the unit component of the diffeomorphism group is contractible. Thus, $\Diff(M)$ and $\Gamma_g$ are the same thing for the eyes of a topologist. Another, simpler, consequence of the contractibility of $\cS(M)$ is the following. The topological group $\Diff(M)$ appears as the structural group of oriented smooth \emph{surface bundles}; and the contractibility of $\cS(M)$ shows that on any surface bundle $E \to B$, we can find complex structures on the fibers which vary continuously. Moreover, all these complex structures are \emph{concordant}.\\
The next important and much deeper step towards the homotopical theory of moduli spaces was Harer's stability theorem. For this, one needs to introduce mapping class groups $\Gamma_{g,n}$ for surfaces of genus $g$ with $n$ boundary components. There are stabilization maps between these mapping class groups given by boundary connected sum (for more details, see \ref{madsenweissthms}). Then the homology groups $H_k (B \Gamma_{g,n})$ and thus $H_k(B \Diff(M_g))$ do not depend on $g$ and $n$, as long as $g$ is large enough compared to $k$.\\
Based on Harer stability, Ulrike Tillmann proved a very strong theorem. The mapping class groups $\Gamma_g$ and $\Gamma_{g,n}$ are perfect; and one can form the infinite mapping class group $\Gamma_{\infty,n}$, as the colimit over all stabilization maps. One can apply Quillen's plus construction to $B \Gamma_{\infty,n}$. The result is a simply connected space $B \Gamma_{\infty,n}^{+}$ and a homology equivalence $B \Gamma_{\infty,n} \to B \Gamma_{\infty,n}^{+}$. Because of Harer stability, the homotopy type of $B \Gamma_{\infty,n}^{+}$ does not depend on $n$ and will be abbreviated by $B \Gamma_{\infty}^{+}$.\\
Tillmann's theorem states that $\bZ \times B \Gamma_{\infty}^{+}$ is an infinite loop space. In other words, there exists a connective spectrum and whose infinite loop space is homotopy equivalent to $ \bZ \times B  \Gamma_{\infty}^{+}$.\\
Later, Madsen and Weiss identified this spectrum. There exists a certain Thom spectrum $\bG_{-2}^{\SO}$ and a homotopy equivalence $\alpha:\bZ \times B \Gamma_{\infty}^{+} \to \Omega^{\infty} \bG_{-2}^{\SO}$, which is a Pontryagin-Thom type construction.\\
The most important ingredient for the proof of the Madsen-Weiss Theorem is Harer's stability theorem.\\
Harer, Bauer and Galatius have proven analogous results for surfaces with \emph{spin structures}, which are one of the main topics of the present work. Now I turn to the results of this dissertation.\\

\paragraph*{Chapter \ref{generalities}:}

Chapter \ref{intro} is this introduction and chapter \ref{generalities} provides the necessary background from the theory of moduli spaces and mapping class groups, except for the modern homotopy theory of moduli spaces, which is described in chapter \ref{thsthomoth}. I describe in detail the connection between the differential topology of surface bundles and the moduli spaces, which are objects of complex geometry (\ref{homotopyanalysis}, \ref{sbabors}, \ref{algeom}). I give a proof of the folklore statement that the classifying space $B \Gamma_g$ of the mapping class group and the moduli space $\fM_g$ of Riemann surfaces are rationally homology equivalent (Theorem \ref{mue}).\\
In section \ref{chclosb}, I give the definition of the most important characteristic cohomology classes of surface bundles. There are the \emph{Morita-Miller-Mumford classes} (or \emph{MMM-classes} for shortness) $\kappa_n \in H^{2n}(B \Diff(M); \bZ)$ (see Definition \ref{defnmum}) and the \emph{symplectic classes} $\gamma_n \in H^{2n}(B \Diff(M); \bZ)$ (see Definition \ref{defnsymp}) for any oriented closed surface $M$.\\
In section \ref{Hodge}, I describe an important series of vector bundles on the moduli spaces of Riemann surfaces, or rather on $B \Diff(M)$, the so-called \emph{Hodge bundles} $V_n$, $ n \in \bN$. They can be viewed as the vector bundles of holomorphic $n$-differentials when one chooses a complex structure on the fibers of the universal surface bundle.\\

\paragraph*{Chapter \ref{spinstrmanibun}:}

We carry out the preparations for the discussion of spin structures on surface bundles. We study more generally spin structures on arbitrary manifold bundles.
After a brief discussion of spin structures on vector bundles (section \ref{spinvector}) and manifolds (section \ref{spinmfd}), we discuss two important questions: 

\begin{quests}
\begin{enumerate}
\item What should the \emph{space} of all spin structures on $M$ be?
\item What is the "correct" space of automorphisms of a manifold with spin structure?
\end{enumerate}
\end{quests}

We address both questions in section \ref{spinspace}, and we do it in a uniform way. Namely, we construct both spaces as classifying spaces of a groupoid. The first groupoid is $\Sspin(M)$ and its classifying space is the space of all spin structures.
The second groupoid is $\Spiff(M)$, which has the same objects as $\Sspin(M)$, but there are much more morphisms, namely, all \emph{spin diffeomorphisms} of $M$.
One should think that $B \Spiff(M)$ is the space of all spin mani\-folds diffeomorphic to $M$ and that $B \Diff(M)$ is the space of all manifolds diffeomorphic to $M$. Thus one expects a fiber sequence

\begin{nam}{}{\ref{rewq}}
$$B \Sspin(M) \to B \Spiff(M) \to B \Diff(M);$$
\end{nam}

which is indeed true. In order to answer the question whether a manifold bundle classified by a map $B \to B \Diff(M)$ admits a spin structure, the fiber sequence \ref{rewq} is not suitable. One cannot apply obstruction theory in a reasonable way, because the situation is too non-Abelian. We need an auxiliary fiber sequence, which is related to the first one by the diagram

\begin{nam}{}{\ref{rewq}}
\xymatrix{ 
\bR \bP^{\infty} \ar[r] \ar[d]& B \Sspin(M) \ar[d]\\
B \Spiff(M, \sigma) \ar[d] \ar[r] & B \Spiff(M) \ar[d]\\
B \Diff(M, \sigma) \ar[r] & B \Diff(M),\\
}
\end{nam}

where $\Diff(M, \sigma)$ is the group of all diffeomorphisms of $M$ which fix the given spin structure $\sigma$. We stress that this group and its classifying space play a purely auxiliary role.
The bottom horizontal map is a finite covering, and the question whether a lift exists is group-theoretical and cannot be addressed by cohomological methods. But the left-hand side fibration is \emph{simple} and the obstruction theory behaves nicely. The fibration is classified by a cohomology class $c  \in h^2 (B \Diff(M, \sigma); \bF_2)$ and we will prove:

\begin{nam}{Corollary}{\ref{obstructionclasas}}
Let $E \to B$ be a smooth $M$-bundle, classified by a map $\lambda:B \to B \Diff(M)$ such that its monodromy $\rho= \pi_1 (B \lambda): \pi_1 (B) \to \pi_0 (\Diff(M))$ fixes $\sigma$. Then the following statement holds.\\
There exists a spin structure on the fiber bundle $E$ extending $\phi^{*} \sigma$ if and only if the mono\-dromy homomorphism takes values in $\pi_0 (\Diff(M, \sigma)$ and $\lambda^{*}c=0 \in H^2 (B; \bF_2)$.
\end{nam}

The two obvious questions which arise from this statement are:

\begin{quests}
\begin{enumerate}
\item
Is the cohomology class $c$ nontrivial, i.e. is the group extension $\bZ/2 \to \Spiff(M, \sigma) \to \Diff(M,\sigma)$ nontrivial?
\item Can one express the cohomology class $c$ in terms of more familiar classes?
\end{enumerate}
\end{quests}

\paragraph*{Chapter \ref{spinstrsurbun}:}

The whole theory in chapter \ref{spinstrmanibun} is valid for manifolds of arbitrary dimension, but in chapter \ref{spinstrsurbun}, we turn to surfaces.
I will review the results of Atiyah and Johnson on spin structures on surfaces in section \ref{spinsurfaces}. The most important thing about spin structures on surfaces is the existence of the \emph{Atiyah-invariant} (see Definition \ref{atiyahinv}) of a spin structure. It takes values in $\bF_2$ and divides the spin structures into the even ones and the odd ones. There are two definitions, one using index theory and the other one using elementary differential topology and linear algebra over $\bF_2$.
We give an easy new proof of Atiyah's theorem that any diffeomorphism of a surface fixes a spin structure (Corollary \ref{Stiefel}).
In the section \ref{spinmcg}, we apply the general theory of chapter \ref{spinstrmanibun} to surfaces. The  main result of this short section is 

\begin{nam}{Proposition}{\ref{haupt}}
Let $M$ be an oriented surface of genus $g \geq 2$.
\begin{enumerate}
\item $B \Spiff(M)$ has two connected components $B \Spiff(M)^{+}$ and $B \Spiff(M)^{-}$ belonging to the different values of the Atiyah invariant.
\item Both components are aspherical if $g \geq2$.
\item The fundamental group $\hat{\Gamma}_{g}^{\epsilon}:=\pi_1 (B \Spiff(M)^{\epsilon})$ is a central $\bZ/2$-extension of $\Gamma(M, \sigma)$, where $\sigma$ is a spin structure with Atiyah invariant $\epsilon \in \{ \pm 1\}$.
\end{enumerate}
\end{nam}

The groups $\hat{\Gamma}_{g}^{\epsilon}$ are the \emph{spin mapping class groups}. We do not study spin mapping class groups of surfaces with boundary. They are studied for example in \cite{bau} and are extremely important for the development of the homotopy theory of spin surface bundles. On the other hand, we shall only \emph{use} the results and for this use, we do not need the mapping class groups with boundary explicitly.\\

\paragraph*{Chapter \ref{someexamples}:}
 
The purpose of this chapter is twofold. I will introduce another theme of this work, the application of \emph{finite group actions} on surfaces to the topological theory of moduli spaces. Finite group actions on surfaces were employed by several authors, who achieved some progress: Glover and Mislin detected torsion classes in the cohomology of $B \Gamma_g$; Kawazumi, Akita, Uemura gave a new and elementary proof of the algebraic independence of the Morita-Mumford classes. Even more recently, Galatius, Madsen and Tillmann used finite group actions to prove that some divisibility properties of the Morita-Mumford classes are optimal (see below).\\
If $G$ is a finite group $M$ a Riemann surface $G \actson M$ a group action, we can study the surface bundle $E(G;M):= EG \times_G M \to BG$ and can ask whether there exists a spin structure on $E(G;M)$. The answer is simple and general:

\begin{nam}{Theorem}{\ref{exisspin}}
Let $G$ be a finite group which acts \emph{faithfully} on a closed surface $M$. Then the induced surface bundle $E(G;M)$ is spin if and only if all $2$-Sylow-subgroups act freely on $M$.
\end{nam}

We use Theorem \ref{exisspin} to show that $c$ is nonzero.\\
\\
The remaining sections \ref{spherebundles} and \ref{torusbundles} of chapter \ref{spinexpls} are devoted to the study of the exceptional genera $g=0$ and $g=1$. It is well-known that any $\bS^2$-bundle is the sphere bundle of a $3$-dimensional oriented vector bundle, which is unique up to isomorphism. In other words, $\Diff(\bS^2) \simeq \SO(3)$. We study the extension \ref{rewq} for $M = \bS^2$:
 
\begin{nam}{Proposition}{\ref{esszwei}}
In the diagram below, all horizontal maps are homotopy equivalences.

$$\xymatrix{
\bZ/2 \ar@{=}[r] \ar[d]& \bZ/2 \ar@{=}[r] \ar[d] & \bZ/2 \ar[d] \\
SU(2) \ar[r] \ar[d]  & \Sl_2 (\bC) \ar[r]\ar[d] & \Spiff(\bS^2, \sigma_0)\ar[d] \\
\SO(3) \ar[r]  & \bP \Sl_2 (\bC) \ar[r]   & \Diff(\bS^2), \\
}$$

\end{nam}

where $\sigma_0$ is the unique spin structure on $\bS^2$. Using this, we can express the spin condition for $\bS^2$-bundles in more familiar terms.

\begin{nam}{Proposition}{\ref{eszwei}}
Let $\pi: E \to B$ be an $\bS^2$-bundle and let $V \to B$ be a $3$-dimensional vector bundle whose sphere bundle is $E$. Then the following conditions are equivalent.
\begin{enumerate}
\item $E$ has a spin structure.
\item $V$ is a spin vector bundle.
\item There exists a $2$-dimensional complex vector bundle $U \to B$, such that $E$ is isomorphic to the projective bundle $\bP U \to B$ and such that $c_1 (U) =0$.
\end{enumerate}
Because of 1 and 2, the class $c \in H^2(B \Diff(\bS^2), \bF_2) \cong H^2(B\SO(3); \bF_2) \cong \bF_2$ agrees with $w_2$.
\end{nam}

Then we calculate the characteristic classes of $\bS^2$-bundles in terms of the characteristic classes of related vector bundles (Propositions \ref{classesgnull}, \ref{projcompl} and \ref{hpunendlich}).\\
We also do calculations for torus bundles. Besides the calculation of the characteristic classes for torus bundles, the main result of section \ref{torusbundles} is

\begin{nam}{Corollary}{\ref{extorspin}}
A torus bundle has an odd spin structure if and only if $\gamma_1 =0 \pmod {2}$. \qed
\end{nam}

\paragraph*{Chapter \ref{thsthomoth}:}
 
This chapter introduces the modern homotopy theory of moduli spaces after Tillmann, Madsen and Weiss. This theory is centered around the Madsen-Tillmann spectrum $\bG_{-2}^{\SO}$ and a map $\alpha$ from the classifying space of the stable mapping class group $\bZ \times B\Gamma_{\infty}$ to the infinite loop space $\Omega^{\infty} \bG_{-2}^{\SO}$. The map $\alpha$ is a homology equivalence by the Madsen-Weiss theorem.\\
Our focus is on a close connection of the map $\alpha$ with the index theory of families of Cauchy-Riemann operators.\\
In the first four sections of chapter \ref{thsthomoth}, we develop the general theory. Section \ref{beckertrans} explains the Becker-Gottlieb transfer and section \ref{pushforward} the relation with pushforward (alias umkehr maps) in cohomology. Both things are classical and the only purpose of these sections is to have the definitions and important properties at hand. Section \ref{unibecgott} introduces the Madsen-Tillmann spectra $\bG_{-d}^{\gF}$ in a general context. In section \ref{relbordism}, we explain the connection between the Madsen-Tillmann spectra and the classical Thom spectra, which is employed for many computations. Section \ref{madsenweissthms} is a short survey on the results of Harer, Madsen, Tillmann and Weiss and a few other people. We also consider the case of spin surfaces. In section \ref{miscellaneous}, we give a few computations which are used at several places in this dissertation.\\
In section \ref{atiyahsinger}, we show how the Atiyah-Singer index theorem for families of elliptic diffe\-rential operators fits into this context. In fact, we concentrate on Dolbeault operators on bundles of Riemann surfaces. Section \ref{indexamples} shows a few examples and in the final section \ref{mumfordtomadsen}, we show how the Madsen-Tillmann map naturally arises from Mumfords conjecture and index theory.

\paragraph{Chapter \ref{divisibility}}

It is a classical topic in the theory of characteristic classes that the presence of spin structures improves the divisibility of characteristic numbers. Let us discuss a simple example. Consider an oriented closed $4$-manifold $M$. By Hirzebruch's signature theorem, the signature of $M$ equals $\frac{1}{3} \langle p_1 (TM), [M]\rangle$. Since the signature is an integer, it follows that $p_1 (TM)$ is divisible by $3$. On the other hand, the signature of $\bC \bP^2 $ is $1$ and the divisibility is optimal for oriented $4$-manifolds. Rochlin's theorem (\cite{LM}, p. 288) gives an improvement of this for spin manifolds. More specifically, if $M$ has a spin structure, then the signature is divisible by $16$ or the first Pontryagin class is divisible by $48$. A proof can be given using the Atiyah-Singer theorem. If $M$ is spin, then there exists the Dirac operator $\not{\!\! D}   $, an elliptic differential operator whose index is $\frac{1}{24}\langle p_1 (TM); [M]\rangle$. Thus the signature is divisible by $8$. The last power of $2$ is encoded in the algebraic symmetries of the Dirac operator. Namely, there exists a parallel quaternionic structure on the spinor bundles, whence the index is an even number (compare \cite{LM}, p. 288). This divisibility of the signature of $4$-dimensional spin manifolds is optimal: a $K3$-surface has signature $16$.\\
In the chapter \ref{divisibility}, we prove divisibility theorems for the characteristic classes $\kappa_n$ of surface bundles of a similar spirit. 

\begin{nam}{Proposition}{\ref{obvious}}
For a spin surface bundle, $\kappa_n$ is divisible by $2^{n+1}$. This holds even integrally.
\end{nam}

For the rest of the discussion, we have to distinguish between even and odd values of $n$. Let us first consider the even classes $\kappa_{2n}$. 

\begin{nam}{Theorem}{\ref{geradteilbar}}
For spin surface bundles, the class $\kappa_{2n}$ is not divisible by any nontrivial multiple of $2^{2n+1}$. This holds in the stable range for spin mapping class groups.
\end{nam}

The condition on the stable range has the following reason. If the genus is not large enough, then all Morita-Mumford classes vanish, and they are divisible by any number. Thus the genus must be large enough to make sense out of this statement.\\
Our proof consists of two steps. First, we use the calculations of section \ref{spherebundles} to prove that $\kappa_{2n}$ of the universal sphere bundle with spin structure is not divisible by $2^{2n+2}$. This is not in the stable range, but can be used to prove results in the stable range. The second step uses the Madsen-Weiss theorem to reduce the statement about the cohomology of $\Omega^{\infty} \bG_{-2}^{\Spin}$. Theorem \ref{geradteilbar} means that the reduction of $\kappa_{2n}$ modulo $2^{2n+2}$ is nonzero and this can be proven using the classifying map $\bH \bP^{\infty} \to \Omega^{\infty} \bG_{-2}^{\Spin}$ of the universal $\bS^2$-bundle with spin structure.\\
For the odd Morita-Mumford classes, there is an old divisibility theorem by Morita (\cite{Mor1}). It says that the denominator of the Bernoulli numbers $\den(\frac{{B_{n}}}{2n})$ divides $\kappa_{2n-1}$ for any oriented surface bundle. Recently, Galatius, Madsen and Tillmann (\cite{GMT}) proved that this is optimal.
Here we prove

\begin{nam}{Theorem}{\ref{Clifford}}
For spin surface bundles, the class $\kappa_{2n-1}$ is divisible by $2^{2n} \den (\frac{B_{n}}{2n}) $.
\end{nam}

The proof of Theorem \ref{Clifford} follows the pattern of the proof of Rochlin's theorem.
The existence of a spin structure implies the existence of a holomorphic square-root $S$ of the vertical cotangent bundle and a differential operator $\delbar_S$. The Grothendieck-Riemann-Roch formula is used to prove that $2^{2n-1} \den (\frac{B_{n}}{2n}) $ divides $\kappa_{2n-1}$, and the last power of $2$ is grasped by a look on internal symmetries of the index bundle of $\delbar_S$.
We do not know whether this divisibility relation is optimal. The construction of surface bundles is difficult and the proof in \cite{GMT} does not generalize, because it relies essentially on the consideration of cyclic group actions on surfaces. The problem is a $2$-primary problem and due to Theorem \ref{exisspin}, there are not enough spin surface bundles given by finite group actions to prove the optimality.

\paragraph{Chapter \ref{icosahedral}} This last chapter is not immediately related to the other chapters of this dissertation in the sense that spin structures do not play any role there. The topic is the following. It is known that $\pi_3 (B \Gamma_{\infty}^{+}) \cong \bZ/24$ and we look for geometric representatives of elements in this group. Usually it is difficult to give a geometric description of maps from a space $Y$ into the plus construction of a space $X$, in particular, if $Y$ is simply-connected (like $\bS^3$) and $X$ is aspherical (like $B \Gamma_{\infty}$), then any map $Y \to X$ is nullhomotopic. There is a better chance for constructing a nontrivial map $\bS^n \to B \Gamma_{\infty}^{+}$ if one starts with a map $X \to B \Gamma_g$ when $X$ is an $n$-dimensional \emph{homology sphere}. Plus construction gives a map $X^+ \to B  \Gamma_{\infty}^{+}$ and $X^+ \simeq \bS^n$.
How do we construct a map $X \to B \Gamma_{\infty}$? It is given by a homomorphism $\pi_1 X \to \Gamma_{\infty}$ which in turn is given by any action $\pi_1 (X) \actson M_g$ if the fundamental group of $X$ is finite. The most popular $3$-dimensional homology sphere is the Poincar\'e sphere, whose fundamental group $\hat{G}$ is the binary icosahedral group. Using the connection between the fundamental group of the Poincar\'e sphere and Euclidean geometry, it is not difficult to construct nontrivial actions of $\hat{G}$ on Riemann surfaces (\ref{kleinsaction}).\\
How can we detect whether the element in $\pi_3 (B \Gamma_{\infty}^{+})$ is nontrivial? The method is the use of the Jones-Westbury formula (see \cite{JonesWestb}). We make use out of the map $B \rho^+: B \Gamma_{\infty}^{+} \to B \Gl_{\infty} (\bZ)^{+}$ into algebraic $K$-theory, which yields a map $\pi_3 (B \Gamma_{\infty}^{+}) \to K_3 (\bZ)$. There exists a homomorphism $e: K_3 (\bZ) \to \bQ / \bZ$, which is injective and maps onto $(\frac{1}{48} \bZ) / \bZ \subset \bQ / \bZ$. It turns out that $e \circ (B \rho^{+})_* : \pi_3 (B \Gamma_{\infty}^{+}) \to \bQ /\bZ$ is injective (subsection \ref{conclusion}). The element in $K_3 (\bZ)$ given by this procedure agrees with the element given by the representation of $\hat{G}$ on the first homology of the surface. The computation of the $e$-invariant of these elements is the subject of the paper \cite{JonesWestb}. To carry out the calculation, one needs information about the fixed points of the actions. In our case, this information is available and the computations yield that we have constructed an element in $\pi_3 (B \Gamma_{\infty}^{+})$ of order $12$ (see Proposition \ref{einvaria}).
Strictly speaking, we only have constructed an element in the stable homotopy of a mapping class group of a fixed genus which does not need to lie in the stable range, but this problem can be overcome by a stabilization (alias connected sum) operation (see section \ref{increasing}).\\
\\
I have included an appendix \ref{appendiy} which contains the proof of two results needed in section \ref{spinspace} and are abstract and technical in nature.\\

\paragraph*{Remark on the notations:}

This work is divided into chapters (numbered like chapter \ref{generalities}), sections (section 1 of chapter 1 is numbered as \ref{homotopyanalysis}). Equations, theorems etc occurring in section \ref{homotopyanalysis} are numbered like Theorem \ref{Earle}). I hope this makes my system of cross-references transparent.\\
Between the appendix \ref{appendiy} and the bibliography I have included a list of the notations used in this work. I tried to use a constant notation throughout the text. The advantage is that I could save a lot of phrases like "Let $\pi:E \to B$ be a surface bundle". The disadvantage of this method is that the reader may sometimes wonder what the meaning of a symbol is. To help in these situations, there is appendix \ref{listof}.\\
 
\pagebreak

\subsection*{Acknowledgements}

It is a great pleasure to thank everyone who helped me during my years as a PhD student.
There is one institution and one person who deserve the warmest thanks. The institution is the Max-Planck-Institut f\"ur Mathematik in Bonn which employed me for three years. Without this financial support, this work would not have been possible. During the last weeks of writing this thesis, I was financially supported by the Graduiertenkolleg 1150 "Homotopie und Kohomologie"\\
I am grateful my advisor Carl-Friedrich B\"odigheimer for all he has done for me in the last years. Almost every week he listened patiently to me when I was developping the ideas of this work and always answered my questions. He always encouraged me to push this work further.\\
I want to thank the Mathematisches Forschungsinstitut Oberwolfach for inviting me three times to very pleasant visits. The most important one for me was the first one: the Arbeitsgemeinschaft in October 2003 on the Madsen-Weiss theorem, organized by Søren Galatius and Michael Weiss. I take the opportunity to thank the organizers and speakers of this workshop, which was the beginning of my study of stable homotopy theory of moduli spaces.\\
I want to thank the American Institute of Mathematics in Palo Alto, California for the invitation to the workshop on the geometry and topology of moduli spaces in March 2005.\\
I had the opportunity to discuss with a lot of people about my research. I want to thank all of them, even if the results of these discussion did not enter the present work. Among them are: Anton Zorich, Maxim Kontsevich, Eduard Looijenga, Friedrich Hirzebruch, John Rognes, Karlheinz Knapp, Jeffrey Giansiracusa, Gerald Gaudens, Christian Ausoni, Holger Reich, Roland Friedrich and many more.\\
Furthermore, I want to thank the organizers and the speakers of the Oberseminar Topologie at the Mathematical Institute of the University of Bonn, which was run by Carl-Friedrich B\"odigheimer, Stefan Schwede, Gerald Gaudens, Christian Ausoni, Elke Markert and Birgit Richter. I could benefit from virtually every talk I heard there in the last years.\\
Furthermore, I want to mention my brother and sister students Balazs Visy, Maria Guadalu\-pe Castillo, Juan Wang, Ulrich Denker and Fridolin Roth for having such a nice time as a PhD student.
During the last weeks, Gerald Gaudens, Christian Ausoni and Elke Markert helped me with LaTeX problems.\\
\\
Bonn, July 2006

\pagebreak

\section{Generalities on surface bundles}\label{generalities}

This chapter is a survey on the theory of surface bundles and the mapping class groups.

\subsection{Complex structures on surfaces and the diffeomorphism
group}\label{homotopyanalysis}

Here, we clarify the relation between the classifying space of the mapping
class groups and the moduli space of Riemann surfaces.

Let $M$ be a smooth, closed, oriented and connected $2$-manifold. It is
well-known (\cite{Hirsch}) that $M$ is determined up to diffeomorphism by a
single invariant, the \emph{genus} of $M$;

$$g:= \frac{1}{2} \dim H_1(M; \bQ) \in \bN.$$

If we want to stress that the genus of $M$ is $g$, we also use the notation
$M_g$. Let $\Diff(M)$ be the group of \emph{orientation-preserving}
diffeomorphisms of $M$. The space of all smooth maps $M \to M$ carries the
$C^{\infty}$-Whitney topology (a sequence of maps converges if and only if all
its derivatives converge uniformly). The diffeomorphism group $\Diff(M)$ is an
open subspace. With the subspace topology, it becomes a topological group.\\
Let a complex structure on $M$ be given, let $\phi: M \supset U \to \bC$ be a holomorphic chart, $p \in U$ and let $T \phi$ be the derivative of
$\phi$. Then $v \mapsto T \phi^{-1} (i T \phi(v))$ does not depend on the choice of $\phi$ and defines a smooth bundle endomorphism $J: TM \to
TM$, such that $J^2 = - \id$ and such that for all nonzero tangent vectors $v \in T_p M$, $ p \in M$, the basis $(v, Jv)$ of $T_p M$ is
positively oriented.\\ A bundle endomorphism with these properties is called an \emph{almost complex structure} on $M$.\\ There is a converse to
this construction. Let $J$ be an almost complex structure on $M$ and let $p \in M$. For any $p \in M$, there exists a neighborhood $U \subset M$
of $p$ and a map $\phi: U \to \bC$, such that $\phi$ is a diffeomorphism $U \to \phi(U)$ and such that $T\phi \circ J = i T \phi$. The
collection of all these maps $\phi$ is a holomorphic atlas on $M$. Both constructions are mutually inverse. The statement that this holomorphic
atlas exists is called the \emph{Newlander-Nirenberg theorem}.\\ To prove the existence of $\phi$, one can use the existence of quasiconformal
mapping with given complex dilatation (see, for example, \cite{ImTan}). In higher dimensions, there is a similar statement, with an additional
integrability condition on $J$, but the proof is much more difficult.\\  Thus complex structures and almost complex structures on oriented
surfaces are equivalent notions. Also, there is a close connection with Riemannian metrics. We say that a Riemannian metric $h$ on a surface
with an almost-complex structure $J$ is \emph{conformal} if $J$ is orthogonal with respect to $h$. It is easy to see that conformal metrics
exist, that two differ by a multiplication with a positive function and that any multiple of a conformal metric with a positive function is
again conformal. Conversely, any Riemannian metric on an oriented surface determines an almost-complex structure: a Riemannian metric plus an
orientation determines the Hodge $*$-operator, which satisfies $* ^2 = - \id$ when acting on $1$-forms. Set $J:= - * ^{\prime}$ (the dual of the
Hodge star operator). This is a complex structure. Both constructions are mutually inverse.\\ We denote the set of all almost complex structures
by

$$\cS(M) \subset C^{\infty} (M, \End(TM)).$$

The vector space $C^{\infty}(M, \End(TM))$ of all smooth sections of
$\End(TM)$ carries the $C^{\infty}$-Whitney topology, which turns it into a
Fr\'echet space.\\ There is a slightly different description of $\cS (M)$. Let
$\Gl(M) \to M$ be the $\Gl_{2}^{+}(\bR)$-principal bundle of oriented frames
of $M$. It is easy to see that $\cS(M)$ can be identified with the space
of smooth sections in the fiber bundle $\Gl(M) \times_{\Gl_{2}^{+}(\bR)}
\Gl_{2}^{+}(\bR)/ \Gl_1(\bC)$. But $\Gl_{2}^{+}(\bR)/ \Gl_1(\bC)$ is
contractible (as a $\Gl_{2}^{+}(\bR)$-space, it is diffeomorphic to the upper
half plane). Thus $\cS(M)$ is contractible.\\ There is an obvious
right-action of $\Diff(M)$ on $\cS(M)$. Let $f \in \Diff(M)$ and $J \in
\cS(M)$. Set $J \cdot f := Tf^{-1} \circ J \circ Tf$. This action is
continuous. As a set, the quotient $\cS(M) / \Diff(M)$ is the same as the set
of conformal equivalence classes of Riemann surfaces of genus $g$. Let us
formulate this more precisely. Let $\Diff_0(M)$ be the subgroup of all
diffeomorphisms which are homotopic to the identity.

\begin{thm}(\cite{EaEe})\label{Earle} Let $g\geq2$.
\begin{enumerate}
\item The action of $\Diff(M)$ on $\cS(M)$ is proper; the action of
$\Diff_0(M)$ is free. The projection map $\cS(M) \to \cS(M) / \Diff_0(M)$
is a principal $\Diff_0(M)$-bundle.  \item There is a natural homeomorphism
$\cS(M) / \Diff_0(M) \to \cT_g$ (the classical Teichm\"uller space).
\end{enumerate} \end{thm}

A new, surprisingly simple proof of the first assertion in this theorem, as
well as a simp\-le proof of the next theorem can be found in the paper
\cite{RobSal}. Theorem \ref{Earle} has many consequences. Because $\cT_g
\cong \bR^{6g-6}$ by Teichm\"uller's theorem, in the fibration $\Diff_0(M) \to
\cS(M) \to \cT_g$ both, basis and total space are contractible. Therefore,
$\Diff_0(M)$ is also contractible. In particular, diffeomorphisms which are
homotopic are also isotopic. \\ For lower genera, Theorem \ref{Earle} needs to
be modified, because the action of $\Diff_0(M)$ is no longer free if $g=0,1$.
The result which is needed in this work is

\begin{thm}\label{Smale}(\cite{EaEe}, see also \cite{RobSal})
\begin{enumerate}
\item If $g=0$, then the inclusion $\SO(3) \to \Diff_0(\bS^2) =
\Diff(\bS^2)$ given by the rotations of $\bS^2$, is a homotopy equivalence
(This was already proved by Smale, \cite{Sma}). \item If $g=1$, then the
inclusion $T \to \Diff_0(T)$ given by the Lie group structure of the torus
$T$, is a homotopy equivalence. \end{enumerate} \end{thm}

The (discrete) group of components of $\Diff(M)$ is the \emph{mapping class group}

$$\Gamma_g := \pi_0 (\Diff(M)).$$

There is a proper action $\Gamma_g \actson
\cT_g$ and the quotient is the \emph{moduli space of Riemann surfaces} of
genus $g$; $\fM_g := \cT_g / \Gamma_g$. There is a complex structure on the
Teichm\"uller space such that the mapping class group acts by biholomorphic
maps. It follows that $\fM_g$ is a complex analytic space.\\  The action of
$\Gamma_g$ on $\cT_g$ is not free, because the stabilizer subgroup of a point
consists of all complex automorphisms of the Riemann surface represented by
this point. This group is always finite (if $g \geq2$), but in general not
trivial. We will say more about the relation between $\fM_g$ and $B \Gamma_g$
in section \ref{algeom}.

\subsection{Surface bundles and bundles of Riemann surfaces}\label{sbabors}

\begin{defn}
Let $M$ be a smooth oriented manifold. An oriented \emph{smooth $M$-bundle}
on a base space $B$ is a fiber bundle $\pi:E \to B$ with fiber $M$ and
structural group $\Diff^{+} (M)$.  \end{defn}

In this work, we will only consider oriented and smooth $M$-bundles with
\emph{compact fibers}. If not stated otherwise, any $M$-bundle is smooth and
oriented. Also, in most cases the dimension of $M$ will be $2$. Then we say
that $E$ is a \emph{surface bundle}. The \emph{genus} of a surface bundle is
the genus of the fiber - we do not consider bundles whose fibers change from
one component of the base space to another one.\\ If $\pi:E \to B$ is a smooth
$M$-bundle, then there is an associated $\Diff(M)$-principal bundle $Q \to B$.
Its fiber over $b \in B$ consists of all orientation-preserving diffeomorphisms
from $M$ to $\pi^{-1} (b)$. By construction, there is a natural isomorphism $Q
\times_{\Diff(M)} M \cong E$.

\begin{defn}
Under the present assumption, the \emph{vertical tangent bundle} of $E \to
B$ is the $d$-dimensional oriented vector bundle

$$T_v E:=Q \times_{\Diff(M)} TM \to Q \times_{\Diff(M)} M = E.$$

\end{defn}

The restriction of $T_v E$ to the fiber $E_b:=\pi^{-1}(b)$ of $E$ is
canonically isomorphic to the tangent bundle of $\pi^{-1}(b)$. It should be
remarked that if $E$ and $B$ are smooth manifolds and $\pi$ is a smooth
submersion, then the vertical tangent bundle is canonically isomorphic to
$\ker T \pi$ and there is a short exact sequence

$$0 \to T_v E \to TE \to \pi^{*} TB \to 0$$

of vector bundles on $E$. If the surface bundle $E$ is clear from the context, we sometimes write $T_v:= T_v E$.\\
If $E \to B$ is a surface bundle and $b \in B$ a point, then we write $E_b:=
\pi^{-1}(b)$ for the fiber.

\begin{defn}\label{defriembun}
Let $ \pi:E \to B$ be a surface bundle and $Q \to B$ the associated
$\Diff(M)$-principal bundle. Then a \emph{complex structure} $J$ on $E$ is a
continuous section of the fiber bundle

$$ Q \times_{\Diff(M)} \cS(M) \to B.$$

\end{defn}

Strictly speaking, we just defined almost-complex structures, and we know that any indivi\-dual fiber has a holomorphic atlas with respect to
$J_b:= J|_{E_b}$. But the Newlander-Nirenberg-Theorem can be improved to a parameterized version, which guarantees the existence of holomorphic
charts depending continuously on the base variable. More precisely, for any $e \in E$, we can find a neighborhood $V \subset E$ and a
homeomorphism $\phi:V \to \pi(V) \times \bD$ which commutes with the projections to $B$ and whose restriction to any fiber of $\pi$ is
biholomorphic. This is a consequence of the theory of quasiconformal mappings, in particular of the continuous dependence of the
solutions of the Beltrami-equation on the coefficients.\\
Because $\cS(M)$ is contractible, all complex structures on $E \to B$ are
essentially equivalent.

\begin{prop}\label{hascstruc}
Any surface bundle admits the structure of a bundle of Riemann surfaces.
Further, this structure is unique up to concordance, i.e. if $J_0,J_1$ are two
complex structures on $\pi$, then there exists a complex structure $J$ on the
surface bundle $\pi \times [0,1]$, such that $J|_{\pi \times \{i\}} = J_i$ for
$i=0,1$. \end{prop}

\textbf{Proof:} Because the fiber $\cS(M)$ of $Q \times_{\Diff(M)} \cS(M) \to B$ is
contractible, the space of all sections is contractible (in particular,
nonempty).\qed\\

For $g \geq2$, the classifying space $B \Gamma_g$ of the mapping class group
classifies isomorphism classes of surface bundles of genus $g$ and also
concordance classes of bundles of Riemann surfaces of genus $g$ - this follows
from the discussion above. The universal surface bundle $E(\Diff(M);M):= E
\Diff(M) \times_{\Diff(M)} M \to B \Diff(M)$ can be viewed as a surface bundle
on $B \Gamma_g$, which by abuse of language\footnote{The group $\Gamma_g$ does
\emph{not} act on the surface.} is denoted by $E(\Gamma_g, M_g)$. This bundle
has a complex structure, unique up to concordance.

\subsection{Characteristic classes of surface bundles}\label{chclosb}

To define the \emph{Morita-Miller-Mumford classes} of a surface bundle, we
need the notion of the cohomological push-forward or umkehr map. Let $\pi:E \to
B$ be an oriented surface bundle. Then there is a map $\pi_{!}:
H^{n} (E; \bZ) \to H^{n-2} (B; \bZ)$. We will discuss the definition and some
of the properties of $\pi_{!}$ in section \ref{pushforward}. If $E$ and $B$
are oriented manifolds, then there is an easy definition of $\pi_{!}$ in terms
of Poincar\'e duality (see \cite{Bred}).

\begin{defn}\label{defnmum}(\cite{Mor1}, \cite{Miller})
The \emph{Morita-Miller-Mumford classes} of $\pi:E \to B$ are by definition

$$\kappa_n(\pi):=\pi_{!}((e(T_v E))^{n+1} \in H^{2n}(B, \bZ).$$

\end{defn}

\begin{rem}
Sometimes, the Morita-Miller-Mumford classes are called only "Mumford classes" or "Morita-Mumford classes". We will prefer the short term \emph{MMM-classes}. Originally, Mumford (\cite{Mum}) defined rational cohomology classes in $H^{*}(\fM_g; \bQ)$, using techniques from algebraic geometry. We will say
more about these rational MMM-classes in section \ref{algeom}. It is not difficult to see that the Morita-Miller-Mumford classes are
natural in $B$. Therefore they define cohomology classes in $H^{2n}(B
\Diff(M); \bZ)$ for any closed Riemann surface $M$. The $0$-th
Morita-Miller-Mumford class is just the element $2-2g \in H^{0} (B; \bZ)$, by
the Gau{\ss}-Bonnet theorem.
\end{rem}

\begin{rem}
There is a slight problem in Definition \ref{defnmum}. We only can define the
MMM-classes for surface bundles over base spaces which are "not too large"
(see \ref{beckertrans} for a more precise condition). In order to know that
this defines a cohomology class in $H^{2n}(B \Diff(M); \bZ)$, we need to know
that $\Diff(M)$ is "not too large", which is not immediate, since it is
infinite dimensional. This is a $ \lim^1$-problem, which is implicitly solved
in \cite{GMT} and \cite{GMTW}. The way in which this problem is solved is not
important in the present work, see also section \ref{beckertrans}. \end{rem}

\begin{mumford}\label{mumconject}
In \cite{Mum}, Mumford conjectured that the homomorphism of graded algebras
$\bQ[\kappa_1, \kappa_2, \ldots ] \to H^{*} (B \Gamma_g, \bQ)$ is an
isomorphism in degrees $* < g/2$. This was quite recently proven by Madsen and
Weiss (\cite{MW}). We will say more about this in chapter \ref{thsthomoth}.\\
\end{mumford}

There is another series of characteristic classes. To define them,
we need the symplectic group $\Sp$. Let $R$ be a commutative ring with unit $1$, let $\eins_n$ be the unit matrix of size $n \times n$ and let

$$I:= \
\begin{pmatrix}
0 & - \eins_n \\
\eins_n & 0
\end{pmatrix}.
$$

The symplectic group is

\begin{equation}\label{sympgroup}
\Sp_{2n} (R):= \{ A \in \Mat (2n \times 2n; R) | A^t I A = I \}.
\end{equation}

It is a fact that $\Sp_{2n} (R) \subset \Sl_{2n} (R)$. Furthermore, $\Sp_{2n}(\bR)$ contains $U (n)$ as a maximal compact subgroup. In particular, $U(n) \subset \Sp_{2n} (\bR)$ is a homotopy equivalence.\\
The mapping class group $\Gamma_g$ acts on $H_1 (M_g; \bZ)\cong \bZ^{2g}$ and
that action preserves the intersection form on $H_1 (M; \bZ)$. If one chooses a symplectic basis of $H_1 (M; \bZ)$, one obtains a homomorphism

\begin{equation}\label{symplecto}
\rho: \Gamma_g \to \Sp_{2g}(\bZ),
\end{equation}

which can be composed with the inclusion $\Sp_{2g} (\bZ) \to \Sp_{2g}
(\bR)$. After passage to classifying spaces, we obtain

\begin{equation}
B \rho: B \Diff(M) \to B\Sp_{2g} (\bR) \simeq U(g).
\end{equation}

There are the integral Chern character classes $s_n \in H^{2n} (BU(g);
\bZ)$, which are as usual defined by $s_n := \sum_{i} x_i^{n}$, where $x_i$ is
the $i$-th Chern root.

\begin{defn}\label{defnsymp}
The \emph{symplectic classes} are defined by

$$\gamma_n := B \rho^{*} (s_n) \in H^{2n}(B \Diff(M); \bZ).$$

\end{defn}

There is another description, which is sometimes more useful. If the surface
bundle $\pi: E \to B$ is classified by a map $f: B \to B \Diff(M)$, then
$\gamma_n (\pi):= f^* \gamma_n$ is the integral Chern character class of the
symplectic vector bundle

$$H_{1}^{v} (\pi; \bR):= \bigcup_{b \in B} H_1 (\pi^{-1}(b); \bR) ,$$

which also can be described as the bundle induced by the maps $B \to B
\Diff(M) \to B \Sp_{2g}(\bR)$. The choice of a complex structure on $E$
determines a complex structure on the vector bundle $H_{1}^{v} (\pi; \bR)$ via
the Hodge decomposition: any element in the cohomology group $H^1
(\pi^{-1}(b); \bR)$ is represented by a harmonic $1$-form by the Hodge
Theorem, and the Hodge $\ast$-operator defines a complex structure on $H^1
(\pi^{-1}(b); \bR)$. Note that $H_{1}^{v} (\pi; \bR)$ is flat as a real vector
bundle and as a symplectic vector bundle, but \emph{not} as a complex vector
bundle.

\subsection{The Hodge bundles}\label{Hodge}

Let $\pi:E \to B$ be a bundle of Riemann surfaces. Then we can define a
series of complex vector bundles $V_n= V_n(\pi)$ on $B$, $n \in \bZ$,
associated to and natural in $\pi$.\\
Our definition starts with the notion of a holomorphic vector bundle on a
bundle of Riemann surfaces, which is not necessarily assumed to be a
holomorphic map. Before, we need to make precise the notion of a smooth vector
bundle on a general manifold bundle. Let $\pi:E \to B$ be a smooth $M$-bundle
and let $W \to E$ be a vector bundle. A \emph{smooth bundle atlas} is a bundle
atlas for $W$, such that all transition functions are smooth in
fiber-direction, with an additional condition: all fiber-direction-derivatives
of the transition functions are assumed to be continuous. An example of this
is the vertical tangent bundle, since we assume that the structural group is
$\Diff(M)$ with the Whitney topology. As usual, we say that two smooth bundle
atlases are equivalent if their union is a smooth atlas, and a smooth vector
bundle on $E$ is a vector bundle,together with an equivalence class of
atlases.\\

Now we can define the notion of a holomorphic vector bundle on a bundle of
Riemann surfaces.

\begin{defn}
Let $\pi: E \to B$ be a bundle of Riemann surfaces and $W \to E$ a smooth
vector bundle. Then a \emph{holomorphic atlas} on $W$ is a bundle atlas,
consisting of fiberwise smooth trivializations, such that all transition
function are fiberwise holomorphic\footnote{Then all derivatives of the
transition functions are automatically continuous, and a holomorphic bundle
atlas is in particular a smooth bundle atlas}. Two holomorphic atlases are
equivalent if their union is again a holomorphic bundle atlas. A \emph{holomorphic
vector bundle} on $E$ is a complex vector bundle, together with an equivalence
class of holomorphic atlases. \end{defn}

The prime example of a holomorphic vector bundle on a bundle of Riemann
surfaces is the vertical tangent bundle. The almost-complex structure $J$
turns $T_v E$ into a complex line bundle; and the complexification $T_v E
\otimes \bC$ decomposes into the $\pm i$-eigenspaces of $J$. The same is true
for the vertical cotangent bundle.

\begin{defn}
The \emph{vertical holomorphic cotangent bundle} $\Lambda_{v} E$ is the
eigenspace of $J$ to the eigenvalue $i$.
\end{defn}

The notion of a smooth section in a smooth vector bundle $W \to E$ on a
manifold bundle $E \to B$ also needs some explanation. Let $s: E \to W$ be a
continuous section. We say that $s$ is smooth if the restriction of $s$ to any
fiber $E_b$ is smooth and if all fiberwise derivatives are continuous on the
whole of $E$. We denote the space of all fiberwise smooth sections by the
symbol

$$C^{\infty}_{v}(E; W).$$

We should interpret the space $C^{\infty}_{v}(E; W)$ as the space of
continuous sections in a bundle $\cC^{\infty}( W)$ of Fr\'echet spaces on $B$.
More precisely, let $E_b := \pi^{-1}(b)$ and $W_b := W|_{E_b}$ and set

$$\cC^{\infty}( W):= \bigcup_{b \in B} C^{\infty} (E_b; W_b)$$

with the obvious projection map to $B$ and the obvious linear structure. We
say that a section $s: B \to \cC^{\infty}( W)$ is continuous if the adjoint
map $E \to W$, $e \mapsto s(\pi(e)) (e) $ is a smooth section of $W \to E$. The
result is a bundle of Fr\'echet spaces. The space of continuous sections in
$\cC^{\infty}( W)$ is isomorphic to the space $C^{\infty}_{v}(E; W)$.\\

If $W \to E$ is a fiberwise holomorphic vector bundle, then the
\emph{Cauchy-Riemann-operator}

$$\delbar_W: \cC^{\infty}( W) \to \cC^{\infty}( W \otimes \overline{\Lambda_v E})$$

is defined. It takes a point $s \in \cC^{\infty}( W)$ alias a smooth section
in $W_b \to E_b$ to $\delbar_{W_b} s $, which is a smooth section in $W_b
\otimes \Lambda_{E_b}$ over $E_b$. The topologies are designed to achieve that
$\delbar_W$ is a continuous bundle map.\\ If the fibers of $\pi$ are compact
and closed, then the regularity theory of elliptic differential operators
implies that Cauchy-Riemann operators are Fredholm operators, i.e. their
kernels and cokernels are finite-dimensional (see \cite{LM} for details).
Moreover, the \emph{Fredholm index} $b \mapsto \ind (\delbar_{W_b}) = \dim \ker
(\delbar_{W_b})- \dim \ker (\delbar_{W_b})$ is locally constant and the
\emph{index spaces}

$$\ind (\delbar_{W_b}) := \ker (\delbar_{W_b }) -  \coker (\delbar_{W_b })$$

fit together to form a virtual vector bundle $\ind_{\delbar_W} \in K^0 (B)$ on
$B$, at least if $B$ is compact. Actually, it is sufficient to assume that the
dimension of the kernel spaces are \emph{bounded}.

\begin{lem}\label{estimate}
If $\pi:E \to B$ is a bundle of Riemann surfaces and if $W \to E$ is a
holomorphic line bundle which has degree $l$ when restricted to any fiber of
$\pi$, then $\dim \ker (\delbar_{W_b})$ is bounded.
\end{lem}

\textbf{Proof:}
We need Clifford's theorem (\cite{GH}, p.251). It states that if $W_0, W_1$
are two holomorphic line bundles on a compact Riemann surface and if both have
nonzero holomorphic sections, then

$$\dim \ker \delbar_{W_0} + \dim \ker  \delbar_{W_1} \leq 1+ \dim \ker
\delbar_{W_0 \otimes W_1}.$$

If $W$ is a line bundle on a surface, then there exists a $q \in \bN$, which
only depends on $l$, such that $\Lambda^{\otimes q} W^{-1} $ has
nonvanishing holomorphic sections, and we apply Cliffords theorem to the pair $W,
\Lambda^{\otimes q} W^{\otimes(-1)}$ to obtain a bound on the dimension.\qed\\

We have seen that families of holomorphic vector bundles on the total space
give rise to a virtual vector bundle.\\
The Serre duality theorem (\cite{GH}, p. 102) has a parametrized version. Fix
a conformal metric on the bundle of Riemann surfaces (this is a contractible
choice) and fix a hermitian metric on the vector bundle $W$ (this is also a
contractible choice). Then we have a bilinear form of vector bundles

$$\cC^{\infty}( W \otimes \overline{\Lambda_v E}) \times \cC^{\infty}(
\overline{W} \otimes \Lambda_v E) \to B \times \bC,$$

which sends a pair of points $(s,t) \in  \cC^{\infty}( W \otimes
\overline{\Lambda_v E})_b \times \cC^{\infty}( \overline{W} \otimes \Lambda_v E)_b$
to

$$\int_{E_b} \beta (s,t) \vol,$$

where $\vol$ is the volume form on the Riemannian manifold $E_b$ and $\beta$
is the bilinear vector bundle map $C^{\infty}(E_b;  W_b \otimes \overline{\Lambda_{E_b}}) \times C^{\infty}(E_b; \overline{W_b} \otimes \Lambda_{E_b}) \to C^{\infty}(E_b; \bC) $. The Serre duality
theorem for a single Riemann surface asserts that this bilinear map
descends to a \emph{perfect pairing} of finite dimensional vector spaces

$$\coker (\delbar_{W_b}) \times \ker(\delbar_{\overline{W_b} \otimes \Lambda_{ E_b}}) \to \bC.$$

This statement extends to families of Riemann surfaces and we obtain

\begin{prop}\label{serredual}
There is a natural (up to contractible choices) perfect pairing

$$\coker (\delbar_{W}) \times \ker(\delbar_{\overline{W} \otimes \Lambda_v E}) \to B \times \bC$$

of finite dimensional families of vector spaces over $B$ and thus an
$\bC$-linear isomorphism

$$\coker (\delbar_{W}) \cong (\ker(\delbar_{\overline{W} \otimes \Lambda_v
E}))^{\prime}$$

of families of vector spaces (or of vector bundles, if the dimensions are
locally constant). \end{prop}

If we want to define a \emph{natural} vector bundle depending on a surface
bundle, then we need natural holomorphic vector bundles on Riemann surfaces.\\
Examples for these are the integral tensor powers of the vertical cotangent
bundle. These are essentially the only examples of natural holomorphic line
bundles. This uniqueness statement can be justified rigorously in the framework
of algebraic geometry, see \cite{HarMor}, p.62. These examples have an
additional feature, namely: the dimensions of the kernel and of the cokernel
only depend on the genus of the Riemann surface, not on the choice of a
complex structure. More precisely, if $M$ is a Riemann surface of genus $g$
and $\Lambda= \Lambda_M$ the canonical line bundle, then

\begin{equation}\label{rieroc}\dim \ker(\delbar_{ \Lambda_{M}^{\otimes m}}) =
\begin{cases}
1 ,&  m=0, g \in \bN,\\
g ,&   m=1, g \in \bN,\\
1-m,& g =0, m \leq 0,\\
0 ,  & g=0,  m >0,\\
1,& g=1, m \in \bZ,\\
0, & g \geq 2, m<0,\\
(2m-1)(g-1), & g \geq 2, m \geq2.\\
\end{cases}
\end{equation}

and $\dim \coker ( \delbar_{\Lambda_{M}^{ \otimes m}}) = \dim \ker(\delbar_{
\Lambda_{M}^{\otimes(1-m)}})$ by Serre duality. The proof can be
found in many books on compact Riemann surfaces, for example \cite{FaKr},
p. 80.\\

\begin{defn}\label{defnhodge}
Let $\pi: E \to B$ be a bundle of Riemann surfaces. Then the $n$-th
\emph{Hodge bundle} is the vector bundle

$$V_n (\pi) := \ker \delbar_{(\Lambda_v E)^{\otimes n  }}.$$

\end{defn}

Sometimes we also write $V_n (E)$ for $V_n(\pi)$. The bundle $H_{1}^{v} (\pi; \bR)$ defined in the preceding section can be viewed as a Hodge bundle. It is $V_{1}^{\prime} = \overline{V_1}$.\\
The Grothendieck-Riemann-Roch theorem computes the rational characteristic classes
of the index bundle $\ind \delbar_{V}$ when $V \to E$ is a family of
holomorphic vector bundles on $E$. The Grothendieck-Riemann-Roch theorem is
the cohomological version of the Atiyah-Singer index theorem for families of
elliptic operators (see \cite{AS4}), when applied to the Cauchy-Riemann
operators. The latter is much stronger and will be discussed in section
\ref{atiyahsinger}. Let $\td(x)=\frac{x}{1-e^{-x}} \in \bQ[[x]]$ be the Todd
power series. We now state the Grothendieck-Riemann-Roch theorem.

\begin{thm}(\cite{AS4})\label{grorieroc}
Let $  \pi:E \to B$ be a bundle of Riemann surfaces and let $V \to E$ be a
family of holomorphic vector bundles on $E$. Let $c \in H^2 (E)$ be the first
Chern class of the vertical tangent bundle. Then the Chern character of the
index bundle is given by the formula

$$ \ch(\ind \delbar_V) = \pi_{!} (\td(c) \ch(V))\in H^{*} (B; \bQ).$$

\end{thm}

However, it is important to note that we missed out an important point. The
Atiyah-Singer family index theorem (\cite{AS4}) is formulated only for bundles
$E \to B$ with a \emph{compact} base space $B$. There are several problems if
one tries to generalize it to arbitrary (at least, paracompact) spaces $B$.
One of the problems is resolved by Lemma \ref{estimate}, but there are other
difficulties. The author of this dissertation does not know whether this
difficulties can be overcome, they are not even easy to isolate. The only
possibility is to check the proof in \cite{AS4} step by step, which can be
expected to be cumbersome and is not done in this work.\\

\subsection{Comparison with algebraic geometry}\label{algeom}

We now construct a comparison map $\mu: \Diff(M) \to  \fM_g$ if $M$ is a
differentiable surface of genus $g$. For $g \geq 2$, the source space is
homotopy equivalent to $B \Gamma_g$. We begin by choosing a complex structure
on the universal surface bundle

$$\pi: E( \Diff(M);M) := E \Diff(M) \times_{\Diff(M)} M \to B \Diff(M),$$

which is possible by the discussion following Definition \ref{defriembun}. If $b \in B \Diff(M)$, then $\pi^{-1} (b)$ is a Riemann surface and we set $\mu: b \mapsto [\pi^{-1}(b)] \in \fM_g$.
This map turns out to be continuous. That is one of the reasons for the
continuity condition in the definition \ref{defriembun}. However, $\mu$ would
certainly be continuous with a weaker notion of continuity, as a consequence of
Teichm\"uller theory.

\begin{cor} to Prop \ref{hascstruc}.
The homotopy class of the moduli map

$$\mu: B \Diff(M) \to \fM_g$$

does not depend on the choice of the complex structure on $E(\Diff(M); M)$.
\end{cor}

The main theorem of this section is

\begin{thm}\label{mue}
If $g \geq 2$, then $\mu$ induces an isomorphism in (co)homology with rational coefficients.
\end{thm}

This theorem is a folklore statement. However, I was not able to find a
proof in the published literature. Thus I present a proof here.
For smaller genera, there is a similar statement. The mapping class group
$\Gamma_1$ of a torus is isomorphic to $\Sl_2(\bZ)$, whose rational homology
is trivial. On the other hand, the Teichm\"uller space of tori is the upper
half plane and the action of the mapping class group is given by M\"obius
transformations. The quotient $\bH / \Sl_2(\bZ)$ is biholomorphic equivalent
to the complex plane $\bC$. Thus the rational homology of the moduli space
$\fM_1$ is also trivial. Note that, however, the space $B \Sl_2(\bZ)$ is
\emph{not} the classifying space of torus bundles, but rather the space of
torus bundles with a section (Any two points on a complex torus can be
translated by a holomorphic automorphism). Morita showed (\cite{Mor4}) that
the space $B \Diff(T)$ is not rationally acycylic.\\ If $ g=0$, then the
moduli space is a point by Riemann's uniformization theorem and the mapping
class group is trivial. Again, the mapping class group is not quite the
correct object to consider, because $\Diff(\bS^2) \simeq \SO(3)$, and $B
\SO(3)$ is certainly not rationally acyclic.\\
The main ingredient for the proof of Theorem \ref{mue} is the following lemma.

\begin{lem}
The mapping class group has a normal torsionfree subgroup of finite index.
\end{lem}

\textbf{Proof:}
We make use of the homomorphism $\rho: \Gamma_g \to \Gl_{2g}(\bZ)$ from
\ref{symplecto}. For $n \in \bN$, the the kernel of the homomorphism
$\Gl_{2g}(\bZ) \to \Gl_{2g}(\bZ/ {n \bZ})$ is the \emph{congruence subgroup}
$\Gamma(n) \subset \Gl_{2g}(\bZ) $. If $3 \leq n$, then $\Gamma(n)$ is
torsionfree (see \cite{Br}, p. 40 f.). Then $\Delta:=\rho^{-1}(\Gamma(n))
\subset \Gamma_g$ is a normal subgroup of finite index. I claim that $\Delta$ is
torsionfree.\\ Let $\gamma \in \Delta$ be a torsion element, $\gamma^r= \id$.
Without loss of generality, we can assume that $r$ is a prime number. Since
$\Gamma(n)$ is torsionfree, $\rho(\gamma)=\id$. Thus, we need to argue that
the kernel of $\rho$, the \emph{Torelli group}, is torsionfree.\\ To this end,
let us show that $\gamma$ (it is a component of the diffeomorphism group)
contains a diffeomorphism $f$ of the same order as $\gamma$. Consider the
bundle

$$\pi:\cV_g:=\cS(M) \times_{\Diff_0(M)} M \to \cT_g.$$

It has an action of the mapping class group $\Gamma_g = \Diff(M)/
\Diff_0(M)$ which covers the action of $\Gamma_g$ on $\cT_g$. If $\gamma$ had
no fixed point on $\cT_g$, then the action of the group $\bZ/r$ generated by
$\gamma$ on $\cT_g$ would be free. Then we would have a finite-dimensional
model for $B \bZ/r$, which is absurd. Thus $\gamma$ has a fixed point $ x \in
\cT_g$. But then $\gamma$ yields an orientation-preserving diffeomorphism $f$
of finite order of the fiber $\pi^{-1} (x)$. By construction, $[f] = \gamma$.
and $f$ is a diffeomorphism of a surface of genus $g \geq 2$, which acts
trivially on $H_1(M, \bZ)$. It is well-known, that $f$ itself must be the
identity. \qed\\

Let $ \Delta \subset \Gamma_g$ be a torsionfree normal subgroup of finite
index, set $G := \Gamma_g /\Delta$. Then we have a commutative diagram

\begin{equation}\label{digramm}
\xymatrix{ B\Delta \ar[r]^-{m}_-{\cong} \ar[d]^{q} & \cT_g / \Delta
\ar[d]^{p}\\ B \Gamma_g \ar[r]^{\mu} & \fM_g\\
}
\end{equation}

In this diagram, $B\Delta$ and $\cT_g/\Delta$ are $G$-spaces and the upper
horizontal map $m$ is an equivariant homotopy equi\-valence (because $\Delta$ is
torsionfree, it must act freely on the finite-dimensional contractible
manifold $\cT_g$). The map $q$ is a finite $G$-covering and $p$ is also a
quotient map by $G$, but with singularities.

\begin{lem}\label{finittt}
Let $X$ be a CW-complex and let $G $ a finite group, which acts cellularly on $X$. Assume that a group element which leaves a cell invariant fixes the cell pointwise. Then there is a natural isomorphism
$H_{*}(X; \bQ) \otimes_{G \bQ} \bQ \cong  H_{*}(X/G; \bQ)$.
\end{lem}

\textbf{Proof:} Under the assumptions on the action, there is a natural
isomorphism

$$C_{*}(X; \bQ) \otimes_{G\bQ} \bQ \cong C_{*}(X/G; \bQ).$$

Because $G$ is finite, $\otimes_{ G \bQ }  \bQ $ is an exact functor. Hence

$$H_{*}(X/G; \bQ) \cong H_{*}(C_{*}(X; \bQ)) \otimes_{G\bQ} \bQ) \cong H_{*}(X;
\bQ) \otimes_{G\bQ} \bQ.$$ \qed\\

\textbf{Proof of Theorem \ref{mue}:}
The assumpions on the action are very mild and can be easily satisfied in our
case. An application of Lemma \ref{finittt} to the vertical maps in the diagram
\ref{digramm} finishes the proof of the theorem.\qed\\

Now I give a construction of the MMM-classes in $H^* (\fM_g, \bQ)$. Recall
the diagram \ref{digramm}

$$
\xymatrix{ B\Delta \ar[r]^-{m}_-{\cong} \ar[d]^{q} & \cT_g / \Delta
\ar[d]^{p}\\ B \Gamma_g \ar[r]^{\mu} & \fM_g\\
}
$$

where $m$ is a homotopy equivalence. $\cT_g/ \Delta$ is a smooth manifold and
$\fM_g$ is not a manifold, but is still an oriented rational homology manifold
of dimension $6g-6$, i.e., for all $x \in \fM_g$, $H_{*}(\fM_g, \fM_g
\setminus \{ x \} ; \bQ) \cong H_{*}(\bR^{6g-6}, \bR^{6g-6} \setminus\{ 0\}; \bQ) $.
Then the proof of Poincar\'e duality in \cite{Bred} applies without any
change, but only for rational coefficients. Hence we can define the
push-forward for proper maps in rational cohomology in the same way as for
smooth manifolds. Now we define $\kappa_{n}^{\prime}:= \frac{1}{[\Gamma_g :
\Delta]} p_{!} (m^{-1})^{*} q^{*} (\kappa_n) \in H^{2n}(\fM_g; \bQ)$. This
does not depend on the choice of the torsionfree normal subgroup $G \subset
\Gamma_g$. One can see that $\mu^{*}(\kappa_{n}^{\prime}) = \kappa_n \in
H^{2n} (B \Gamma_g; \bQ)$. However, the integral MMM-classes do not lie in
the image of $H^{2n}(\fM_g; \bZ) \to H^{2n}(B \Gamma_g; \bZ)$. \\ The reason is
that there exist Riemann surfaces $M$ and cyclic groups $G$ acting on $M$,
such that $\kappa_n(EG \times_G M \to BG) \neq 0$. Moreover, such actions exist
for many values of the genus and $n$. This is shown in the work by Akita,
Kawazumi and Uemura \cite{AKU}.\\ The moduli map $\mu: BG \to \fM_g$ of such a
surface bundle is constant, because all fibers are conformally equivalent by
construction. If there would be a class $c \in H^{2n} (\fM_g; \bZ)$ with
$\mu^{*}(c) = \kappa_n$, then we would have a contradiction.\\ Thus we see
that the integral cohomology of the moduli space is a quite irregular object.
It is usually not studied, because it does not seem to be accessible.

\pagebreak

\section{Spin structures on manifold bundles}\label{spinstrmanibun}

\subsection{Spin structures on vector bundles}\label{spinvector}

\begin{defn}
Let $1 \leq d \in \bN$. The group $\Spin(d)$ is defined to be the unique connected twofold covering group of $\SO(d)$. For $d=0$, we say that $\Spin (0)= \bZ/2$.
\end{defn}

A more explicit definition can be given in terms of Clifford algebras (\cite{LM}). Via the covering map $\Spin(d) \to \SO(d)$, $\Spin(d)$ acts on $\bR^d$.
The group $\Gl^{+}_{d}(\bR)$ also has a unique connected two-fold covering group $\spingroup$. It is well-known that in the commutative diagram 

$$\xymatrix{
\Spin(d) \ar[r] \ar[d] & \spingroup \ar[d] \\
\SO(d) \ar[r]   & \Gl^{+}_{d}(\bR)
}$$

the horizontal maps are inclusions and homotopy equivalences.

\begin{defn}\label{defspin}
Let $E \to X$ be an oriented $d$-dimensional real vector bundle. A \emph{Spin structure} $\sigma$ on $E$ is a pair $\sigma=(P, \alpha)$, where $P$ is a $\spingroup$-principal bundle on $X$ and $\alpha$ is an isomorphism $P \times_{\spingroup} \mb{R}^d \cong E$ of oriented vector bundles on $X$. An \emph{isomorphism} of two spin structures $\sigma_0 =(P_0, \alpha_0)$ and $\sigma_1 =(P_1, \alpha_1)$ is an isomorphism $\beta:P_0 \to P_1$ of $\spingroup$-principal bundles over $\id_X$ , such that $\alpha_1 \circ (\beta \times_{\spingroup} \id_{\bR^d} ) = \alpha_0$. We also say that $\sigma_0$ and $\sigma_1$ are \emph{equivalent}.
\end{defn}

If $f:Y \to X$ is a continuous map and $(P, \alpha)$ a spin structure on $E$ as above, then $(f^{*}P, f^{*}\alpha)$ is a spin structure on $f^{*} E$.\\
The notion of equivalence can be stated in a slightly different way.

\begin{lem}\label{equspin}
Two spin structures $(P_0,\alpha_0)$ and $(P_1, \alpha_1)$ are equivalent, written $(P_0,\alpha_0) \cong (P_1, \alpha_1)$ if and only if they are concordant, i.e. there exists a spin structure $(P, \alpha)$ on $E \times [0,1] \to X  \times [0,1]$, such that $j_{i}^{*}(P,\alpha)=(P_i,\alpha_i)$ for $i=0,1$, where $j_i: X \to X \times [0,1]$ is the inclusion $x \mapsto (x,i)$.
\end{lem}

\textbf{Proof:}
If the two spin structures are equivalent, then we can use the isomorphism $\beta:P_0 \to P_1$ to glue the bundles $P_0 \times[0, \frac{1}{2}]$ and $P_1 \times[ \frac{1}{2},1]$ together. The result is a concordance between the two spin structures.\\
Conversely, if a concordance is given, then an isomorphism is constructed using the covering homotopy theorem.\qed\\

The usual definition of a spin structure assumes the existence of a Riemannian metric on the vector bundle $E$. Then a spin structure is a $\Spin(d)$-principal bundle with analogous properties as in our definition.
Of course, a vector bundle has a spin structure in our sense if and only if, after choosing a Riemannian metric, the resulting $\SO(d)$-vector bundle has a spin structure in the usual sense.
For most topological problems, one can always assume that a vector bundle comes with a metric, because the space of scalar products on $\bR^d$ is contractible. In our applications, this assumption would be too restrictive, see \ref{spinmfd} and \ref{spinspace}.\\
In differential geometry, there is a more serious reason why one defines spin structures only for Riemannian vector bundles. It is a fact from basic Lie group theory that any homomorphism $\spingroup \to \Gl_n(\bC)$ factors through $\Gl_d (\bR)$. In contrast, there are important representations of $\Spin(d)$ which do not factor through $\SO(d)$ and they are used to construct the spinor bundles. It follows that the spinor bundles cannot be defined without specifying a metric - this is the main psychological difficulty in spinor geometry.\\
There is another description of spin structures, which is sometimes convenient. If $Q \to X$ is the bundle of oriented frames in $E$, then we can also say that a spin structure is a $\spingroup$-principal bundle, together with a $\spingroup$-equivariant map $P \to Q$ over $X$. Yet another way of saying this is that $P$ is a twofold covering of $Q$, whose restriction to any fiber of $Q$ is nontrivial.\\ 
The next theorem is a well-known standard result (\cite{LM}, p.78-85). Its proof is an easy application of obstruction theory.

\begin{thm}\label{wzwei}
An oriented vector bundle $E$ has a spin structure if and only if $w_2 (E)=0$. If a vector bundle has a spin structure, then the set of equivalence classes of spin structures is in bijection with $H^{1}(X;\mb{Z}/2)$. More precisely, it is a $H^{1}(X;\mb{Z}/2)$-torsor (i.e. a principal homogeneous space).\\
Furthermore, if $E_2=E_0 \oplus E_1$ as oriented vector bundles, then the choice of spin structures on any of the three bundles determines a spin structure on the third.
\end{thm}

Because $H^1(X,  \bF_2)$ is in bijection with the set of equivalence classes of real line bundles, there is an implicit simply-transitive action of the set of real line bundles on the set of spin structures. It is worth to make this explicit. I do this in two ways, one algebro-topological and one geometrical.
Let $E \to X$ be an oriented vector bundle and let $Q\to X$ be the associated frame bundle. For simplicity of notation, we choose a Riemannian metric and talk about $\SO(d)$- and $\Spin(d)$-principal bundles.\\
Recall that $H^1(\SO(d); \bF_2)\cong \bF_2$ (in the sequel, all cohomology groups are with $\bF_2$-coeffi\-cients) and that $H^1(B\SO(d))=0$; $H^2(B\SO(d)) = \bF_2 \langle w_2 \rangle$. Consider the Leray-Serre spectral sequence for $\bF_2$-cohomology of the universal bundle $\SO(d) \to E\SO(d) \to B\SO(d)$. The differential $d_2: H^1(\SO(d)) = E_{2}^{0,1} \to E_{2}^{2,0} = H^2 (B\SO(d))$ is an isomorphism and sends the nontrivial element to the universal Stiefel-Whitney class. By naturality of the spectral sequence it follows that for an arbitrary $\SO(d)$-bundle $Q  \to X$, $X$ connected, the differential 

$$d_2: H^1 (SO(d)) = E_{2}^{0,1} \to E_{2}^{0,1}= H^2 (X)$$
 
sends the nontrivial element in $H^1 (SO(d))$ to $w_2 (Q)$. From the Leray-Serre spectral sequence, one deduces the Serre exact sequence

$$0 \to H^1 (X) \to H^1 (Q) \to H^1 (\SO(d)) \to H^2 (X),$$

where the last map is $d_2$ and $d_2 =0$ if $Q$ is spin. From the description of spin structures as coverings of $Q$, we see that the set of equivalence classes of spin structures can be identified with the nontrivial coset in $H^1 (Q)$. The $H^1 (X)$-action is via addition.\\
Now I turn to the geometrical description. If $\sigma_i= (P_i, \alpha_i)$, $i=0,1$ are two spin structures on $E$, then we can construct a real line bundle (or rather a $\bZ/2$-principal bundle) $L=L_{P_1, P_{0}^{-1}}$ on $X$ as follows. 
Let $x \in X$. Then the fiber $L_x$ is the set of all $\spingroup$-equivariant maps $a:(P_0)_x \to (P_1)_x$ which give the identity on $E_x$ in the sense that the induced map 

$$(\alpha_1) \circ (a \times_{\spingroup} \bR^d) \circ (\alpha_0)^{-1} :E_x \to (P_0)_x \times_{\spingroup} \bR^d \to (P_1)_x \times_{\spingroup} \bR^d \to E_x$$

is the identity. It is easy to check that $L$ is a $\bZ/2$-principal bundle; the notation is self-explanatory.\\
Conversely, let $\sigma=(P, \alpha)$ be a spin structure and let $L \to X$ be a $\bZ/2$-principal bundle.
To define the spin structure $L \sigma$, we consider the spin structure as a twofold covering of $Q$. This determines a real line bundle on $Q$, which is nontrivial when restricted to any fiber. If $L \to X$ is a line bundle and $q:Q \to X$ is the bundle projection, then we can tensorize the spin structure with the line bundle $q^* L$. One can check that this defines an action of the set of line bundles on the set of equivalence classes of spin structures and that $L_{P_1 P_{0}^{-1}} P_0 = P_1$. \\
For complex line bundles, there is yet another description of spin structures. If $L$ is a complex line bundle, then a spin structure can be viewed as another complex line bundle $S$, together with a complex isomorphism $S \otimes_{\bC} S \to L$. This is easy to see; one uses that $\Spin(2) \cong \SO(2)$ and that the covering map is the squaring map.

\subsection{Spin structures on manifolds}\label{spinmfd}

Now we let $M$ be a smooth oriented $d$-dimensional manifold. If not stated otherwise, all diffeomorphisms are assumed to be orientation-preserving. We say that a \emph{spin structure on} $M$ is a spin structure on the tangent bundle $TM$ of $M$. A \emph{spin manifold} is an oriented manifold, equipped with a spin structure.\\
Usually, one defines a spin structure on a manifold to be spin structure of the \emph{cotangent} bundle. This is very convenient in differential and algebraic geometry. Of course, spin structures on the tangent and on the cotangent bundle are closely related. If $(P, \alpha)$ is a spin structure for the tangent bundle, then $(P, (\alpha^{-1})^{\prime}  )$ is a spin structure for the cotangent bundle, where $(\alpha^{-1})^{\prime} $ is the dual of the inverse of $\alpha$.\\
Orientation-preserving diffeomorphism or, more generally, orientation-preser\-ving local diffeomorphisms can be used to pull-back spin structures on manifolds. Let $f: N \to M$ be an orientation-preserving local diffeomorphism. Since $f$ is regular, its differential defines an orientation-preserving isomorphism $df: TN \to f^{*}TM$ of bundles on $N$. If $\sigma:=(P, \alpha)$ is a spin structure on $M$, then we define $f^{*} \sigma:= (f^{*}P, (df)^{-1} \circ f^{*}\alpha)$. It is a spin structure on $N$.\\
If $f_0$ and $f_1$ are regularly homotopic, i.e. they are homotopic through local diffeomorphisms, then the spin structures $f_{0}^{*}\sigma$ and $f_{1}^{*}\sigma$ are equivalent. This follows immediately from Lemma \ref{equspin} and from the definitions.\\
Thus we see that the group $\Gamma(M):=\pi_0 (\Diff(M))$ of isotopy classes of diffeomorphisms acts from the right on the set of equivalence classes of spin structures. Although the set of equivalence classes of spin structures is in bijection with $H^1 (M; \bF_2)$, there is no $\Gamma(M)$-equivariant identification in general, compare the remarks after \ref{wzwei}.\\
We call two spin structures $\sigma_1$ and $\sigma_2$ \emph{conjugate} if there exists a diffeomorphism $f$ such that $\sigma_1 \cong f^{*} \sigma_2$.\\
If $\sigma$ is a spin-structure on $M$, then we denote

\begin{equation}\label{spindiffff}
\Diff(M,\sigma):= \{ f \in \Diff(M)| f^{*}\sigma \cong \sigma  \}.
\end{equation}

Here the symbol $\cong$ denotes equivalence of spin structures (and not conjugacy). 
More generally, if $(M, \sigma)$ and $(N, \tau)$ are spin manifolds, then $\Diff((M, \sigma),(N, \tau))$ denotes the set of all (orientation-preserving) diffeomorphisms $f:M \to N$ such that $f^{*} \tau \cong \sigma$. The preceding discussion, together with the fact that regular homotopy is an open equivalence relation on $\Diff(M)$ shows

\begin{lem}\label{spinnlemma}
$\Diff(M,\sigma) \subset \Diff(M)$ is an open subgroup.
\end{lem}

At the first glance, it might seem that the group $\Diff(M, \sigma)$ is the correct automorphism group of the spin structure $\sigma$. This is not the case, as I will explain now with an easy example.\\
Consider the $2$-sphere with the standard metric. As an oriented surface, $\bS^2$ is a spin manifold. The spin structure is unique up to equivalence, because $H^1(\bS^2, \bF_2)=0$. Thus, we certainly have $\Diff(\bS^2, \sigma) = \Diff(\bS^2)$.
A concrete description of the spin structure is as follows. Consider the Hopf fibration $\bS^3 \to \bS^2$ which is a $\SO(2)=\Spin(2)$-principal bundle. Then there is an isomorphism $\bS^3 \times_{\SO(2)} \bR^2 \cong T\bS^2$, where $\SO(2)$ acts on $\bR^2 = \bC$ by the \emph{square} of the standard representation. Thus the Hopf fibration is a spin structure for $\bS^2$. But the group $\SO(3) \subset \Diff(\bS^2)$ does \emph{not} act on the total space of this spin structure.\\
Thus $\Diff(M, \sigma)$ is not the correct thing to consider, especially if we attempt to study \emph{families} of spin manifolds.
This requires that we do not merely consider the discrete \emph{set} of all equivalence classes of spin structures, but a \emph{space} of all spin structures on $M$, which has a nontrivial topology. This will be done in the next section.

\subsection{Spaces of Spin structures}\label{spinspace}

\begin{defn}\label{spiccc}
Let $M$ be an oriented manifold. Then the category $\Sspin (M)$ is defined as follows. Its objects are all spin structures on $M$. A morphism between two spin structures $\sigma_0$ and $\sigma_1$ is an isomorphism in the sense of Definition \ref{defspin}.
\end{defn}

Because we will consider the classifying space of this category, one should take a small skeleton of this category instead. Such a small skeleton certainly exists. The classifying spaces of any two skeleta of a category are homotopy equivalent through an equivalence which is canonical up to homotopy. The question whether the space of all these canonical homotopy equivalences is contractible is irrelevant in this work.
As usual, a (topological) group is identified with a (topological) groupoid with one object.

\begin{prop}\label{easything}
If $M$ is connected, the classifying space $B\Sspin(M)$ of the category $\Sspin(M)$ is homotopy equivalent to $H^1(F; \bF_2) \times \bR \bP^{\infty}$.
\end{prop}

\textbf{Proof:} We have seen that the set of equivalence classes of spin structures on $M$ is in bijection with $H^1(M; \bF_2)$. We choose one spin structure in each equivalence class. This gives a set $R$ of representatives for the equivalence classes of spin structures. It was remarked above that the automorphism group of a spin structure is isomorphic to $\bZ/2$. An application of Lemma \ref{groupoid} from the appendix \ref{appendiy} finishes the proof.\qed\\

There is another category of spin structures which also involves the action of diffeomorphisms on spin structures. The definition needs some preparation.\\
Let $M$ and $N$ be smooth oriented $d$-manifolds, let $\sigma=(P,\alpha)$ be a spin structure on $TM$ and $\tau=(Q,\beta)$ be a spin structure on $N$.
If $\hat{f}:P \to Q $ is a map of $\spingroup$-principal bundles over a base map $f:M \to N$, then $\hat{f}$ induces a vector bundle map $\hat{f}^{\prime}: P \times_{\spingroup} \bR^d \to Q \times_{\spingroup} \bR^d$ over $f$, which can be compared to the differential of $f$. The vector bundle map $\hat{f}^{\prime}$ has a chance to agree with the differential only if $f$ is a local diffeomorphism.

\begin{defn}
A \emph{local spin diffeomorphism} from $(M, \sigma) \to (N, \tau)$ is a pair $(f, \hat{f})$, where $f$ is an orientation-preserving local diffeomorphism (alias codimension zero immersion) and $\hat{f} : P \to Q$ is an isomorphism of $\spingroup$-principal bundles over $f$, such that the diagram
 
$$\xymatrix{
P \times_{\spingroup} \bR^d \ar[r]^{\hat{f}^{\prime}} \ar[d]^{\alpha} & Q \times_{\spingroup} \bR^d \ar[d]^{\beta} \\
TM \ar[r]^{Tf} & TN
}
$$

of vector bundle maps commutes.
Local spin diffeomorphisms can be composed in the obvious way: $(g, \hat{g}) \circ (f, \hat{f}):= (g \circ f,  \hat{g} \circ \hat{f})$.\\
A local spin diffeomorphism $(f, \hat{f})$ is called a \emph{spin diffeomorphism} if the underlying map $f$ is a diffeomorphism.
\end{defn}

We can say shortly that a spin diffeomorphism is a diffeomorphism $f$ which preserves the spin structure $\sigma$, together with a choice of an isomorphism of spin-structures $\sigma \cong f^{*} \sigma$. It is easy to see that the set of spin diffeomorphisms of a spin structure $\sigma$ on $M$ forms a \emph{group}.

 \begin{defn}\label{spinnnnnn}
The group of spin diffeomorphism of $(M, \sigma)$ is denoted by $\Spiff(M, \sigma)$.\\
More generally, let $M$ be an oriented manifold. Then $\Spiff(M)$ is the following groupoid:
\begin{itemize}
\item The objects are the spin structures $\sigma$ on $M$.
\item The morphisms $\sigma_0 \to \sigma_1$ are the spin diffeomorphisms $(M, \sigma_0 ) \to (M, \sigma_1)$.
\end{itemize}
The morphism sets between different spin structures are denoted by $\Spiff((M,\sigma_0) ,(M,\sigma_1))$.
\end{defn}

If $M$ does not admit a spin structure, then $\Spiff(M)$ is the empty groupoid, but we will not consider this case. Note that the group $\Spiff(M, \sigma)$ is the automorphism group of\\ $\sigma \in \Ob (\Spiff(M))$. 
There is an obvious forgetful map 

\begin{equation}\label{psidefn}
\psi:\Spiff(M, \sigma) \to \Diff(M, \sigma); \psi: (F, \hat{f}) \mapsto f
\end{equation}
 
which is surjective.
I will define a topology on the morphism sets of $\Spiff(M)$ and thus on $\Spiff(M, \sigma)$ below.\\
By definition, the group $\Spiff(M, \sigma)$ acts on the principal bundle $P$; and the action covers the action of the diffeomorphism group. More precisely, we have:

\begin{lem}\label{spifflemma}
If $M$ is a connected manifold, then the kernel of $\psi:\Spiff(M, \sigma) \to \Diff(M, \sigma)$ has two elements.
\end{lem}

\textbf{Proof:}
If $(f, \hat{f}) \in \ker \psi$, then $f=\id$. Let $(\id_M,\hat{f})$ be a spin diffeomorphism, and let $c$ be the nontrivial element in the kernel of $\spingroup \to \Gl_{d}^{+} (\bR)$. Then $\hat{f}$ is an automorphism of $P$ over $\id_M$, and the condition that $(\id_M \hat{f})$ is a spin diffeomorphism means that either $\hat{f}(p)=p$ or $\hat{f}(p)= p \cdot c$ for all $p \in P$ ($M$ is connected). In the first case $\hat{f}= \id_P$, in the second case, $\hat{f}$ is right-multiplication with $c$.\qed\\

Now I will explain the topology on $\Spiff(M, \sigma)$. Any space of the type $C^{\infty}(X;Y)$ of smooth mappings between smooth manifolds $X$ and $Y$ is endowed with the weak Whitney-topology (\cite{Hirsch}, p. 34 ff.). This is the topology of uniform convergence of all derivatives on all compact subsets of $X$. \\
If $\sigma=(P,\alpha)$ and $\tau =(Q, \beta)$ are spin structures on $M$, $N$, respectively, then there is an injection $\Spiff((M,\sigma), (N, \tau)) \subset C^{\infty}(P;Q)$. We let $\Spiff((M,\sigma), (N, \tau)$ carry the subspace topology. With this definition, the composition and inversion maps in the groupoid $\Sdiff(M)$ become continuous and the forgetful maps 

$$\Spiff((M,\sigma),(N, \tau)) \to \Diff((M, \sigma),(N,\tau)$$

become covering maps. There is a forgetful functor of topological groupoids

\begin{equation}\label{Psidefn}
\Psi:\Spiff(M) \to \Diff(M).
\end{equation}

\begin{thm}\label{homofiber}
The homotopy fiber of $B \Psi :B \Spiff(M) \to B \Diff(M)$ is homotopy equivalent to $B \Sspin (M) \simeq H^1(M; \bF_2) \times \bR \bP^{\infty}$.
\end{thm}

\textbf{Proof:} We apply Proposition \ref{groupoid2} from the appendix \ref{appendiy} to the groupoid homomorphism $\Psi: \Spiff(M) \to \Diff(M)$. Note that $\Sspin(M)$ is the kernel (in the sense of groupoids, see appendix \ref{appendiy}) of $\Psi$ by definition. Furthermore, $\Phi$ is surjective (in the sense of groupoids): If $f \in \Diff(M)$ and $\sigma \in \Ob (\Spiff(M))$, then there is a spin diffeomorphism $h:\sigma \to (f^{-1})^* \sigma$ with $\Psi(h)=f$. The map on morphism spaces is a Serre fibration, because it is even a covering. Apply Proposition \ref{groupoid2}.\qed\\

Fix a spin structure $\sigma$ on $M$. We can compare the two fibrations

\begin{equation}\label{rewq}
\xymatrix{ 
\bR \bP^{\infty} \ar[r] \ar[d]& B \Sspin(M) \ar[d]\\
B \Spiff(M, \sigma) \ar[d] \ar[r] & B \Spiff(M) \ar[d]\\
B \Diff(M, \sigma) \ar[r] & B \Diff(M).\\
}
\end{equation}

Their difference is not as much as it seems at the first glance. The spin structure $\sigma$ can be viewed as a point in the zero-skeleton of both classifying spaces, $B \Sspin (M)$ and $B \Spiff(M)$. If $\sigma$ is taken as a base-point, then all maps in the diagram \ref{rewq} are pointed. We write for short $\Gamma(M):= \pi_0 (\Diff(M))$ and $\Gamma(M, \sigma):= \pi_0 (\Diff(M, \sigma))$. Note that $\Gamma(M, \sigma)$ is a subgroup of $\Gamma(M)$ because of \ref{spinnlemma}. \\
Also because of \ref{spinnlemma}, $\pi_r (B \Diff(M, \sigma)) \to\pi_r ( B \Diff(M))$ is an isomorphism if $r \geq2$.\\ Thus, for $r \geq 3 $, the diagram \ref{rewq} gives rise to a commutative square of isomorphisms (the base-point is omitted from the notation)

$$\xymatrix{
\pi_r (B \Spiff(M, \sigma)) \ar[d] \ar[r] & \pi_r ( B \Spiff(M) ) \ar[d] \\
\pi_r (B \Diff(M, \sigma)) \ar[r] & \pi_r (B \Diff(M)).\\
}
$$

It remains to analyze the low-dimensional terms. They are

\begin{equation}\label{ewq}
\xymatrix{ 
0 \ar[r] \ar[d]& 0 \ar[d] \\
\pi_2 (B \Spiff(M, \sigma)) \ar[d] \ar[r] &\pi_2( B \Spiff(M)) \ar[d] \\
\pi_2(B \Diff(M, \sigma)) \ar[r] \ar[d] & \pi_2( B \Diff(M)) \ar[d] \\
\bZ/2 \ar[r] \ar[d]& \bZ/2 \ar[d]\\
\pi_1 (B \Spiff(M, \sigma)) \ar[d] \ar[r] &\pi_1( B \Spiff(M)) \ar[d]\\
\Gamma(M, \sigma) \ar[r] \ar[d] & \Gamma(M) \ar[d] \\
\pt \ar[r] \ar[d]& \pi_0(B \Sspin(M)) \ar[d]\\
\pt \ar[d] \ar[r] & \pi_0(B \Spiff(M))=\pi_0 (B  \Sspin (M)) / \Gamma(M) \ar[d]\\
\pt \ar[r] & \pt.\\
}
\end{equation}

\subsection{Spin structures on fiber bundles}\label{spinfiber}

Next, we study the question: What does the classifying space of the groupoid $\Spiff(M)$ classify?\\
We consider smooth $M$-bundles $(E,B, \pi)$, where the oriented manifold $M$ admits a spin structure. To a smooth $M$-bundle, there is an associated $\Diff(M)$-principal bundle $Q \to B$ and a natural isomorphism $Q \times_{\Diff(M)} M \to E$ of fiber bundles. The base space $B$ is not assumed to be a manifold, but a "reasonable" space. With this term, I mean that $B$ is paracompact (such that $G$-principal bundles on it are classified by maps into $BG$), locally arcwise connected and locally relatively simple connected (such that there is a well-behaved theory of covering spaces, \cite{Bred}, p.155).\\

\begin{defn}
A spin structure on an oriented smooth $M$-bundle $ \pi:E \to B$ is a spin structure on the vertical tangent bundle $T_v E$.
A smooth $M$-spin bundle is a smooth $M$-bundle endowed with a spin structure.
\end{defn}

Note that any fiber of a smooth $M$-spin bundle becomes a spin manifold.

\begin{defn}
Let $M$ be a manifold with $w_2 (TM)=0$ and let $\sigma$ be a spin structure on $M$. For a topological space $B$, let $\Bun_{M, \Spin} [B]$ be the set of all concordance classes of smooth $M$-bundles with spin structure. Let $\Bun_{M, \sigma} [B]$ be the set of all concordance classes of smooth $M$-bundles with spin structure such that the restriction of the spin structure to any fiber is conjugate to $\sigma$.
\end{defn}
 
\begin{thm}
For any reasonable space $B$, there are bijections, natural in $B$:
\begin{enumerate}
\item $[B; B \Spiff(M, \sigma) ] \cong \Bun_{M, \sigma} [B]$ and\\
\item $[B; B \Spiff(M) ] \cong \Bun_{M, \Spin} [B]$.
\end{enumerate}
\end{thm}

\textbf{Proof:}
Let us prove the first claim. Because $\Spiff(M, \sigma)$ is a group, we only need to relate $\Spiff(M, \sigma)$-principal bundles on $B$ with $\Bun_{M, \sigma} [B]$.
If $\sigma=(P, \alpha)$, then $\Spiff(M, \sigma)$ acts on $P$, as well as on $\Gl^{+} (M)$ and the map $P \to \Gl^{+} (M)$ is equivariant. Let $Q \to B$ be a $\Spiff(M, \sigma)$-principal bundle. Then 

$$E:=Q \times_{\Spiff(M, \sigma)} M \to B$$

is a smooth $M$-bundle and 

$$Q \times_{\Spiff(M, \sigma)} P,$$

together with the isomorphism

$$Q \times_{\Spiff(M, \sigma)} P \times_{\spingroup} \bR^d \to Q \times_{\Spiff(M, \sigma)} TM = T_v E$$

is a spin structure on $E$. Conversely, if $E \to B$ is a smooth $M$-bundle with spin structure such that for all $b \in B$, $E_b $ is spin diffeomorphic to $(M,\sigma)$, then we define a $\Spiff(M, \sigma)$-principal bundle $Q \to B$ as follows. The fiber $Q_b$ over $b \in B$ consists of all spin diffeomorphisms $(M, \sigma) \to E_b$. Both constructions are mutually inverse. This shows the first claim.\\
For the proof of the second claim, we assume without loss of generality that $B$ is connected. We choose a system of representatives $R$ for the set of conjugacy classes of spin structures on $M$ and \emph{define} a functorial map $[B; B \Spiff(M) ] \to \Bun_{M, \Spin} [B]$ by the commutativity of the diagram

$$\xymatrix{
\coprod_{\sigma \in R} [B; B \Spiff(M, \sigma) ] \ar[r]^{\cong} \ar[d]  &  [B; B \Spiff(M) ] \ar[d]\\
\coprod_{\sigma \in R}  \Bun_{M, \sigma} [B]  \ar[r] & \Bun_{M, \Spin} [B] .\\
}
$$

The upper horizontal arrow is bijective because of Lemma \ref{groupoid}. The right-hand side vertical arrow is provided by the first claim and bijective by part 1 of this proof. The lower horizontal arrow is bijective because $R$ is a system of representatives for the set of conjugacy classes of spin structures.\\
This shows the existence of the functorial bijection. Using the fact that the map $BG \to BG$ induced by an inner automorphism of the topological group $G$ is homotopic to the identity, one can show that the map
$[B; B \Spiff(M) ] \to \Bun_{M, \Spin} [B]$ defined above does not depend on the choice of the system of representatives $R$.\qed\\

\subsection{When does a smooth fiber bundle have a spin structure?}\label{wdafbhass}

Now, let us consider a smooth $M$-bundle $\pi:E \to B$ on a reasonable space. We try to find necessary and sufficient conditions for the existence of a spin structure on $E \to B$. First of all, $M$ needs to satisfy $w_2(TM)=0$, and we assume that this is the case.\\ 
One is tempted to apply obstruction theory to the fibration 

$$B \Sspin (M) \to B\Spiff(M) \to B\Diff(M)$$

from the last section. Unfortunately, this does not work, because the situation is highly non-abelian. There is no group structure on $\pi_0 (B \Sspin (M))$ which has a relation to spin structures (even if $B \Sspin (M)$ is abstractly homotopy-equivalent to an abelian topological group).\\
The situation is more tractable if we consider the fibration sequence 

\begin{equation}\label{erweierubg}
\bR \bP^{\infty} \to B\Spiff(M, \sigma) \to B\Diff(M, \sigma).
\end{equation}

To start the discussion, we assume (without loss of generality) that $B$ is path-connected and that we have a base point $b_0 \in B$. Further, we fix an orientation-preserving diffeomorphism $ \phi:M \to \pi^{-1}(b_0)$. Such a triple $(E, \pi, \phi)$ should be called a \emph{pointed smooth $M$-bundle}.\\
Concordance classes of pointed $M$-bundles are in bijection with pointed homotopy classes of pointed maps $(B,b_0) \to (B \Diff(M), \pt)$.\\
On fundamental groups, the classifying map of a pointed $M$-bundle induces a homomorphism 

$$\rho_{(E, \pi, \phi)}: \pi_1 (B, b_0) \to \pi_0 \Diff(M),$$

the \emph{monodromy} of the manifold bundle. There is the following explicit description of the monodromy homomorphism. Let $\gamma: \bS^1 \to B$ be a pointed loop. Then $\gamma^{*} E=:E^{\prime}$ is a bundle of manifolds on $\bS^1= [0,1]/(0 \sim 1)$ and the fiber at the base point $ 0 \in \bS^1$ is canonically diffeomorphic to $M$. Let $q: [0,1] \to \bS^1$ be the identification map. Then choose a trivialization $ a:q^{*}E^{\prime} \cong [0,1] \times M$, extending the identity over $0$. At the point $1$, $a$ defines a diffeomorphism $f:M \to M$. Its isotopy class does not depend on the choice of $ \gamma$ within its homotopy class and on the choice of $q$. The isotopy class of $f$ is $\rho([\gamma])$. Note that the monodromy does \emph{not} lift to a homomorphism $\pi_1(B) \to \Diff(M)$ in general; this only happens for flat fiber bundles.\\

\begin{lem}\label{monodromya}
If there is a spin structure on $E$ which extends the spin structure $\phi^{*} \sigma$ on $\pi^{-1}(b_0)$, then $\rho$ takes values in the subgroup $\pi_0 (\Diff(M, \sigma))$.
\end{lem}

\textbf{Proof:} Restrict the bundle $E$ and the spin structure on $E$ to a circle representing $a \in \pi_1(B, b_0)$. Then the spin structures $\sigma$ and $\rho(a)^{*} \sigma $ are concordant. By Lemma \ref{equspin}, they are equivalent.\qed\\

Consider the fibration sequence \ref{erweierubg}. It is necessarily a principal fibration in the sense of homotopy theory: the action of $\pi_1 (B\Diff(M, \sigma))$ on the homotopy type of $\bR \bP^{\infty}$ is trivial. Thus the next definition is meaningful.

\begin{defn}
Let $c \in H^2(B \Diff(M, \sigma); \bF_2)$ be the primary obstruction class to the existence of a section of the fibration \ref{erweierubg}.
\end{defn}

The class $c$ also can be described as the cohomology class of the (topological) group extension $\bZ/2 \to \Spiff(M, \sigma) \to \Diff(M, \sigma)$.\\
Now we can reformulate the question of the existence of a spin structure as a lifting problem.

\begin{prop}\label{mainthm} Let $\pi: E \to B$ be a smooth $M$-bundle, let $\phi: \pi^{-1} (b_0) \to M$ be a fixed diffeomorphism and let $\lambda: B \to B \Diff(M)$ be a classifying map. Assume that there exists a spin structure $\sigma$ on $M$ such that the monodromy homomorphism takes values in the subgroup $\pi_0(\Diff(M, \sigma))$. Then there exists a spin structure on $E$ which extends $\phi^{*} \sigma$ if and only if there exists a lift

$$\xymatrix{
 & B   \Spiff(M, \sigma) \ar[d]\\
B \ar[ur] \ar[r]_-{\lambda} & B\Diff(M, \sigma).
}$$

Homotopy classes of lifts (as lifts) are in $1-1$ correspondence with equivalence classes of spin structures.
\end{prop}

The proof is clear.

\begin{cor}\label{obstructionclasas}
With the notation as above, the following statement holds.\\
There exists a spin structure on the fiber bundle $E$ extending $\phi^{*} \sigma$ if and only if the mono\-dromy homomorphism takes values in $\pi_0 (\Diff(M, \sigma)$ and $\lambda^{*}c=0 \in H^2 (B; \bF_2)$.
\end{cor}

\textbf{Proof:}
This follows from Theorem \ref{mainthm} and Lemma \ref{spifflemma} by an easy application of obstruction theory.\qed\\

In the section \ref{spinexpls}, we will show that the cohomology class $c$ is nonzero if $M$ is an oriented surface and $\sigma$ an arbitrary spin structure on $M$.

\pagebreak

\section{Spin structures on surface bundles}\label{spinstrsurbun}

\subsection{Spin structures on surfaces}\label{spinsurfaces}

In this section, we recollect the facts about spin structures on oriented closed surfaces. They are mostly taken from the articles by Atiyah (\cite{At}) and Johnson (\cite{John}). An oriented closed surface $M$ of genus $g$ has even Euler number $2-2g$. Thus, the second Stiefel-Whitney class of the tangent bundle is tri\-vial, and the tangent bundle has a spin structure. Thus an oriented surface is spin. The number of inequivalent spin structures on a closed surface of genus $g$ is $2^{2g}= \sharp H^1(M; \bF_2)$.\\
Spin structures on surfaces have different descriptions. If $M$ has a Riemannian metric, hence a complex structure with a hermitian metric, then $TM$ is a hermitian line bundle in a canonical way. It was said in section \ref{spinvector} that a spin structure on a hermitian line bundle $L \to X$ is nothing else than another hermitian line bundle $S \to X$, together with an isomorphism $S \otimes_{\bC} S \to L$.\\
Now, if we choose a complex structure on $M$, then $TM$ is even a \emph{holomorphic} line bundle. Then a spin structure is also a holomorphic line bundle, in a natural way: call a local section $s$ of $S$ holomorphic if its square $s \otimes s$ is a holomorphic section of $TM$. By the existence of holomorphic square-roots, local holomorphic sections of $S$ without zeroes always exist. Thus $S$ can be considered as a holomorphic line bundle, and the reference isomorphism $S \otimes S \to TM$ is holomorphic. The same discussion applies also to the cotangent bundle, and for the purposes of complex analysis it is more common to consider the holomorphic cotangent bundle $\Lambda_M$ of $M$ as the basic object.\\
Furthermore, if $M$ is compact, then a holomorphic isomorphisms between two isomorphic holomorphic line bundles is unique up to multiplication with a constant in $\bC^{\times}$. Thus the space of possible reference isomorphisms is connected. It follows that, as long as we consider spin structures on a single surface, we can forget about the reference isomorphism. Hence we can say:

\begin{prop}
A spin structure on a compact Riemann surface $M$ is nothing else than a square-root $S$ of the element $\Lambda_M \in \Pic (M)$ (the holomorphic Picard group of the Riemann surface $M$, see \cite{GH}, p. 133) representing the canonical line bundle (alias cotangent bundle).\qed
\end{prop}

However, if we attempt to study families of Riemann surfaces, the nontriviality of the topology of the space of these reference maps becomes important. But we already set up the machine to handle this situation in the last sections. \\
There is an important topological invariant of spin structures on surfaces, its \emph{parity}, which we will call \emph{Atiyah-invariant}, because it was studied by Atiyah in \cite{At}.
Let $M$ be a Riemann surface and let $S$ be a spin structure in the sense above, i.e. a square root of the canonical line bundle $\Lambda_M$ of $M$. Then we have the Cauchy-Riemann operator $\delbar_{S}: C^{\infty}(M; S) \to C^{\infty}(M; S \otimes \overline{\Lambda_M})$. By the Riemann-Roch-theorem, its Fredholm index is zero.\\

\begin{defn}\label{atiyahinv}
The \emph{Atiyah-invariant} of the spin structure $S$ is 

$$At (M;S):= \dim \ker \delbar_{S} \pmod{2} \in \bZ/2.$$

A spin structure is called \emph{even} (\emph{odd}) if its Atiyah-invariant is even (odd).
\end{defn}

It is not at all obvious that $At(M;S)$ is a topological invariant, i.e. invariant under deformations of the spin structure. However, Atiyah proved an even stronger assertion:

\begin{prop}\label{atiyahstheorem}(\cite{At})
$At (M;S)$ is invariant under spin-bordism.
\end{prop}

I will sketch the proof. First, one needs to identify the quantity $\dim \ker \delbar_{S} \pmod{2}$ as a mod 2 index of the Dirac operator associated to the spin structure (see \cite{LM}). This index lies in the group $KO^{-\dim_{\bR} M}= KO^{-2}(*) \cong \bZ/2$. The Atiyah-Singer index theorem for families of real elliptic operators (\cite{AS5}) expresses this index as a $KO$-characteristic number of the spin manifold $M$ (Recall that a spin vector bundle is $KO$-oriented). As such a characteristic number, it is an invariant under spin bordism.\\

\begin{example}\label{spinstrtours} Let us consider the spin structures on a holomorphic torus. The holomorphic cotangent bundle of a torus is the trivial holomorphic line bundle. Thus the trivial holomorphic line bundle can also be considered as a spin structure $S_0$. It has Atiyah-invariant $1$.\\
There are three other spin structures. All of them have Atiyah-invariant zero, because a holomorphic line bundle of degree zero cannot have nonzero holomorphic sections unless it is trivial.
\end{example}

There is another description of the Atiyah-invariant which does not involve the choice of a complex structure. Also, the treatment of the properties does not need index theory, but only elementary homotopy constructions. We start with the following observation:

\begin{prop}
The spin bordism group $\Omega^{\Spin}_{1}$ of $1$-manifolds has two elements. The nontrivial element is represented by the circle, together with the trivial $O(1)$-bundle on it, which is a spin structure.
\end{prop}

For bordism-theoretical constructions, it is much more convenient to consider a spin structure of a manifold to be a spin-structure on its stable \emph{normal} bundle. Both concepts are equivalent. A normal spin structure on a manifold $M$ can be restricted to any submanifold $N \subset M$, as long as the normal bundle of $N$ in $M$ is trivialized (or, at least, endowed with a spin structure).\\
Now let $M$ be a surface with a spin structure $\sigma$. For any embedded oriented $1$-dimensional submanifold $N \subset M$, the normal bundle has a preferred homotopy class of trivializations and so we can restrict the spin structure of $M$ to $N$. The bordism class of this spin structure does not depend on the orientation of $N$.
We define a map 

$$Q= Q^{\sigma}: H_1(M; \bF_2) \to \bF_2$$
 
by the following procedure. Let $x \in H_1 (M; \bF_2)$. Choose a $1$-dimensional submanifold $N \subset M$ representing $x$ and an orientation of $N$. Then set 

$$Q(x):=\begin{cases}
0     & \text{if the spin manifold $N$ is nullbordant, }\\
1      & \text{otherwise.}
\end{cases}
$$

The main theorem in the article \cite{John} by Johnson is

\begin{thm}(\cite{John})\begin{enumerate}
\item The map $Q$ is well-defined, i.e. it does not depend on the choices involved.
\item $Q$ is a quadratic form, i.e. $Q(a+b)= Q(a) +Q(b) + a \cdot b $, where $(a,b) \mapsto a\cdot b$ denotes the intersection form on the homology group.
\end{enumerate}
\end{thm}

The next property is obvious from the definitions. If $f$ is a diffeomorphism, then it acts on the spin structures as well as on the first cohomology group.

\begin{prop}\label{naturalitie}
If $f \in \Diff(M)$ and $\sigma$ is a spin structure on $M$, then $f^{*} Q^{\sigma} = Q^{f^* \sigma}$.
\end{prop}

For quadratic forms on a symplectic $\bF_2$-vector space, the \emph{Arf-invariant} is defined. A \emph{symplectic form} on an $\bF_2$-vector space $V$ is a nondegenerate symmetric bilinear form $(x,y)\mapsto x \cdot y \in \bF_2$ such that $x \cdot x =0$ for all $x \in V$. The vector space $V$, together with the symplectic form, is called a \emph{symplectic vector space}. A \emph{symplectic basis} of $V$ is a basis $(a_1, \ldots a_g , b_1, \ldots b_g)$ such that $a_i \cdot a_j = b_i \cdot b_j =0$ for all $i,j$ and $a_i \cdot b_j = \delta_{ij}$. Any symplectic vector space $V$ has a symplectic basis and it has even dimension.

\begin{defn}
Let $V$ be a symplectic $\bF_2$-vector space and $q: V \to \mb{F}_2$ be a quadratic form (i.e. $q(a+b) = q(a)+ q(b) + a \cdot b$). Let $a_1, \ldots, a_g, b_1, \ldots, b_g$ be a symplectic basis for $V$. Then the \emph{Arf-invariant} is defined to be $ \Arf(q):= \sum_{i=1}^{g} q(a_i) q(b_i) \in \bF_2$. In multiplicative notation, we set $\arf(q) := (-1)^{\Arf(q)} \in\{ \pm 1\}$.
\end{defn}

\begin{prop}\label{addition}
If $V= V_1 \oplus V_2$ is an orthogonal decomposition of a symplectic vector space and if symplectic bases of both subspaces are chosen, then $\Arf (q|_{V_1}) + \Arf (q|_{V_2})= \Arf(q)$.
\end{prop}

The proof is trivial. 
\begin{cor}\label{39}
\begin{enumerate}
\item $\arf(q)= 2^{\dim V /2 } \sum_{x \in V} (-1)^{q(x)}$.
\item $\Arf(q)$ does not depend on the basis used to define it.
\end{enumerate}
\end{cor}

\textbf{Proof:}
The second statement is an immediate consequence of the first one, because it does not involve the basis.\\
The first formula is proven by induction on the dimension of $V$, using Proposition \ref{addition} and an explicit consideration for $\dim V =2$.\qed\\

\begin{expls}
It is worth to give the examples of quadratic forms on a $2$-dimensional vector space. Let $x,y$ be a symplectic basis. Then a form $Q_0$ with $\Arf(Q_0)=0$ is given by 

$$Q_0(x)= Q_0(y)=Q(0)=0; Q_0(x+y)=1,$$

while another form $Q_1$ with $\Arf(Q_1)=1$ is given by 

$$Q_1(x)= Q_1(y)=Q_1(x+y)=1; Q_1(0)=0.$$

There are two other forms on a two-dimensional space with Arf-invariant $0$, which are both isomorphic to $Q_0$.\\
One can check that the quadratic form $Q_0$ is the quadratic form of the spin structure $S_0$ from \ref{spinstrtours}. The other three quadratic forms belong to the three even spin structures on the torus.
\end{expls}

\begin{cor}
The Arf-invariant does not change under symplectic isomorphisms. In particular, if $\sigma$ is a spin structure and $f \in \Diff(M)$, then $\Arf (Q^{f^* \sigma}) = \Arf (Q^{\sigma})$.
\end{cor}

\begin{cor}\label{nullstelle}
Let $n >0 $. The numbers of zeroes of a quadratic form $q$ on a $2n$- dimensional symplectic vector space equals $2^{n-1}(2^n -1)$ if $\arf(q) =1$ and $2^{n-1}(2^n +1)$ if $\arf(q) =-1$. In particular, any quadratic form has a nontrivial zero $(x\neq 0)$ except in the case that $n=1$ and $\arf(q)=-1$.
\end{cor}

\textbf{Proof:} Immediate from the first statement of Corollary \ref{39}.\qed\\

\begin{prop}\label{agree}
Let $(V_0,q_0)$ and $(V_1, q_1)$ be symplectic vector spaces with quadratic forms. Then $(V_0,q_0) \cong (V_1, q_1)$ if and only if the dimensions and the Arf-invariants agree.
\end{prop}

The proof is also very easy. First, the proof for $\dim V =2$ is clear by inspection. For higher dimensions, decompose $V$ into the direct sum of $2$-dimensional spaces. The proposition then follows from the fact that $(\bF_{2}^{2},Q_1) \oplus (\bF_{2}^{2},Q_1) \cong (\bF_{2}^{2},Q_0) \oplus (\bF_{2}^{2},Q_0)$, which can be checked directly, and from Proposition \ref{addition}.\\

\begin{prop}
Let $V$ be a symplectic vector space of dimension $2g$ and $V^{\prime}$ its dual space. Denote by $n_{g}^{-}$ the number of quadratic forms on $V$ with Arf-invariant $-1$ and by $n_{g}^{+}$ the number of quadratic forms with Arf-invariant $+1$. Then
$n_{g}^{-}= 2^{g-1}(2^g -1)$ and $n_{g}^{+}= 2^{g-1}(2^g +1)$.
\end{prop}

\textbf{Proof:}
Since the difference of two quadratic forms is a linear form, the set of all quadratic forms is a $V^{\prime}$-torsor, hence $n_{g}^{+}+n_{g}^{-}=2^{2g}$. The statement is true for $g=1$ by inspection. Proceed by induction, using the addition formula \ref{addition}.\qed\\

The algebraic theory of the Arf-invariant can be applied to surfaces. Another result from Johnson's paper \cite{John} is 

\begin{prop}\label{mjiop}
The assignment of the quadratic form $Q^{\sigma}$ to a spin structure $\sigma$ on a surface is a bijection.
\end{prop}

\textbf{Proof:} Because the number of spin structures and the number of quadratic forms are both equal to $2^{2g}$, it is sufficient to show that any quadratic form is the quadratic form of a spin structure. First, we show this on the torus. There are two different spin structures on $\bS^1$. Taking cartesian products yields four spin structures on $\bS^1 \times \bS^1$. It is straightforward to check that these four spin structures realize the four different quadratic forms.\\
If $g \geq 2$, then we decompose the surface $M$ as the connected sum of $g$ tori. This gives a direct sum decomposition of $H^1( M; \bF_2)$ into $2$-dimensional summands. Each summand has a quadratic form, which can be realized by a spin structure. They can be glued together to give the desired spin structure.\qed\\

\begin{cor}\label{dfghjk}
Two spin structures $\sigma_0, \sigma_1$ on a surface $M$ are conjugate under the action of the mapping class group if and only if $\Arf (Q^{\sigma_0}) = \Arf (Q^{\sigma_0})$.
\end{cor}

\textbf{Proof:}
The "only if" direction is trivial after all we have said above.\\
For the "if" direction, note that by \ref{agree}, there exists a symplectic isomorphism $\phi$ of $H_1(M;\bF_2)$ with $\phi^* Q^{\sigma_1} = Q^{\sigma_0}$. 
We use the fact that $\Gamma_g \to \Sp_{2g}(\bZ)$ is surjective. This can be seen by studying an explicit system of generators of $\Sp_{2g}(\bZ)$, which turn out to realizable by mapping classes. Thus $\Gamma_g \to \Sp_{2g}(\bF_2)$ is also surjective.\\
Hence there exists an $f\in \Diff(M)$ whose action on cohomology is $\phi$. It follows that $Q^{\sigma_0} = f^{*} Q^{\sigma_1} =   Q^{f^{*} \sigma_1}$. By Proposition \ref{mjiop}, it follows that $\sigma_0 = f^{*} \sigma_1$, as desired.\qed\\

\begin{prop}\label{atiarf}
Let $M$ be a Riemann surface and let $\sigma$ be a spin-structure. Then $\At(M, \sigma)=\Arf(Q^{\sigma})$.
\end{prop}

\textbf{Proof:} 
First one proves that $\Arf(Q^{\sigma})$ is spin-bordism invariant. The method is the following. If $(M_0, \sigma_0)$ and $(M_1, \sigma_1)$ are spin-bordant, then $M_1$ is obtained from $M_0$ by a sequence of surgeries which respect the spin structures. Using the algebraic properties of the Arf-invariant, it follows that the Arf-invariant is preserved under all spin-surgeries.\\
Next, one proves that a spin surface with Arf-invariant zero is null-bordant. Using \ref{nullstelle}, one sees that one can decrease the genus by one with a spin-surgery, except in the case that $(g, \Arf)=(1,-1)$ and in the case $(g, \Arf)=(0,1)$. The latter is nullbordant, while the other one is not. So, any spin surface with odd Arf-invariant is bordant to the torus with the odd spin structure.\\
Thus, the Atiyah-invariant as well as the Arf-invariant define homomorphisms $\Omega^{\Spin}_{2} \to \bZ/2$. The Arf-invariant is an isomorphism as we have seen before. The proposition follows, because the Atiyah-invariant is nontrivial by the discussion in \ref{spinstrtours}.\qed\\

\begin{cor}\label{spinbrod}
The second spin bordism group $\Omega_{2}^{\Spin}$ has precisely $2$ elements. \qed
\end{cor}

\begin{prop}\label{stiefel}
Any surface bundle on $\bS^1$ admits a spin structure.
\end{prop}

In view of Lemma \ref{monodromya}, this proposition has an immediate corollary.

\begin{cor}\label{Stiefel}
Any diffeomorphism of $M$ fixes at least one spin structure.
\end{cor}

The corollary was proven by Atiyah (\cite {At}) with algebraic methods. We present an easy topological proof.\\

\textbf{Proof of Proposition \ref{stiefel}:}
Let $p:M \to \bS^1$ be an oriented surface bundle. The total space $M$ is a smooth oriented $3$-manifold. A well-known theorem (due to Stiefel) says that $M$ is parallelizable.
Now we see that $TM \cong M \times \bR^3  \cong T_v M \oplus p^* T \bS^1 \cong T_v M \oplus( M \times \bR)$. Thus the Stiefel-Whitney classes of $T_v M$ are trivial. In particular, the surface bundle $p$ is spin.\qed\\

\subsection{The spin mapping class group}\label{spinmcg}

We now consider the general theory of smooth spin $M$-bundles in the case of (connected, closed, oriented) surfaces of genus $g \geq 2$. The smaller genera will be treated below in the sections \ref{spherebundles} and \ref{torusbundles}. Recall the diagram \ref{rewq} of homotopy fiber sequences 

\begin{equation}
\xymatrix{ 
\bR \bP^{\infty} \ar[r] \ar[d]& B \Sspin(M) \simeq H^1(M; \bF_2) \times \bR \bP^{\infty} \ar[d]\\
B \Spiff(M; \sigma) \ar[d] \ar[r] & B \Spiff(M) \ar[d]\\
B \Diff(M; \sigma) \ar[r] & B \Diff(M).\\
}
\end{equation}

\begin{prop}\label{haupt}
Let $M$ be an oriented surface of genus $g \geq 2$.
\begin{enumerate}
\item $B \Spiff(M)$ has two connected components $B \Spiff(M)^{+}$ and $B \Spiff(M)^{-}$ belonging to the different values of the Atiyah-invariant.
\item Both components are aspherical if $g \geq2$.
\item The fundamental group $\pi_1 (B \Spiff(M)^{\pm})$ is a central $\bZ/2$-extension of $\Gamma(M, \sigma)$, where $\sigma$ is a spin structure with Atiyah-invariant $\pm$.
\end{enumerate}
\end{prop}

\textbf{Proof:} The diagram \ref{ewq} becomes 

\begin{equation}
\xymatrix{ 
0 \ar[r] \ar[d]& 0 \ar[d] \\
\bZ/2 \ar[r]^{\cong} \ar[d]& \bZ/2 \ar[d]\\
\pi_1 (B \Spiff(M; \sigma)) \ar[d] \ar[r]^{\cong} & \pi_1( B \Spiff(M)) \ar[d]\\
\Gamma(M; \sigma) \ar[r] \ar[d] & \Gamma(M) \ar[d] \\
\pt \ar[r] \ar[d]& \pi_0(B \Sspin(M)) \ar[d]\\
\pt \ar[d] \ar[r] & \pi_0(B \Spiff(M))=\pi_0 (B  \Sspin (M)) / \Gamma(M) \ar[d]\\
\pt \ar[r] & \pt\\
}
\end{equation}

because $\Diff_0(M)$ is contractible. All other homotopy groups are zero. By Corollary \ref{dfghjk}, $\pi_0 (B  \Sspin (M)) / \Gamma(M)$ has exactly two elements, which are distinguished by their Atiyah-invariants.
This is enough to show the propostion.\qed\\

\begin{defn}
Let $g \geq 0$ and let $\epsilon \in \bZ/2$. Let $M$ be a surface of genus $g$ and let $\sigma$ be a spin structure on $M$ of Atiyah-invariant $\epsilon$. Set $\Gamma_{g}^{\epsilon} := \Gamma(M, \epsilon)$. The \emph{spin mapping class groups} $\hat{\Gamma}_{g}^{\epsilon}$ is the fundamental group of $B \Spiff(M)^{\epsilon}$.
\end{defn}

\begin{rem}
These groups were defined in other terms by Gregor Masbaum in \cite{Mas}. The definition given is more conceptual.
\end{rem}

The extension $\bZ/2 \to \hat{\Gamma}_{g}^{\epsilon} \to \Gamma_{g}^{\epsilon}$ is always nontrivial. This is Theorem \ref{nontriviall} in the next section.

\pagebreak

\section{Some examples}\label{someexamples}

\subsection{Finite group actions and the nontriviality of the obstruction class}\label{spinexpls}

Let $G$ be a finite group and let $M$ be a closed oriented surface of genus $g$. We will study actions of $G$ on $M$ by orientation-preserving diffeomorphisms. Such an action induces a surface bundle

$$E(G;M):=EG \times_{G} M \to BG.$$

We deal with the question under which condition the surface bundle $E(G;M)$ admits a spin structure. The answer is:

\begin{thm}\label{exisspin}
Let $G$ be a finite group which acts \emph{faithfully} on a closed surface $M$. Then the induced surface bundle $E(G;M)$ is spin if and only if all $2$-Sylow-subgroups act freely on $M$.
\end{thm}

Since all $2$-Sylow-subgroups of a finite group are conjugate, it suffices to check the criterion for one of them.\\
Before I prove Theorem \ref{exisspin}, I will give a few preparatory remarks. First of all, if $G$ acts on $M$ by diffeomorphisms and if we choose a Riemannian metric $h$ on $M$, then we can average $h$ over $G$:

$$h^G := \frac{1}{\sharp G} \sum_{g \in G} g^{*} h.$$

The new metric $h^G$ is $G$-invariant. A metric defines an almost-complex structure on $M$, which is the same as a complex structure. So we can assume without loss of generality that $M$ is a Riemann surface and $G$ a group of biholomorphic automorphisms.\\
Another important fact is that for any $x \in M$, the isotropy subgroup $G_x$, is \emph{cyclic}. This can be seen as follows. Choose an invariant Riemann metric. It is a well-known fact that if an isometry of a connected Riemannian manifold fixes a point and if the differential at this point is the identity, then the isometry must be the identity. Thus the fixed-point representation $G_x \to \SO(2)$ at $x$ is faithful and $G_x$ is cyclic.\\
Now we can prove Theorem \ref{exisspin}. It is a combination of Proposition \ref{propeins} and Proposition \ref{propzwei} below.

\begin{prop}\label{propeins}
Let $G$ be a finite group acting on a closed oriented surface $M$. Let $H \subset G$ be
a $2$-Sylow-subgroup and $j:H \to G$ the inclusion map. Then $E(G;M)$ is spin if and only if $E(H;M)= (Bj)^*E(G;M)$ is spin.
\end{prop}

\textbf{Proof:}
It is clear that $E(H;M)$ is spin if $E(G;M)$ is. To prove the converse, observe that if we choose $BH:= EG/H$, then $Bj: BH \to BG$ is a covering of odd finite degree. Furthermore, we have a covering of total spaces $p: E(H;M) \to E(G;M)$ of the same degree over $Bj$, which is nothing else than $(EG \times M )/H \to (EG \times M )/G$. Note that $p^{*} T_v E(G;M) \cong T_v E(H;M)$. Thus our proposition follows immediately from

\begin{claim}
Let $p:X \to Y$ be a covering map of finite odd degree and let $E \to Y$ be an oriented vector bundle
such that $p^{*} E$ is spin. Then $E$ is spin.
\end{claim}

To prove the claim, we use the transfer $p_{!}:H^2(X, \bF_2) \to H^2(Y, \bF_2)$. Then $p_{!} p^{*}$ is the multiplication with the degree, hence the identity on $H^2(Y; \bF_2)$. Thus $p^{*}$ is injective, which shows that $E$ is spin.\qed\\

\begin{cor}
Let $G$ be a finite group of odd order acting on a surface $M$. Then $E(G;M)$ is spin.\qed\\
\end{cor}

\begin{prop}\label{propzwei}
Let $G$ be a $2$-group acting faithfully on a Riemann surface $M$.
Then $E(G;M)$ is spin if and only if the $G$-action on $M$ is free.
\end{prop}

For the proof of Proposition \ref{propzwei}, I will need two lemmata. If $\rho:G \to \bC^{\times}$ is a representation,
let $\bC_{\rho}$ denotes the complex numbers with the given $G$-representation.

\begin{lem}\label{hzweig}
Let $G$ be a finite group. Then the map $\Hom(G, \bC^{\times}) \to H^2(G, \bZ)$, $ \rho \mapsto c_1(EG  \times_{G} \bC_{\rho})$
is an isomorphism.
\end{lem}

\textbf{Proof:}
Since $G$ is finite, $H^s(BG, \bC)=0$ if $s>0$ and the coefficient sequence $\bZ \to \bC \to \bC^{\times}$ yields an isomorphism $H^{1}(BG, \bC^{\times}) \to H^2(BG, \bZ)$. Further, there is an isomorphism $\Hom(G, \bC^{\times}) \cong H^1(BG, \bC^{\times})$. It is easy to see that the composition of both isomorphisms is the Chern class homomorphism as above.\qed\\

\begin{lem}\label{zweihochn}
Let $\rho: C \to \bC^{\times}$ be a faithful representation of a (necessarily cyclic) group $C$ of order $2^k >1$. Then the second Stiefel-Whitney class of $\rho$ is nonzero.
\end{lem}

\textbf{Proof:}
We must show that $EC \times_{C} \bC_{\rho}$ is not a spin vector bundle.\\
Due to Lemma \ref{hzweig}, it suffices to show that the representation $\rho$ does not have
a square root. This is immediate: $C$ contains a nontrivial element $h$ of order $2$. If $\rho$ had a square root $\sigma$,
then $\rho(h)=(\sigma(h))^2= 1$ and $\rho$ is not faithful.\qed\\

\textbf{Proof of Proposition \ref{propzwei}:}
Let $G$ act freely on $M$. Then there exists an unbranched covering map
$f: E(G;M) \to BG \times M/G$ over $BG$. In addition, $f^{*}(BG \times T(M/G)) \cong T_v E(G;M)$.
Because $M/G$ is a spin manifold, $E(G;M)$ is spin.\\
Conversely, if the action is not free, then there is an action of a nontrivial cyclic subgroup $C \subset G$
which is not free. We prove that $E(C;M)$ is not spin if $C$ is a cyclic subgroup of order $2^k$ for some $k \in \bN$ which does not act freely.
By passing to a smaller nontrivial subgroup, we can assume that a generator of $C$ has a fixed point $x \in M$.\\
The fixed point $x$ yields a section $s:BC \to E(C;M)$ of the surface bundle.
The pullback $s^{*}T_v E(C;M)$ is isomorphic to the complex line bundle $EC \times_{C} \bC_{\rho}$, where $\rho$ the fixed point representation at $x$. The fixed point representation must be faithful.
By Lemma \ref{zweihochn}, $s^{*}T_v E(C;M)$ is not a spin vector bundle and $E(C;M)$ cannot be spin.\qed\\

As an application of Theorem \ref{exisspin}, we can prove that the obstruction cohomology class $c \in H^2(B \Diff(M; \sigma); \bF_2)$ from \ref{obstructionclasas} is nonzero if $M$ is an oriented closed surface and $\sigma$ a spin structure on $M$. Equivalently, the central extensions

$$\bZ/2 \to \Spiff(M) \to \Diff(M)$$

and (if $g \geq 1)$

$$\bZ/2 \to \hat{\Gamma}_{g}^{\epsilon} \to \Gamma_{g}^{\epsilon}$$

are nontrivial.

\begin{thm}\label{nontriviall}
If $M$ is a closed oriented surface and $\sigma$ is a spin structure on $M$, then the obstruction class $c \in H^2 (B  \Diff(M, \sigma); \bF_2)$ is nontrivial.
\end{thm}

\textbf{Proof:}
We need to show that there exist surface bundles whose monodromy fixes $\sigma$, but which do not admit any spin structure.
We have seen that any automorphism of a Riemann surface fixes a spin structure (\ref{Stiefel}). So any cyclic group action fixes a spin structure. However, we need a sharper statement. Let $f \in \Diff(M)$ an element which acts trivially on $H_1 (M; \bF_2)$ and let $\sigma$ be a spin structure. The formula \ref{naturalitie} shows that $Q^{\sigma} = Q^{f^* \sigma}$. By \ref{mjiop}, it follows that $\sigma = f^* \sigma$. Thus $f$ fixes any spin structure.\\
On any surface $M$, there exist \emph{hyperelliptic involutions}. These are involutions, $f \circ f = \id$ on $M$, which are characterized by the following equivalent properties:

\begin{enumerate}
\item\label{fdsa} $f$ has $2g+2$ fixed points.
\item\label{gfds} $f$ acts as $- \id $ on $H_1 (M; \bZ)$.
\item The quotient $M/\{\id, f\}$ is homeomorphic to $\bC \bP^1$.
\end{enumerate}

Let $f$ be such an involution, determining an action $\bZ/2 \actson M$. Let $\sigma$ be an arbitrary spin structure on $M$. Then $f$ acts trivially on $H_1 (M; \bF_2)$ by \ref{gfds}, thus it fixes $\sigma$. Thus the monodromy of the bundle $E(\bZ/2; M)$ fixes $\sigma$. Nevertheless, $f$ has fixed points and so, by Theorem \ref{exisspin}, $E(\bZ/2; M)$ does not have a spin structure.\qed

\subsection{Sphere bundles}\label{spherebundles}

Recall Smales theorem (Theorem \ref{Smale}) which says that the inclusion $\SO(3) \to \Diff(\bS^2)$ is a homotopy equivalence. In more geometric terms, this means that any bundle of surfaces of genus $g=0$ is isomorphic, as an oriented $\bS^2$-bundle, to the unit sphere bundle of a 3-dimensional Riemannian vector bundle. Alternatively, since $\SO(3)$ is a maximal compact subgroup of $\bP \Sl_2(\bC)$, any surface bundle of genus zero can be viewed as a $\bC \bP^1$-bundle with structural group $\bP \Sl_{2}(\bC)$. It is not true, however, that any $\bC \bP^1$-bundle is the projective bundle of a $2$-dimensional complex vector bundle.

Since $\bS^2$ is simply-connected, there is a unique spin structure $\sigma_0$ on it. This is fixed by any diffeomorphism. The spin diffeomorphism group fits into an exact sequence

\begin{equation}\label{extension}
\bZ/2 \to  \Spiff(\bS^2, \sigma_0) \to \Diff(\bS^2).
\end{equation}

\begin{prop}\label{esszwei}
In the diagram below, all horizontal maps are homotopy equivalences.

$$\xymatrix{
\bZ/2 \ar@{=}[r] \ar[d]& \bZ/2 \ar@{=}[r] \ar[d] & \bZ/2 \ar[d] \\
SU(2) \ar[r] \ar[d]  & \Sl_2 (\bC) \ar[r]\ar[d] & \Spiff(\bS^2, \sigma_0)\ar[d] \\
\SO(3) \ar[r]  & \bP \Sl_2 (\bC) \ar[r]   & \Diff(\bS^2)  .\\
}$$

\end{prop}

\textbf{Proof:} The left half of this diagram does not need explanation, and we only consider the right half. We have already seen that the extension on the right-hand-side is nontrivial (\ref{nontriviall}). However, there is a more direct and explicit way of seeing this.\\
Let us regard $\sigma_0$ as a complex line bundle (a square-root of the cotangent bundle), or better, as a $\bC^{\times}$-principal bundle. Then its total space is $\bC^2 \setminus \{(0,0)\}$, and the bundle map to $\bC \bP^1$ is the defining map (alias the Hopf fibration). The group $\Sl_2 (\bC)$ acts in the usual way on this total space, and this shows that $\Sl_2(\bC)$ is a subgroup of $\Spiff(\bS^2, \sigma_0)$. Certainly, the element $-1$ acts as the element in the kernel of $\Spiff(\bS^2, \sigma_0) \to \Diff (\bS^2)$. This provides the diagram; and the proposition now follows from Smales theorem.\qed\\

\begin{prop}\label{eszwei}
Let $\pi: E \to B$ be an $\bS^2$-bundle and let $V \to B$ be a $3$-dimensional vector bundle whose sphere bundle is $E$. Then the following conditions are equivalent.
\begin{enumerate}
\item $E$ has a spin structure.
\item $V$ is a spin vector bundle.
\item There exists a $2$-dimensional complex vector bundle $U \to B$, such that $E$ is isomorphic to the projective bundle $\bP U \to B$ and such that $c_1 (U) =0$.
\end{enumerate}
Because of 1 and 2, the class $c \in H^2(B \Diff(\bS^2), \bF_2) \cong H^2(B\SO(3); \bF_2) \cong \bF_2$ agrees with $w_2$.
\end{prop}

\textbf{Proof:}
We have seen before that the existence of a spin structure on a $\bC \bP^1$-bundle is equivalent to the existence of a lifting of the structural group from $\bP \Sl_2 (\bC)$ to $\Sl_2(\bC)$. This shows the equivalence of $1$ and $3$. Because the maximal compact subgroup of $\Sl_2(\bC)$ is $\Spin(3) = \bS^3$, the equivalence of these with $2$ also follows immediately.\qed\\

The class of projective bundles of two-dimensional complex vector bundles is slightly larger than the class of spin-$\bS^2$-bundles. There is an obstruction for the lifting of bundles of pro\-jective spaces to vector bundles which lies (universally) in $H^3 (B\bP \Sl_2 (\bC); \bZ) \cong \bZ/2$. It can be identified with $\beta c$, the image of $c \in H^2(B \bP \Sl_2(\bC); \bF_2)=H^2(B \Diff(\bS^2); \bF_2) $ under the Bockstein homomorphism for the short exact coefficient sequence $0 \to \bZ \to \bZ \to \bF_2 \to 0$.\\
Now let $E \to B$ be an oriented surface bundle of genus zero. Let $V\to B$ a 3-dimensional oriented real Riemannian vector bundle whose unit sphere bundle is $E$. We can express the characteristic classes of the surface bundle $\pi: E \to B$ by the characteristic classes of the vector bundle $V$.

\begin{prop}\label{classesgnull}
\begin{enumerate}
\item The symplectic classes of $\pi:E \to B$ are trivial.
\item The odd MMM-classes vanish, $\kappa_{2k+1}=0$.
\item The even MMM-classes are related to the Pontryagin classes of $V$ by the formula: $\kappa_{2k} = 2 p_1(V)^k$.
\end{enumerate}
\end{prop}

\textbf{Proof:}
The first statement is trivial, since the first homology group of a sphere is zero.\\
We prove the other statements for the universal case, that is: $V:= E(\SO(3); \bR^3 ) = E\SO(3) \times_{\SO(3)} \bR^3 \to B \SO(3)$ and $E:= E(\SO(3); \bS^2)$.
The total space of the universal surface bundle on $B\SO(3)$ is $E\SO(3) \times_{\SO(3)} \SO(3)/\SO(2) = E\SO(3)/ \SO(2) \simeq \cpinf$.\\
The unit tangent bundle of $\bS^2$ is diffeomorphic to $\SO(3)$; it is a free transitive $\SO(3)$-space under the action of the isometries of $\bS^2$. Hence the total space of the vertical unit tangent bundle of the universal $\bS^2$-bundle is contractible. The only $\bS^1$-bundles over $\cpinf$ whose total spaces are contractible are the universal $\bS^1$-bundle and its dual, and we can identify the vertical tangent bundle of the universal $\bS^2$-bundle with the universal complex line bundle or its dual. This ambiguity will be resolved at the end of the proof. \\
If $z $ is the standard generator of $H^2(\cpinf; \bZ)$, namely, the first Chern class of the universal line bundle $(E(\bS^1; \bC))$, then the Euler class $e$ of the vertical tangent bundle is $\pm z$.\\
The projection map $q:\cpinf \to B\SO(3)$ is homotopic to the map induced by the inclusion $\SO(2) \to \SO(3)$ of the standard maximal torus in $\SO(3)$. If we take the real standard representation of $\SO(3)$ on $\bR^3$, complexify it and restrict it to the maximal torus, we obtain a sum of the trivial representation, the standard representation of $\bS^1$ and its dual. Let $p_1 \in H^4(B\SO(3); \bZ)$ be the universal first Pontryagin class. Then $q^{*}(p_1 (V))=-c_2(q^{*}V \otimes \bC)$ by the definition of the Pontryagin classes. But $-c_2(q^{*}V \otimes \bC) = -c_2((\gamma_1 \oplus \bR) \otimes \bC)= -c_2 (\gamma_1 \oplus \bar{\gamma_1} \oplus \bC) = -c_1(\gamma_1) c_1(\bar{\gamma_1}) = -z (-z) = z^2$. Thus $q^{*}(p_1 (V))= z^2$.\\
We are ready to compute the MMM-classes: \\
$\kappa_{n+2} = q_{!}((\pm z)^{n+3}) =(\pm)^{n+3} q_{!}(q^{*}(p_1(V) z^{n+1}) = (\pm )^{n+1} p_1(V) q_{!}(z^{n+1} )= p_1(V) \kappa_{n}$.
From the facts that $\kappa_0 =+2 = \chi(\bS^2)$ and $\kappa_1 =0$ (since $H^2(B\SO(3), \bZ)=0$) the claim follows by induction.\qed\\

If the $\bS^2$-bundle is given as a projective bundle $\bP W$ of a complex $2$-dimensional vector bundle $W$, then we can give an expression of its MMM-classes in terms of the Chern classes of $W$.
By Proposition \ref{classesgnull}, it suffices to compute the first Pontryagin class of a certain three-dimensional vector bundle $V$, depending on $W$.

\begin{prop}\label{projcompl}
Let $W \to B$ be a two-dimensional complex vector bundle. Then the MMM-classes of the projective bundle $\bP W$ (it is a surface bundle of genus $0$) are given by the formulae $\kappa_{2k+1}=0$ and $\kappa_{2k}= 2 (c_1(W)^2-4c_2(W) )^k$.
\end{prop}

\textbf{Proof:}
We can assume that the structural group of $W$ is reduced to $U(2)$. Let $Q$ be a principal $U(2)$-bundle for $W$. Then, as an $\bS^2$-bundle, $\bP W$ is the same as $\bS V$, where the $3$-dimensional vector bundle $V$ is obtained by $V:= Q \times_{U(2)} \bR^3$. Here $U(2)$ acts on $\bR^3$ via the map $U(2) \to \bP U(2) \cong \SO(3)$ or, more geometrically by linearly extending the action on $\bS^2$ by M\"obius transformations. For the calculation of $p_1$, we only need the complexification of $V$, which is $Q \times_{U(2)} \bC^{3}_{\rho}$, where $\rho: U(2) \to \Gl_3(\bC)$ is the complexification of the representation above.\\
A character calculation, carried out on the standard maximal torus of $U(2)$, shows that $\rho$ is equivalent to $\Sym^2(D) \otimes (\det (D))^{-1}$, where $D$ is the defining representation of $U(2)$. In other words, we see that $V_{\bC} = \Sym^2(W) \otimes (\det (W))^{-1}$, whence $p_1 (V) = -c_2 (V_{\bC} ) = -c_2 (\Sym^2 (W) \otimes \det(W)^{-1})$.\\
By the splitting principle, we can assume that $W$ is a sum of line bundles, $W= E_0 \oplus E_1$. Then

$$\Sym^2 (W) \otimes \det(W)^{-1} \cong \bC \oplus E_{0} \otimes E_{1}^{*} \oplus E_{1} \otimes E_{0}^{*}.$$

It follows that $p_1 (V)= -c_2 (\bC \oplus E_{0} \otimes E_{1}^{*} \oplus E_{1} \otimes E_{0}^{*}) =  -c_2 ( E_{0} \otimes E_{1}^{*} \oplus E_{1} \otimes E_{0}^{*})= - c_1( E_{0} \otimes E_{1}^{*}) c_1 ( E_{1} \otimes E_{0}^{*}) = c_1 (E_{0} \otimes E_{1}^{*})^2 $. This is the same as $(c_1(W))^2-4c_2(W) $.
\qed\\

We have seen that $\Spiff(\bS^2, \sigma_0) \simeq SU(2)= \bS^3$. But the classifying space $B \bS^3$ is the same as the infinite quaternionic projective space $\bH \bP^{\infty}$, whose integral cohomology is

$$H^{*} (\bH \bP^{\infty}; \bZ) \cong \bZ[ u]; \deg (u) =4.$$

The generator $u$ can be chosen to be the second Chern class of the universal $2$-dimensional vector bundle $ V:= E SU(2) \times_{SU(2)} \bC^2 \to BSU(2)$.

\begin{prop}\label{hpunendlich}
With the choice of $u$ specified above, the MMM-classes of the universal spin-$\bS^2$-bundle on $\bH \bP^{\infty}$ are $\kappa_{2n-1}=0$ and $\kappa_{2n}= (-1)^{n} 2^{2n+1} u^n$.
\end{prop}

\textbf{Proof:}
This is immediate from Proposition \ref{projcompl}, because the universal $\bS^2$-bundle on $\bH \bP^{\infty}$ is the projective bundle of the complex vector bundle $V$.\qed\\

There is a more direct and easier way to see this, at least up to sign: Use the fact that the total space of $\bP V$ is homotopy equivalent to $\bC \bP^{\infty}$ and use the Gysin-sequence of the oriented sphere bundle $\bS^2 \to \bC \bP^{\infty} \to \bH \bP^{\infty} $.\\
Now let us calculate the analytic MMM-classes (see \ref{analytical}) of sphere bundles.

\begin{prop}
Let $E$ be the unit sphere bundle in a 3-dimensional oriented real vector bundle $V$. Then for $n \geq 1$, the analytic MMM-classes are given by $\lambda_{n}(E) = s_n (V \otimes_{\bR} \bC)$.
\end{prop}

\begin{proof}
The only nontrivial summand in the index bundle defining the analytic classes is the vector bundle of holomorphic tangent vector fields to the fibers of $E$. Recall that $H^0(\cpeins, \cO_{T\cpeins})= H^0(\cpeins, \cO(2)) \cong \Sym^2(\bC^2)$, the vector space of homogeneous polynomials of degree $2$ in $2$ variables. The $\bP \Sl_2(\bC)$-action is given by substituting the variables.\\
The restriction of this representation to the maximal compact subgroup $\SO(3) \subset \bP \Sl_2 (\bC)$ is the complexification of the defining representation of $\SO(3)$. If $Q$ is a $\bP \Sl_2(\bC)$-principal bundle for $E$, it follows that $\lambda_n (E) = s_n (Q \times_{\bP \Sl_2(\bC)} \Sym^2(\bC^2)) = s_n (V \otimes_{\bR} \bC)$.
\end{proof}

Because $H^2(B\SO(3), \bZ)=0$, we have $\lambda_1 =0$. Note that the class $c_3(V \otimes \bC)$ is $2$-torsion, but does not vanish. We can compute $\lambda_2 = -2 c_2 = \kappa_2$. For the determination of the values of the higher $\lambda_n$, we use the Newton formula (\cite{MSt}, p.195). In our case it states that $\lambda_n + c_2 \lambda_{n-2} -c_3 \lambda_{n-3}=0$. Inductively, one can see that $\lambda_n =  \kappa_n$ for $n=0,1,2,4$ and that the difference is annihilated by $2$. However, the difference vanishes only if $c_3(V \otimes \bC)=0$.
The argument in the proof of Proposition \ref{projcompl} shows that this is the case for example if the surface bundle is the projective bundle of a $2$-dimensional complex vector bundle.

\subsection{Torus bundles}\label{torusbundles}

For torus bundles, things behave completely different. Let $T = \bC / (\bZ \oplus i \bZ)$ be the standard torus. The starting point for our considerations is the short exact sequence of topological groups

$$\Diff_0(T) \to \Diff(T) \to \Sl_2(\bZ),\label{homotp}$$

where $\Diff_0(T)$ is the group of diffeomorphisms homotopic to the identity and the homomorphism $\Diff(T) \to \Sl_2(\bZ)$ is induced by the action on the first homology group of $T$. There is a splitting $\Sl_2(\bZ) \to \Diff(T)$ given by the linear action of $\Sl_2 (\bZ)$ on $\bC^2$. Thus there is semidirect product decomposition $\Diff(T) \cong \Diff_0(T) \rtimes \Sl_2(\bZ)$, where $\Sl_2(\bZ)$ acts on $\Diff_0(T)$ by conjugation. Since $T$ is a connected Lie group, there is a homomorphism $j:T \to \Diff_0 (T)$, which is a homotopy equivalence by the work of Earle and Eells (\ref{Smale}). Moreover, $j$ is $\Sl_2 (\bZ)$-equivariant. Hence, the short exact sequence above is homotopically equivalent to the right-split exact sequence

$$T \to T \rtimes \Sl_2 (\bZ) \to \Sl_2(\bZ).$$

It follows that the homotopy groups of $B \Diff (T)$ are given by

\begin{equation}
\pi_{r}( B \Diff(T)) \cong
\begin{cases}
\Sl_2(\bZ), &  r=1,\\
\bZ^2 ,  &  r = 2,\\
0, & \text{ otherwise.}\\
\end{cases}
\end{equation}

Although $\Diff(T) \simeq \Sl_2(\bZ) \times T$, the classifying space is \emph{not} a product. Instead we have:

\begin{prop}
The action of $\pi_1 (B\Diff(T)) = \Sl_2 (\bZ) $ on the second homotopy group $\bZ^2$ is given by the natural action of $\Sl_2(\bZ)$.\\
More precisely, $B\Diff(T) \simeq E\Sl_2(\bZ) \times_{\Sl_2(\bZ)} B T$ and the total space of the universal surface bundle is equivalent to $E\Sl_2(\bZ) \times_{\Sl_2(\bZ)} E T$.
\end{prop}

\textbf{Proof:}
It follows from the general construction of classifying spaces for topological groups $G$, that $EG$ is not only a right $G$-space, but also a left $\Aut(G)$-space, both actions commute and $EG \to BG$ is $\Aut (G)$-equivariant. See \cite{Seg68} or any other source on the simplicial construction of $EG$ and $BG$ for details.\\
Because $\Sl_2(\bZ) = \Aut(T)$, there is a torus bundle

$$E\Sl_2(\bZ) \times_{\Sl_2(\bZ)} E T \to E\Sl_2(\bZ) \times_{\Sl_2(\bZ)} B T.$$

I claim that this is a universal one, or, equivalently, that its classifying map

$$E\Sl_2(\bZ) \times_{\Sl_2(\bZ)} BT  \to B\Diff(T)$$

is a homotopy equivalence. This is clear, since if $\Gamma$ is a discrete group acting on another group $A$, then the classifying space of the semidirect product $G = A \rtimes \Gamma$ is given by $E \Gamma \times_{\Gamma} BA$.\qed

\begin{lem}\label{verticaltorus}
If $\pi:E \to B$ is a bundle of complex tori, then $\pi^{*}V_1 (E) \cong \Lambda_v (E)$ (the notations are explained in section \ref{Hodge}).
\end{lem}

\textbf{Proof:}
$\pi^{*}V_1 (E)$ is the complex vector bundle on $E$ whose fiber at a point $e \in E$ consists of all holomorphic $1$-forms $\omega $ on the Riemann surface $E_{\pi(e)}$. Since the vector space of holomorphic $1$-forms on a torus is $1$-dimensional and since a holomorphic 1-form on a torus does not have a zero, the evaluation map $\pi^{*}V_1 (E) \to \Lambda_v E$; $(\omega,e)\mapsto \omega(e)$ is an isomorphism.\qed\\

\begin{cor}
For torus bundles, all topological MMM-classes $\kappa_n$ vanish (in integral cohomology).
\end{cor}

\textbf{Proof:}
Lemma \ref{verticaltorus} implies that $c_{1} (T_v E)= \pi^{*} c_{1}(V_1)$. Thus\\
$\kappa_n (\pi)= \pi_{!}( (c_{1} (T_v E))^{n+1}) = \pi_{!} (\pi^{*}(c_1 (V_1))^n c_1(T_v E) = (c_1 (V_1))^n \pi_{!}(c_1(T_v E))$. \\
But $\pi_{!}(c_1(T_v E)) = \kappa_0 =0$ (which is the Euler number of the fibers).\qed\\

For the computation of the symplectic classes and of the analytic MMM-classes, we need a description of the Hodge bundles of the universal torus bundle.\\
For tori, one has $\dim V_n=1$ for all $n$, and the the natural homomorphisms

\begin{equation}\label{tutut}
V_k \otimes V_l \to V_{k+l}
\end{equation}

are isomorphisms of $1$-dimensional vector bundles. Thus it is sufficient to work out the case of the bundle $V_1$.\\
Let $\phi: B\Diff(T) \to BU(1)$ be the classifying map for the first Hodge bundle. Because it is induced by the action of $\Diff(T)$ on $H^1 (T; \bZ)$, we conclude that $\phi$ factors through the map $B\Sl_2(\bZ) \to B \Sl_2 ( \bR) \simeq BU(1)$.
Thus it suffices to consider torus bundles with a section, i.e. with structural group $\Sl_2(\bZ)$. The cohomology of $\Sl_2(\bZ)$ is well-known. In the fanciest language, the result is

\begin{prop}\label{bousfield}
The abelianization of $\Sl_2(\bZ)$ is a map $\Sl_2(\bZ) \to \bZ/12$. The induced map $B\Sl_2(\bZ) \to B (\bZ /12)$ is a $H \bZ$-localization in the sense of homotopy theory (see \cite{Bous}). The compositions $\Sl_2 (\bZ) \to \bZ/12 \to \Sl_2 (\bR)$ (take the standard representation of $\bZ/12$ on $\bR^2$) and $\Sl_2 (\bZ) \to Sl_2 (\bR)$ give homotopic maps on the level of classifying spaces.
\end{prop}

The proof relies on the well-known fact that $\Sl_2 (\bZ)$ is isomorphic to the amalgamated product $\bZ/4 \ast_{\bZ/2} \bZ/6$ and the Mayer-Vietoris sequence for group cohomology.\\
As a consequence, $H^{*}(B\Sl_2 (\bZ); \bZ ) \cong H^{*}(B \bZ /12 ; \bZ ) \cong \bZ[u] / (12 u)$, $\deg (u) = 2$. The generator $u$ can be chosen as the Euler class of the representation $\Sl_2 (\bZ) \to \Gl_{2}^{+}  (\bR)$.\\
Let us compute the first symplectic class $\gamma_1 \in H^2 (B\Sl_2(\bZ); \bZ)$. The bundle of first cohomology groups $E\Sl_2(\bZ) \times_{\Sl_2(\bZ)} H^1 (T, \bR) = E\Sl_2(\bZ) \times_{\Sl_2(\bZ)} \bR^2$ has Euler class $u$, by definition. The choice a complex structure on the fibers of the universal torus bundle gives rise to a complex structure on this real oriented vector bundle. With this complex structure, the resulting complex line bundle is isomorphic to the first Hodge bundle. Thus

\begin{cor}\label{nichtnull}
We have the equality $\gamma_{1} = u \in H^2(B\Sl_{2}(\bZ), \bZ)$.
\end{cor}

Now we are able to compute the analytic MMM-classes (see \ref{analytical}) of torus bundles. The result is

\begin{prop}
For a torus bundle, we have $\lambda_n = (-1)^n (1-2^n) u^n$.
\end{prop}

This is a trivial consequence of the definitions, of \ref{tutut}, \ref{bousfield} and \ref{nichtnull}.\\
Next, we turn to the question of the existence of a spin structure on a torus bundle.
The odd spin structure $\sigma$ on a torus is fixed by any diffeomorphism, because it is the unique spin structure with Atiyah invariant $1$. It can be described as the trivial holomorphic line bundle or as the spin structure which corresponds to the Lie group framing on $T$ under the isomorphism $\Omega^{fr}_{2} \cong \Omega^{\Spin}_{2}$. This spin structure is called the \emph{odd spin structure}.\\

\begin{cor}\label{torspin1}
A torus bundle $ \to B$ has a spin structure if the bundle $V_1 (E)$ is a spin vector bundle. Thus, a torus bundle has a spin structure if the first symplectic class $\gamma_1$ is divisible by two.
\end{cor}

\textbf{Proof:} Apply Lemma \ref{verticaltorus}.\qed\\

\begin{prop}\label{gammadurchzwei}
If a torus bundle has an odd spin structure, then $\gamma_1$ is divisible by $2$.
\end{prop}

\textbf{Proof:}
The space of holomorphic sections of an odd spin structure $S$ on a complex torus $T$ is one-dimensional and the map $H^0(T,\cO_S) \otimes H^0(T,\cO_S) \to H^0(T, \Lambda_T)$ is an isomorphism. Therefore the line bundle $V_1$ on the base of the surface bundle is spin.\qed\\

\begin{prop}
For the odd spin structure $\sigma$ on $T$, we can identify the classes $c$ and $ \gamma_1 \pmod{2} \in H^2 (B \Diff(T; \sigma); \bF_2)$.
\end{prop}

\textbf{Proof:}
This is not completely obvious from Propositions \ref{torspin1} and \ref{gammadurchzwei}, because one needs to guarantee that the spin structure provided by \ref{torspin1} is odd. This follows from the computation of $H^2 (B \Diff(T); \bF_2) \cong H^2 (B \Sl_2 (\bZ); \bF_2) \cong \bZ/2$. The universal class $\gamma_1$ is nonzero (see \ref{nichtnull}) and the class $c$ is also nonzero (see \ref{nontriviall}). \qed\\

\begin{cor}\label{extorspin}
A torus bundle has an odd spin structure if and only if $\gamma_1 =0 \pmod{2}$. \qed
\end{cor}

The case of \emph{even} spin structures seems to be more difficult.

\pagebreak

\section{The stable homotopy theory of Riemann surfaces}\label{thsthomoth}

\subsection{The Becker-Gottlieb-Transfer}\label{beckertrans}

In this section, $\pi:E \to B$ will always denote a smooth bundle of $d$-dimensional closed mani\-folds. The fiber is denoted by $M$; and we denote by $Q$ the associated $\Diff(M)$-principal bundle. Thus $E = Q \times_{\Diff(M) } M$ and $T_v E = Q \times_{\Diff(M) } TM$. In this section, we make no assumptions on the orientability of $M$; this will enter the picture later. Our first goal is to recapitulate the definition of the Becker-Gottlieb transfer (see \cite{BG}). Becker and Gottlieb gave the definition only in the special case of a compact structural group; this restriction is unnecessary.
We do not assume that $B$ is a manifold, we merely want to assume that $B$ is not too large. More prescisely:

\begin{assumption}\label{embed}
There exists an exhaustion\footnote{We cannot assume that all $B^{(i)}$ are nonempty.} $B^{(0)} \subset B^{(1)} \subset \ldots $ of $B$ and there exist continuous maps $j_n: E^{(n)} := \pi^{-1} (B^{(n)}) \to \bR^n$, such that the following conditions are satisfied.
\begin{itemize}
\item For any $x \in B^{(n)}$, $j_{n}|_{\pi^{-1}(x)}$ is a smooth embedding and all derivatives of $j_n$ in fiber-direction are continuous.
\item The maps $j_n$ are compatible, i.e. $j_{n+1}|_{E^{(n)}} = j_n$.
\end{itemize}
\end{assumption}

\begin{lem}
If $B$ is a countable CW-complex, then \ref{embed} is satisfied.
\end{lem}

\textbf{Proof:} Consider the space $\Emb (M, \bR^n)$ of smooth embeddings, endowed with the weak Whitney $C^{\infty}$-topology. It has a continuous action of the diffeomorphism group $\Diff(M)$. By Whitney's embedding theorem, $\Emb (M, \bR^n)$ is $r$-connected if $2d+2+r \leq n$. Thus there exists a section in this bundle over the space $B^{(n-2d-2)}$, which can be chosen to extend a given section over smaller skeleta.
On the other hand, a section over $B^{(n)}$ is nothing else than a map $j_n$ whose existence is asserted in the lemma.\qed\\

\begin{rem}\label{kriegl}
The space $\Emb(M; \bR^{\infty}):= \colim_{n \to \infty} \Emb(M; \bR^n)$ is weakly contractible by Whitney's embedding theorem (i.e.: all homotopy groups are trivial) and it carries a free action of $\Diff(M)$. \\
One can show that the quotient map $\Emb(M; \bR^{\infty}) \to \Emb(M; \bR^{\infty})/\Diff(M)$ is a $\Diff(M)$-principal bundle. Thus $\Emb(M; \bR^{\infty})$ is a model for $E \Diff(M)$. It follows that the universal smooth $M$-bundle satisfies Assumption \ref{embed}.\\
\end{rem}

With the given exhaustion of $B$ and $E$, we have models for the suspension spectrum of $B$ and $E$. Namely, we set

$$(\widetilde{\Sigma^{\infty}B_+})_n:= \Sigma^n B^{(n)}_+$$

and similarly for $E$. There is a map of spectra $\widetilde{\Sigma^{\infty}B_+} \to \Sigma^{\infty} B_+$, which is a homotopy equivalence in the sense of \cite{Adams}. Consequently, our notation will not distinguish between these two spectra.
Next, we define the \emph{Thom spectrum} of the vertical normal bundle. The embedding $j_n: E^{(n)} \to \bR^n$ realizes the vertical tangent bundle $T_v E^{(n)}$ as a $d$-dimensional subbundle of $E^{(n)}  \times \bR^n$. We now define the \emph{stable vertical normal bundle} $-T_v E$ of $\pi$. It is not an honest vector bundle, but only a stable vector bundle. 

\begin{defn}
Let $X=(X_0 \subset X_1 \subset X_2 \subset \ldots )$ be an infinite sequence of topological spaces. Then a \emph{stable real vector bundle} of dimension $r$ on $X$ is a sequence of real vector bundles $V_n \to X_n$ of dimension $n+r$, together with a choice of isomorphisms $V_n \oplus \bR \cong V_{n+1}|_{X_n}$.
\end{defn}

\begin{defn}\label{defthomspec}
Let $X$ be a sequence of spaces as above and let $V$ be a stable vector bundle. Then the \emph{Thom spectrum} $\bTh(V)$ is the spectrum whose $n$th space is the Thom space $\Th (V_n)$ of $V_n$ and whose structural maps are given by the inclusions 

$$\Sigma \Th (V_n) \cong \Th (V_n \oplus \bR) \cong  \Th (V_{n+1}|_{X_{n+1}  }) \to \Th (V_{n+1}).$$

\end{defn}

We set $X_n := E^{(n)}$ and 

$$(-T_v E)_n := (T_v E^{(n)})^{\bot} \subset E^{(n)} \times \bR^n.$$

This is an $(n-d)$-dimensional vector bundle; and due to the compatibility condition in Assumption \ref{embed}, there are obvious bundle isomorphisms $\epsilon_n : (-T_v E)_n \oplus \bR \cong (-T_v E)_{n+1}|_{E^{(n)}}$. Thus the vector bundles $(-T_v E)_n $ and the isomorphisms define a stable vector bundle $-T_v E$ of dimension $-d$.\\
The inclusion of $(-T_v E)_n$ into $E^{(n)} \times \bR^n$ as a subbundle yields a map of spectra

$$\inc: \bTh (-T_vE) \to \Sigma^{\infty} E_+.$$

Now we choose tubular neighborhoods of the submanifold $j_n (E_x)$ for all $x \in B$. This can be done in a continuous way, since the space of all tubular neighborhoods is a contractible space (this follows from \cite{Hirsch}, Theorem 5.3 on p. 112). The result is an open embedding of $(-T_v E)_n$ into $B^{(n)} \times \bR^n$.
The Pontryagin-Thom collapse construction now gives a map 

$$\Sigma^n B^{(n)}_+ \to \bTh (-T_v E)_{n},$$

which yields a map of spectra

$$\prt:\Sigma^{\infty} B_+ \to \bTh (-T_vE),$$

the \emph{pretransfer}.

\begin{defn}
The \emph{Becker-Gottlieb-transfer} is the map of spectra

$$\tr= \tr_{\pi}:=\inc \circ \prt: \Sigma^{\infty} B_+ \to \Sigma^{\infty} E_+ .$$

\end{defn}

It is not difficult to show that different choices of the exhaustion and of the embeddings in Assumption \ref{embed} leads to homotopy equivalent Thom spectra and to homotopic transfers.

\subsection{The pushforward}\label{pushforward}

In this section, we let $A$ be a multiplicative cohomology theory alias a ring spectrum. 
Let $u \in \tilde{A}^r (\bS^r)$ be the $r$th suspension of the class $1 \in \tilde{A}^0 (\bS^0)$. Let $V \to X$ be a vector bundle of rank $r$. For $x \in X$, we can identify $V_x$ with $\bR^r$. This gives a map of Thom spaces $j_x : \bS^r \to \Th (V)$. An \emph{$A$-orientation} or an \emph{$A$-Thom class} of $V$ is a class $t \in \tilde{A}^{r} (\Th (V))$, such that for any point $x \in X$, $j_{x}^{*} t = \pm u$ (see \cite{Rud}, Chapter V, for information about this).
An $A$-oriented vector bundle $V \to X$ has an $A$-theory \emph{Euler class}: Let $s: X  \to \Th (V)$ be the zero section and set 

$$\chi_A(E) := s^* t \in A^d (X).$$

The notion of an $A$-orientation can be generalized to stable vector bundles, but this requires another notion which is important in the present context, the \emph{spectrum cohomology}. Let $\gZ = (Z_n, z_n)_{n \in \bN}$ be a spectrum. Then the structural maps $z_n : \Sigma Z_n \to Z_{n+1}$ induce maps $A^{r+n+1} (Z_{n+1} ) \to A^{n+d+1} (\Sigma Z_n)$, which can be composed with the suspension isomorphism $A^{n+d+1} (\Sigma Z_n)   \cong A^{n+d} ( Z_n)$ to give maps $A^{r+n+1} (Z_{n+1} )\to A^{r+n} (Z_{n} )$. The \emph{spectrum cohomology} of $\gZ$ is the inverse limit of this inverse system. For suspension spectra $\Sigma^{\infty} X_+$, the spectrum cohomology agrees with the cohomology of the space $X$, and we will not distiguish between them.\\
There are maps (see \cite{Rud}, p. 70 ff)

\begin{equation}\label{suspenss}
\sigma: A^r (\gZ) \to A^r (\Omega^{\infty} \gZ).
\end{equation}

If $X$ is a space, $\gZ$ a spectrum and $\phi: \Sigma^{\infty} X_+ \to \gZ$ a map, then denote by $\phi^{\flat}: X \to \Omega^{\infty} \gZ$ the adjoint of $\phi$. It is trivial to check from the definitions that the diagam

\begin{equation}\label{specko}
\xymatrix{
A^{*} (\Sigma^{\infty} B_+) \ar@{=}[d]  & A^{*}(\gZ) \ar[l]^-{\phi} \ar[d]^{\sigma} \\
A^{*} (B) & A^{*} (\Omega^{\infty} \gZ) \ar[l]^-{\phi^{\flat}}\\
}
\end{equation}

is commutative.\\
We say that an $A$-orientation of the $d$-dimensional stable vector bundle $V$ on $X$ is a spectrum cohomology class $t \in \tilde{A}^d (\bTh (V))$, which restricts to an $A$-orientation of the vector bundle $V_n$ for every $n \in \bN$.\\
A stable vector bundle does \emph{not} have an Euler class, because there is no spectrum map $\Sigma^{\infty} X_+ \to \bTh (V)$ induced by the zero-section.\\
$A$-orientations give rise to \emph{Thom isomorphisms} $A^{*} (\Sigma^{\infty} X_+) \to \tilde{A}^{*+d} (\bTh (V))$.\\
We say that a manifold bundle $ \pi: E \to B$ is $A$-oriented if its vertical tangent bundle $T_v E$ is $A$-oriented. By \cite{Rud}, Proposition 1.10, p. 309, this is equivalent to the existence of an $A$-orientation of the stable vertical normal bundle of $E$.\\
Let $\pi:E \to B$ be an $A$-oriented bundle of $d$-dimensional smooth manifolds. Then the Thom class defines an isomorphism 

$$\Phi:A^{*} (\Sigma^{\infty} E_+) \to A^{*-d} (\bTh (-T_v E)).$$

\begin{defn}\label{defpusj}
If the vertical tangent bundle of $\pi$ is orientable with respect to the cohomology theory $A$, then the \emph{$A$-theory pushforward} (or umkehr map) is 

$$\pi_{!} := \prt^* \circ \Phi : A^* (E) \to A^{*-d} (B).$$

\end{defn}

The relation between the pushforward map and the Becker-Gottlieb transfer can be stated as follows:

\begin{prop}\label{umkehr} (\cite{BG})
For all $x \in A^{*}(E)$, one has

$$\trf^* (x) = \pi_{!} (\chi_A(T_v E) x).$$

\end{prop}

The proof in \cite{BG} is only formulated for bundles with compact structural groups, but the proof easily generalizes if one uses the definition of the Becker-Gottlieb transfer given in this chapter.\\
The issue of the multiplicative behavior of the transfer is slightly subtle. In general, the cohomology $A^*(E)$ does not admit a cup-product, although $A$ is a ring spectrum. However, if there is a diagonal $E \to E \wedge E$, then there is a multiplication on $A^{*} (E)$. This happens if $E$ is the suspension spectrum of a space; and then the product in spectrum cohomology agrees with the usual cup-product. If $f:E \to F$ is a map of spectra with diagonals, then $f$ induces a ring homomorphism if $f$ commutes with the diagonals. This happens, for example, if $f$ is of the form $f = \Sigma^{\infty} g$ for a map $g$ between spaces. The transfer is not of this form and does not commute with the diagonals. So we cannot expect that $\trf^* (xy) = \trf^* (x) \trf^* (y)$. In particular, we cannot expect that $\trf^* (1) = 1$.\\

\subsection{The universal Becker-Gottlieb transfer}\label{unibecgott}

Now we discuss the universal Becker-Gottlieb transfer. Let $r:\gF \to B \Gl_d (\bR)$ be a Serre fi\-bra\-tion, for example $\gF= B\SO(d)$ or $\gF = B \Spin(d)$. Then we define a spectrum $\bG_{-d}^{\gF}$ as follows.\\
Let $\Gr_d (\bR^{n})$ be the Grassmann manifold of $d$-dimensional (unoriented) subspaces of $\bR^{n}$. Let $\Gr_d (\bR^{n}) \times \bR^{n} \supset U_{d,n} \to \Gr_d (\bR^{n})$ be the tautological $d$-dimensional vector bundle, which is classified by a map $\lambda:\Gr_d (\bR^{n}) \to B \Gl_{d}(\bR)$ and let $U_{n,d}^{\bot}$ be the orthogonal complement of $U_{d,n}$, which is $(n-d)$-dimensional.
Define $B_{n,d}$ to be the pullback in the diagram 

\begin{equation}\label{Fstruk}
\xymatrix{
B_{n,d} \ar[r] \ar[d]^{\theta_{n,d}} & \gF \ar[d]^{r} \\
\Gr_{d}(\bR^n) \ar[r]^-{\lambda} & B\Gl_d (\bR).\\
}
\end{equation} 

It is obvious that there are inclusion maps $B_{n,d} \to B_{n+1,d}$ and that the vector bundles $\theta_{n,d}^{*} U_{d,n}^{\bot} \to B_{n,d}$ form a stable vector bundle $\theta^{*} U_{d}^{\bot}$ of dimension $-d$. Of course, replacing the general linear groups by the orthogonal groups does not change anything essential.

\begin{defn}(\cite{GMTW})\label{defmadtil}
Let $\gF$ be a fibration as above. The spectrum $\bG_{-d}^{\gF}$ is the Thom spectrum of the stable vector bundle $\theta^{*} U_{d}^{\bot}$.
\end{defn}

The maps $\theta_{n,d}^{*} U_{d,n}^{\bot} \to \theta_{n,d}^{*} U_{d,n}^{\bot} \oplus \theta_{n,d}^{*} U_{d,n} = B_{n,d} \times \bR^{n} \to \gF \times \bR^n$ define a map of spectra

\begin{equation}\label{omega}
\omega: \bG_{-d}^{\gF} \to \Sigma^{\infty} \gF_+.
\end{equation}

Now we consider a manifold bundle $\pi: E \to B$ with $d$-dimensional fibers. A \emph{tangential $\gF$-structure} on $E$ is a lift $s$

$$\xymatrix{
   & \gF \ar[d]^{r_{d}} \\
 E \ar[ur]^{s} \ar[r]_-{T_v E} & B\Gl_{d} (\bR)\\
}$$ 

of the classifying map of the vertical tangent bundle of $E$. Examples are ordinary orientations of the vertical tangent bundle and spin structures on the vertical tangent bundle. A tangential $\gF$-structure defines a map of spectra $\eta:\bTh (-T_v E) \to \bG_{-d}^{\gF}$ as follows. The Gau{\ss} map $T_v E: E^{(n)} \to \Gr_d (\bR^{n})$ defined by the embedding $\Gr_d (\bR^{n})$ from Assumption \ref{embed} together with the lift $s$ gives us a map $E^{(n)} \to B_{n,d}$, which is covered by a bundle map $(-T_v E)_n \to \theta_{n,d}^{*} U_{d,n}^{\bot}$. This is a map of stable vector bundles. Thomification yields the desired map of spectra.\\
Moreover, the tangential $F$-structure on $E$ defines a spectrum map $\Sigma^{\infty} E_+ \to \Sigma^{\infty} \gF_{+}$, also denoted by $s$. The diagram 

\begin{equation}\label{defmaddd}
\xymatrix{
\Sigma^{\infty} B_+ \ar[r]^-{\prt} & \bTh (-T_v E) \ar[r]^-{\inc} \ar[d]^{\eta} & \Sigma^{\infty} E_+ \ar[d]^{s} \\
    & \bG_{-d}^{\gF} \ar[r]^{\omega} & \Sigma^{\infty} \gF_{+}\\
}
\end{equation}

is obviously commutative. The adjoint of the map $\alpha^{\sharp}:=\eta \circ \prt:\Sigma^{\infty} B_+ \to \bG_{-d}^{\gF}$ will be denoted by 

$$\alpha=\alpha_{(E, \pi, B)} : B \to \Omega^{\infty} \bG_{-d}^{\gF}.$$

We call this map the \emph{Madsen-Tillmann map} of the manifold bundle with $\gF$-structure $\pi$.\\
Now we assume that $\gF$-bundles are $A^{*}$-oriented, in other words, there exists an $A^{*}$-orientation of the universal bundle $r^{*} E(\Gl_d (\bR), \bR^d)$ on $\gF$. If $\pi: E \to B$ is a $d$-manifold bundle with a tangential $\gF$-structure and if $c \in A^{l}(\gF)$ is a characteristic class, then we can define \emph{generalized MMM-classes} by the formula

$$\pi_{!} (\chi_A (T_v E) \cdot c (T_v E)) \in A^{l-d}(B).$$

The relation between the generalized MMM-classes and the transfer is now easy to explain.

\begin{prop}\label{mummorid}

$$\pi_{!} (\chi_A (T_v E) \cdot c (T_v E))= \eta^{*} \circ \prt^{*} \circ  \omega^{*} (c).$$

\end{prop}

This follows from a trivial diagram chase in the diagram \ref{defmaddd} and from Proposition \ref{umkehr}. The case $d=2$ and $\gF= \SO(2)$ is most important for us. If we set $A = H \bZ$ and $c = z^n= c_{1}^{n}$, we obtain

\begin{equation}\label{gaga}
\kappa_n =(\alpha^{\sharp})^* (\omega^* z^n).
\end{equation}

We emphasize again that $(\alpha^{\sharp})^*  \circ \omega^*$ is \emph{not} multiplicative.

\subsection{Relation to bordism theory}\label{relbordism}

A particularly important case of a fibration $\gF \to B \Gl_d (\bR)$ arises from \emph{stable fibrations} over $B \Gl (\bR)$. By this term, we mean a sequence of spaces $F_d$, fibrations $r_d: F_d \to B \Gl_d (\bR)$ and maps $i_d: F_d \to F_{d+1}$, making the following diagram strictly commutative

\begin{equation}\label{gleichnug}
\xymatrix{
F_d \ar[d]^{r_d} \ar[r]^{i_d} & F_{d+1} \ar[d]\\
B \Gl_d (\bR) \ar[r] & B \Gl_{d+1}(\bR).\\
}
\end{equation}

The examples one should have in mind are $F_d := B \SO(d)$ and $F_d := B \Spin (d)$. 
Although it is inprecise notation, we will write $\bG_{-d}^{\gF} := \bG_{-d}^{F_d}$, $\bG_{-d}^{\SO}:=\bG_{-d}^{B \SO(d)}$ and $\bG_{-d}^{\Spin}:=\bG_{-d}^{B \Spin(d)}$, because there is no danger of confusion.
One can define maps of spectra

$$\iota_d :  \bG_{d}^{\gF} \to \Sigma\bG_{-(d+1)}^{\gF}$$

by the following construction. Adding the span of the last unit vector $e_{n+1} \in \bR^{n+1}$ gives a map $\Gr_{d} (\bR^n) \to \Gr_{d+1} (\bR^{n+1}) $, which is covered by a bundle map $U_{d,n}^{\bot} \to U_{d+1,n+1}^{\bot}$. Using the commutative diagrams \ref{gleichnug} and 

\begin{equation}\label{nichtnutz}
\xymatrix{
\Gr_{d}(\bR^{n}) \ar[d] \ar[r] & B \Gl_{d} (\bR) \ar[d]\\
\Gr_{d+1}(\bR^{n+1})  \ar[r] & B \Gl_{d+1} (\bR),\\
}
\end{equation}

one constructs maps 

$$(\bG_{-d}^{\gF})_n = \Th (\theta^{*}_{n,d} U_{d,n}^{\bot}) \to \Th (\theta^{*}_{n+1,d+1} U_{d+1,n+1}^{\bot}) =  (\bG_{-(d+1)}^{\gF})_{n+1} $$

which are compatible with the structural maps of the spectra and give the requested map.
\\
\\
\begin{prop}\label{cofibseq}(\cite{GMTW})
\begin{enumerate}
\item The homotopy colimit of the sequence $\bG_{-d}^{\gF} \to \Sigma \bG_{-(d+1)}^{\gF} \to  \Sigma^{2}  \bG_{-(d+2)}^{\gF} \ldots $ is the universal Thom spectrum $\Sigma^{-d} M \gF$.
\item There is a cofiber sequence of spectra

$$
\xymatrix{
\Sigma^{-1} \bG_{-(d-1)}^{\gF} \ar[r]^-{\iota_{d}} &  \bG_{-d}^{\gF} \ar[r]^-{\omega} & \Sigma^{\infty} \gF_{d,+},\\
}
$$

where $\omega$ is the map described in \ref{omega}.
\end{enumerate}
\end{prop}

The proof can be found (modulo details) in \cite{GMTW}. It is also instructive to make the composite of all the $\iota$'s explicit. This is a map $\lambda: \bG_{-d}^{\gF} \to \Sigma^{-d}M \gF$. The bundle $\theta_{n,d}^{*} U_{n,d}^{\bot} \to B_{n,d}$ is an $n-d$-dimensional vector bundle which has an $\gF$-structure, that is a lift

$$
\xymatrix{
 &   F_{n-d} \ar[d]^{r_{n-d}}\\
 B_{n,d} \ar[ur] \ar[r] & B \Gl_{n-d} (\bR).\\
 }
 $$

On the level of Thom spaces, we obtain a map 

$$(\bG_{-d}^{\gF})_n:=\Th (\theta_{n,d}^{*} U_{n,d}^{\bot}) \to \Th (r^{*}_{n-d} E( \Gl_{n-d}(\bR); \bR^{n-d})) = (M \gF)_{n-d} = (\Sigma^{-d}M \gF)_n,$$
 
which define the spectrum map $\lambda$.\\
Let us consider the whole construction for the case that the base space $B$ is a point. Let $M$ be a $d$-dimensional manifold with a normal and tangential $\gF$-structure and let $\pi:M \to \pt$ be the constant map. Then we have maps

\begin{equation}\label{thompont}
\xymatrix{
\Sigma^{\infty} B_+ \ar@{=}[r] & \Sigma^{\infty} \bS^0 \ar[r]^-{\alpha^{\sharp}_{(M, \pi, *)}} & \bG_{-d}^{\gF} \ar[r]^-{\lambda} & \Sigma^{-d} M\gF,\\ 
}
\end{equation}

or in other words; an element in $\pi_d (M \gF)$, which is the bordism group  $\Omega^{\gF}_{d}$ of $d$-dimensional manifolds with a normal $\gF$-structure by the Pontryagin-Thom theorem. The following is not surprising at all; and the proof is merely an unwinding of all the definitions.

\begin{prop}\label{bordismclass}
The element in $\pi_d (M \gF) \cong \Omega^{\gF}_{d}$ represented by the map \ref{thompont} corres\-ponds to the $\gF$-bordism class of the manifold $M$ under the Pontryagin-Thom isomorphism.
\end{prop}

Finally, one can consider $\bG_{-2}^{\SO}$ as $\bG_{-2}^{U}$, because $U(1) =SO(2)$. Thus one also has a map $\lambda: \bG_{-2}^{\SO} \to \Sigma^{-2} \gM \gU$.\\
It is clear that all constructions in this and the preceding section are natural with respect to maps of fibrations over $B \Gl (\bR)$.

\subsection{The theorems of Tillmann, Madsen and Weiss}\label{madsenweissthms}

From now on, we let $d=2$. There are two examples for $\gF$ which play a role in this work. The first one is $\gF = BSO(2) = \bC \bP^{\infty}$ where the fibration $\gF \to B\Gl_2 (\bR)$ is the classifying map for the universal complex line bundle, considered as a real vector bundle. We denote the resulting spectrum by $\bG_{-2}^{SO}$. In \cite{MW}, it is denoted by $\bC \bP^{\infty}_{-1}$ for some reasons. 
With this choice of $\gF$, an $\gF$-structure on a surface bundle $E \to B$ is nothing else than an orientation of the vertical tangent bundle.
If $M$ is an oriented closed smooth surface of genus $g$, then smooth $M$-bundles on $B$ are classified by maps $B \to B \Diff^+(M)$, which is homotopy equivalent to $B \Gamma_g$ for $g \geq 2$. Because of the remark \ref{kriegl}, the constructions above can be applied to the universal surface bundle, and we obtain a map 

$$\alpha:B \Diff(M) \to  \Omega^{\infty} \bG_{-2}^{SO},$$

which was first studied in \cite{MT}. The group of components $\pi_0 (\Omega^{\infty} \bG_{-2}^{SO})$ is isomorphic to $\bZ$; and the map $\alpha$ maps $B \Diff(M)$ into the component belonging to $g-1$ (or $1-g$, depending on the choice of the isomorphism), see \cite{GMT}. The Madsen-Weiss theorem is the following

\begin{thm}\label{madsenweiss} (\cite{MW}, \cite{GMTW})
If $2k < g$, then the map $\alpha_*:H_k( B \Diff(F) )\to H_k(\Omega^{\infty}_{g-1} \bG_{-2}^{SO})$ is an isomorphism.
\end{thm}

The range of dimensions in which $\alpha$ is an equivalence is precisely the range of stability in Harer's stability theorem (\cite{Harer}), which is a crucial ingredient for the proof of Theorem \ref{madsenweiss} and the only part of the proof that does not apply to more general situations. I will now explain Harer's theorem.\\
For the proof of both, the Harer stability and the Madsen-Weiss theorem, it is crucial to consider surfaces with boundary. Let $M_{g,n}$ be a connected surface of genus $g$ with $n \geq 0$ boundary components and let $\Diff(M_{g,n}; \partial M_{g,n})$ be the group of all orientation-preserving diffeomorphisms of $M_{g,n}$ which fix the boundary pointwise. If $g \geq 2$ or $n \geq1$, then the unit component of these groups are contractible (\cite{EaEe}, \cite{Iv}) and the groups $\Gamma_{g,n} := \pi_0 (\Diff(M_{g,n}; \partial M_{g,n}))$ of components are the \emph{mapping class groups}. There are stabilization maps 

\begin{enumerate}
\item $\Gamma_{g,n} \to \Gamma_{g,n-1}$ for $n \geq 1$ (glue in a disc in one boundary component and extend diffeomorphisms by the identity);
\item $\Gamma_{g,n} \to \Gamma_{g,n+1}$ for $n \geq 1$ (glue in a pair of pants along one boundary component);
\item $\Gamma_{g,n} \to \Gamma_{g+1, n-1}$ for $n \geq 2$ (glue in a pair of pants along two boundary components).
\end{enumerate}

Harer's stability theorem asserts that all these maps are $H_k ( ; \bZ)$-isomorphisms as long as $2k < g$. 
The composition of the last two stabilization maps is a map $ \Gamma_{g,n} \to \Gamma_{g+1, n}$ (which is obtained by glueing in a torus with two boundary components) and one can define the infinite mapping class group

$$\Gamma_{\infty,n} := \colim_{g \to \infty} \Gamma_{g,n}.$$ 

The map $\alpha$ can be constructed also for surfaces with boundary and this can be done in a way that it extends to a map $B \Gamma_{\infty,n} \to \Omega^{\infty}\bG_{-2}^{\SO}$.\\
It is a well-known theorem that the mapping class groups are perfect if $g \geq 3$, i.e. $H_1 (B \Gamma_{g,n}; \bZ) =0$. It was first proven by Powell (\cite{Pow}). A simpler proof was given by Harer \cite{Harer2}. It also follows from the Madsen-Weiss theorem. Thus one can apply Quillen's plus construction (see \cite{Ads}). The result is a simply connected space $B \Gamma_{g,n}^{+}$ and an integral homology equivalence $B \Gamma_{g,n} \to B  \Gamma_{g,n}^{+}$. Quillen's plus construction is functorial (up to homotopy). The glueing map $\Gamma_{\infty, n} \to \Gamma_{\infty, n-1}$ are homology equivalences by Harer's theorem and consequently, we obtain homotopy equivalences 

$$B \Gamma_{\infty,n}^{+} \to B \Gamma_{\infty,n-1}^{+},$$

and so we can write $B \Gamma_{\infty}^{+} := B \Gamma_{\infty,n}^{+}$ without danger of confusion.\\
By elementary obstruction theory, any map $B \Gamma_{\infty,n} \to Y$, where $Y$ is a simple space (the fundamental group is abelian and acts trivially on all higher homotopy groups), extends uniquely (up to homotopy) to a map $B \Gamma_{\infty}^{+} \to Y$. Any grouplike $H$-space $Y$ (i.e. $\pi_0 (Y)$ is a group) is simple, in particular $\Omega^{\infty }\bG_{-2}^{\SO}$ and $\Omega_{0}^{\infty }\bG_{-2}^{\SO}$, and the map $\alpha:B \Gamma_{\infty,n} \to \Omega^{\infty}_{0}\bG_{-2}^{\SO}$ yields a map $\alpha^+ : B \Gamma_{\infty}^{+} \to \Omega^{\infty}_{0} \bG_{-2}^{\SO}$. The Madsen-Weiss theorem can be reformulated as the statement that $\alpha^+ $ is a homotopy equivalence.
It follows that $\bZ \times B \Gamma_{\infty}^{+} \simeq \Omega^{\infty} \bG_{-2}^{\SO}$. In particular, $\bZ \times B \Gamma_{\infty}^{+}$ has the homotopy type of an infinite loop space. This was proven by Ulrike Tillmann (\cite{Till}) before the spectrum $\bG_{-2}^{\SO}$ and the Madsen-Tillmann map was constructed. Her proof relies on abstract infinite loop space machinery (\cite{Ads},\cite{Seg}). Tillmann's proof gives a connective spectrum $B^{\infty}(\bZ \times B \Gamma_{\infty}^{+})$ (this is Adams' notation, \cite{Ads}) and a homotopy equivalence $\Omega^{\infty} (B^{\infty}(\bZ \times B \Gamma_{\infty}^{+})) \simeq \bZ \times B \Gamma_{\infty}^{+}$.\\
Moreover, there is a map $\Sigma^{\infty} (\bZ \times B \Gamma_{\infty,+}^{+})
\to B^{\infty} (\bZ \times B \Gamma_{\infty}^{+})$ adjoint to the identity on $\bZ \times B \Gamma_{\infty}^{+}$.\\
In the paper \cite{MT}, it is also shown that the map $\alpha: \bZ \times B \Gamma_{\infty} \to \Omega^{\infty} \bG_{-2}^{\SO}$ is an infinite loop map. Thus Tillmanns infinite loop space structure coincides with the infinite loop space structure provided by the Madsen-Weiss theorem.\\
Also in \cite{Till}, Tillmann showed that the map $\bZ \times B \Gamma_{\infty}^{+} \to \bZ \times B \Gl_{\infty} (\bZ)^{+}$, which is given by the action of the mapping class on the first homology of the surface, is an infinite loop map. Thus it comes from a map of spectra $\rho:B^{\infty}(\bZ \times B \Gamma_{\infty}^{+}) \to K \bZ$ (the latter is the algebraic $K$-theory spectrum of $\bZ$, see section \ref{indexamples} for details).\\

\paragraph*{The spin case}
The second example of $\gF$-structures is $\gF = B \Spin (2) = \bC \bP^{\infty}$, and the map $\gF \to B \Gl_2 (\bR)$ is the classifying map of the \emph{square} of the universal complex line bundle. An $\gF$-structure on a surface bundle is the same as a spin structure. The universal Madsen-Tillmann map in this case is 

$$\alpha:B \Spiff(M) \to \Omega^{\infty} \bG_{-2}^{\Spin}.$$

It can be shown that the group of components $\pi_0 (\Omega^{\infty} \bG_{-2}^{\Spin})$ is isomorphic to $\bZ \oplus \bZ/2$ and that $\alpha$ maps $B \Spiff(M)^{\epsilon}$ to the component indexed by $(g-1, \epsilon)$ (see \ref{arfmad}).
There are spin mapping class groups for surfaces with nonempty boundary and glueing maps as above (although they are more subtle). Harer proved a stability theorem for the groups $\Gamma_{g}^{\epsilon}$ (\cite{Harerspin}), which was later enhanced by Tilman Bauer (\cite{bau}) to a stability theorem for the (really relevant) groups $\hat{\Gamma}_{g}^{\epsilon}$. The known range of homological stability is much smaller in the spin case: one needs that $2k^2 +6k -2 \leq g$ for Bauer's theorem. The Madsen-Weiss theorem in the spin case is:

\begin{thm}\label{madsenweissspin} (\cite{Gal}, \cite{GMTW})
If $2k^2 +6k -2 \leq g$, then $H_k (\alpha): H_{k} (B \Spiff (M)^{\epsilon}; \bZ) \to H_k (\Omega_{(g-1,\epsilon)}^{\infty} \bG_{-2}^{\Spin} ; \bZ)$ is an isomorphism.
\end{thm}  

\subsection{Miscellaneous calculations}\label{miscellaneous}

\paragraph*{Mumfords conjecture:}
We show how the Madsen-Weiss theorem leads to a computation of the rational cohomology of $B \Gamma_g$ in the stable range. This is based on the cofiber sequence 

$$
\xymatrix{
\bG_{-2}^{SO} \ar[r]^-{\omega}&  \Sigma^{\infty} \bC \bP^{\infty} \ar[r] & \bG_{-1}^{SO}\\
}
$$

derived from \ref{cofibseq}. The spectrum $\bG_{-1}^{SO}$ is homotopy equivalent to $\Sigma^{\infty- 1} \bS^0$: its $n$th space is the Thom space $\Th (U_{1,d})$ of the orthogonal complement of the tautological real line bundle on $\Gr_{1}^{SO} (\bR^n)$, the Grassmannian of oriented lines in $\bR^n$, which is homotopy equivalent to $\bS^{n-1}$. Thus the inclusion of a fiber $\bS^{n-1} \to \bG_{-1}^{SO}$ is $(n-2)$-connected; and the sequence of all these inclusions give the desired map of spectra $\Sigma^{\infty -1} \bS^0 \to \bG_{-1}^{\SO}$, which is a homotopy equivalence.\\
Thus $\omega$ as well as $\Omega^{\infty} \omega$ induce isomorphisms in rational cohomology. It is a classical fact (\cite{Segcp}) that 
$H^{*} (Q (\bC \bP^{\infty}_{+}); \bQ)$ is a polynomial algebra in the generators $\sigma z^n$, where $\sigma: H^{*} (\bC \bP^{\infty}) \to H^{*} (Q(\bC \bP^{\infty}_+))$ is the suspension map (which is \emph{not} multiplicative). Write $y_i:=   \Omega^{\infty}(\omega)^* \sigma^* z^i$. It follows that

\begin{equation}\label{computt}
H^{*}(\Omega_{0}^{\infty} \bG_{-2}^{\SO}; \bQ) = \bQ [y_1, y_2,y_3, \ldots]. 
\end{equation}

We can identify the MMM-classes as cohomology classes coming from the spectrum cohomology of $\bG_{-2}^{SO}$. 
By \ref{gaga} and \ref{specko}, it follows that

\begin{equation}\label{laberei}
\alpha^* (\Omega^{\infty} \omega)^* \sigma z^n = \kappa_n.
\end{equation}

In view of \ref{computt} and \ref{laberei}, Mumfords conjecture \ref{mumconject} follows from the Madsen-Weiss theorem.\\

\paragraph*{The spin case}
It is not difficult to see that the map $\bG_{-2}^{\Spin} \to \bG_{-2}^{SO}$ is a rational homotopy equivalence (it is even an equivalence after localization away from $2$), see \cite{Gal}.\\
Let us study the component groups of the spectra $\bG_{-2}^{\Spin}$ and $\bG_{-2}^{SO}$. More generally, we study the component group of $\Omega^{\infty} \bG_{-d}^{\gF}$ if $\gF$ is a stable fibration which has \emph{connected} total spaces $F_d$. Start with the suspension of the cofiber sequence \ref{cofibseq}:

$$\bG_{-d}^{\gF} \to \Sigma \bG_{-(d+1)}^{\gF} \to \Sigma^{\infty +1} (F_{d+1,+}) .$$

From the long exact homotopy sequence of this cofiber sequence we see that $\pi_i (\bG_{-d}^{\gF} )\to\pi_i ( \Sigma \bG_{-(d+1)}^{\gF} )$ is an epimorphism if $i \leq 0$ and an isomorphism if $i  <0$. Inductively, it follows that 

$$\pi_0 (\Sigma \bG_{-(d+1)}^{\gF} )\cong \pi_0 (\Sigma^2 \bG_{-(d+2)}^{\gF} ) \cong \ldots \pi_0 (\Sigma^{-d} M \gF) = \pi_d (M\gF).$$

The epimorphism $\pi_0 (\bG_{-d}^{\gF} )\to \pi_0 ( \Sigma \bG_{-(d+1)}^{\gF} )$ is part of an exact sequence
\begin{equation}\label{shortext}
\bZ \cong \pi_1 (\Sigma^{\infty +1} \gF_{d+1,+}) \to \pi_0 (\bG_{-d}^{\gF}) \to \pi_0 ( \Sigma \bG_{-(d+1)}^{\gF} ) \cong \pi_d (M\gF) \to 0.
\end{equation}

We are interested in the case $d=2$ and in the cases $\gF =B \SO$ and $\gF = B \Spin$. It is well-known that $\pi_0 (\bG_{-2}^{SO}) \cong \bZ$ (see \cite{GMT}) and that $\pi_2 (MSO)=0$. We have seen that $\pi_2 (M \Spin) = \bF_2$ (Corollary \ref{spinbrod}). There is an obvious map of stable fibrations $B \Spin \to B \SO$ over $B \Gl(\bR)$. If we compare the two short exact sequences \ref{shortext} for these cases, we obtain

\begin{equation}\label{pinullrechnung}
\xymatrix{
\bZ \ar[r] \ar[d]^{\cong} & \pi_0 (\bG_{-2}^{\Spin}) \ar[d] \ar[r] & \bF_2 \ar[d] \ar[r] & 0 \\
\bZ \ar[r]^-{\cong}  & \pi_0 (\bG_{-2}^{SO}) \ar[r] & 0, & \\
}
\end{equation}

from which we see by a diagram chase that $\pi_0 (\bG_{-2}^{\Spin}) \to \pi_0 (\bG_{-2}^{SO}) \oplus \pi_2 (M\Spin) \cong \bZ \oplus \bF_2$ is an isomorphism.\\
In the case $\gF= B \SO$ and $d=2$, the element in $\pi_0(\bG_{-2}^{SO}) \cong \bZ$ represented by $\Sigma^{\infty} \bS^0 \to \bTh (\nu_M) \to \bG_{-2}^{SO}$ corresponds to $\pm \frac{1}{2} \chi (M) =\pm ( g(M)-1)$ (see \cite{GMT}). It follows

\begin{prop}\label{arfmad}
If $M$ is a surface with spin structure, then the element in $\pi_0 (\bG_{-2}^{\Spin}) \cong \bZ \oplus \bF_2$ represented by $\Sigma^{\infty} \bS^0 \to \bTh (\nu_M) \to \bG_{-2}^{\Spin}$ corresponds to $( \pm(g-1); \At(M))$.
\end{prop}

\paragraph*{The third homotopy group}
We let $d=2$ and $G = SO$ and turn to the computation of the homotopy group $\pi_3 (\bG_{-2}^{SO})$. We will make use of this computation at the very end of this work, in section \ref{conclusion}. For a pointed space $X$, denote by $X_+ $ the space $X  \cup \{ * \}$; $* \notin X$. There are obvious pointed  maps $X_+ \to X$ and $\projection:X_+ \to \bS^0$. They induce a splitting 

\begin{equation}\label{spaltt}
Q(X_+) \simeq Q(X) \times Q (\bS^0).
\end{equation}

\begin{prop}\label{drittehomotopie}
The composition 

$$
\xymatrix{
\bG_{-2}^{SO} \ar[r]^-{\omega} & \Sigma^{\infty} \bC \bP^{\infty}_{+} \ar[r]^{\projection} & \Sigma^{\infty} \bS^0\\
}
$$ 

induces an isomorphism on the third homotopy groups. Thus, $\pi_3 (\bG_{-2}^{SO}) \cong \bZ / 24$.
\end{prop}

\textbf{Proof:}
The map $\omega$ belongs to the cofiber sequence of spectra

$$\bG_{-2}^{SO} \to \Sigma^{\infty} \bC \bP^{\infty}_{+} \to \Sigma^{\infty-1} \bS^0,$$

from \ref{cofibseq}. Because $\pi_4 (\Sigma^{\infty} \bS^0) = \pi_5 (\Sigma^{\infty} \bS^0) =0$, the map $\pi_3(\omega)$ is an isomorphism.\\
A computation of Juno Mukai (\cite{Muk}) shows that the third homotopy group of the infinite loop space of $ \Sigma^{\infty} \bC \bP^{\infty}$ (without additional base point) is zero. Thus $\pi_3(\projection)$ is an isomorphism in view of \ref{spaltt}.\qed\\

\subsection{The Atiyah-Singer index theorem}\label{atiyahsinger}

We now restrict to the case of a bundle $\pi:E \to B$ of Riemann surfaces and show how the Atiyah-Singer index theorem for families of elliptic operators (\cite{AS4}) fits in the present context. A similar discussion is valid for families of complex manifolds of arbitrary dimension and for other classical operators such as the Atiyah-Singer Dirac operator for families of Spin manifolds.\\
Denote by $\beta: K^{*} \to K^{*-2}$ the Bott isomorphism. Then the Atiyah-Singer index formula for families can be stated as follows.
If $V \to E$ is a smooth complex vector bundle, then one can "twist" the Cauchy-Riemann operator with the vector bundle $V$ (\cite{Atibott}). The result is only determined up to operators of lower order, and the $K$-theory class of the index index bundle only depends on the $K$-theory class of $V$. If $V$ is holomorphic, then one can choose the Cauchy-Riemann operator of $V$ as the twisted operator. 

\begin{thm}\label{indexsatz}
Let $\pi:E \to B$ be a bundle of Riemann surfaces over a compact base space $B$. Let $V \to E$ be a complex vector bundle. Then the index of the twisted Dolbeault operator $\ind ( \delbar_V) \in K^0 (B)$ equals $\beta^{ -1} \pi_{!} ([V])$.
\end{thm}

Atiyah and Singer formulated the index theorem using a different push-forward map, which we denote by $\pi_{\sharp}$. In the sequel, the $0$-th $K$-theory group with compact support of a locally compact space $X$ will be denoted simply by $K(X)$.\\
For the definition of the push-forward, Atiyah and Singer proceed in three steps. If $f: X \to Y$ is a proper embedding with a complex structure on the normal bundle, then let $U \subset X$ be a tubular neighborhood. They define $f_{\sharp} $ as the composition $K(X) \to K(U) \to K(Y)$, where the first map is the Thom homomorphism for the normal bundle and the second map is induced by the inclusion.\\
The second case for which they define $f_{\sharp}$ is when $f:W \to X$ is the projection map of a complex vector bundle. In this case, the pushforward is defined to be the inverse of the Thom-isomorphism $K^0 (X) \to K^0(W)$.\\
In general, the pushforward of a map $f:X \to Y$ is defined by factoring $f$ into $ X \to W \to Y$, where $W$ is a complex vector bundle on $Y$. One needs the assumption that the stable normal bundle $\nu_f = f^* TY - TX$ has a complex structure (the choice of the complex structure influences the definition of $f_{\sharp}$) and that such a factorization of $f$ exists with the first map a \emph{proper} embedding. 
If $E \to B$ is a manifold bundle and if $j: E \to B \times \bR^n$ is an embedding of the fibers, then the tangential map $T_v j : T_v E \to B \times \bR^{2n}$ is a proper embedding with a canonical complex structure on the normal bundle. The symbol of any elliptic pseudo-differential operator $P$ defines a class $\sigma_p \in K^0 (T_v E)$. The index theorem in the formulation of \cite{AS4} says that the index bundle $\ind (P) \in K^{0}(B)$ is the same as $\projection_{\sharp}^{B} \circ (T_v j)_{\sharp} (\sigma_P)$.\\
The principal symbol of the Dolbeault operator is the $K$-orientation of $T_v E$ defined by the Bott-class. It is now easy to derive the formulation \ref{indexsatz}, using naturality properties of the push-forward.\\

Sometimes the index factors through the transfer.

\begin{prop}\label{transfact}
Let $E \to B$ is a bundle of Riemann surfaces. Let $V \to E$ be a complex vector bundle and let $\Lambda$ be the vertical cotangent bundle. Then 

$$\ind (\delbar_{V \otimes (1- \Lambda)} ) = \trf^{*}([V]) \in K^0 (B).$$

\end{prop}

\textbf{Proof:}
The $K$-theory Euler class of a complex line bundle $L$ is $\beta^{-1}(1- L^{\prime})$. Further, it follows from the definitions that the pushforward commutes with Bott periodicity. Thus

$$\ind (\delbar_{V \otimes (1- \Lambda)} ) = \beta^{-1}(\pi_{!}(V \beta(\chi_K(T_v E))))
 = \pi_{!}(V \chi_K(T_v E)) = \trf (V ).$$
 
\qed

\subsection{Examples for the index theorem}\label{indexamples}

I will discuss several examples of Proposition \ref{transfact}, which are "natural" in the sense that the operators exist on the universal surface bundle. We explained after \ref{grorieroc}, we can apply the index theorem with certainty only in the case that the base $B$ is compact. If we say that two maps out of the universal base space $B \Gamma_{\infty}$ are homotopic, then we mean that they are homotopic after restriction to any compact subspace.\\
Before we can study these examples, we need to say a few words about algebraic $K$-theory. Let $R$ be an associative ring with unit. Then there are obvious homomorphisms $\Gl_n (R) \times \Gl_m (R) \to \Gl_{n+m} (R)$, which turn the space

$$\coprod_{n \in \bN} B \Gl_n (R)$$

into a topological monoid. The \emph{group completion} of this monoid is

$$ \Omega B (\coprod_{n \in \bN} B \Gl_n (R)).$$

The group completion theorem (\cite{Ads}) shows that there is a natural map $c: \bZ \times B \Gl_{\infty} (R) \to  \Omega B (\coprod_{n \in \bN} B \Gl_n (R))$, which is a homology equivalence. Furthermore, the infinite loop space machinery (\cite{Seg}, \cite{Ads}) shows that there is a connective spectrum $K(R)$ and a homotopy equivalence $\Omega^{\infty} K(R) \simeq  \Omega B (\coprod_{n \in \bN} B \Gl_n (R))$. This is the \emph{algebraic $K$-theory spectrum of $R$}. Because $ \Omega B (\coprod_{n \in \bN} B \Gl_n (R))$ is an $H$-space, the map $c$ is a $H \bZ$-localization in the sense of homotopy theory (see \cite{Bous}).\\
It is well-known that the commutator subgroup $[\Gl{\infty}(R),\Gl{\infty}(R)]$ is perfect. Thus we can apply Quillen's plus construction to the space $B \Gl_{\infty}(R)$ with respect to the commutator subgroup. The result is a simple space $B \Gl_{\infty}(R)^{+}$ and a homology equivalence $B \Gl_{\infty}(R) \to B \Gl_{\infty}(R)^{+}$. Thus $\Omega K (R) \simeq \bZ \times  B \Gl_{\infty}(R)^{+}$.\\
If $R$ is a topological ring, then $\Gl_n(R) $ is a topological group and we can apply the construction of $KR$ literally. When $R = \bC$, the result is the usual connective $K$-theory spectrum $\bold{k}$.\\
Similarly, one can construct symplectic $K$-theory, replacing $\Gl_n (R)$ by $\Sp_{2n} (R)$ everywhere (we only need the case $R = \bZ$). More precisely, there is a connective spectrum $K \Sp (\bZ)$ and a homology equivalence $\bZ \times B \Sp_{\infty} (\bZ) \simeq \Omega^{\infty} K \Sp (\bZ)$. There is an obvious map $\tau: K \Sp (\bZ) \to K (\bZ)$. Note: under the isomorphisms $\pi_0 K(\bZ) \cong \bZ$ and $\pi_0 K \Sp(\bZ)\cong \bZ$, the map induces multiplication by $2$.\\
The maps $\Sp_{2n}(\bZ) \to \Sp_{2n} (\bR) \simeq U(n)$ induce a map $\zeta: K \Sp (\bZ) \to K(\bC) = \bold{k}$.
\\
Now we can describe our examples for the index theorem.
The first example we want to study is the example $V=1$, the trivial line bundle. As explained above, $\trf^{*} (1)$ is in general not $1 \in K^0 (B)$.\\
Consider the unit in $K^0 (\bC \bP^{\infty})$. It is a map of spectra $\Sigma^{\infty} \bC \bP^{\infty}_+ \to \bold{k}$, which factors through the projection $\Sigma^{\infty} \bC \bP^{\infty}_+ \to \Sigma^{\infty} \bS^0 $ and the unit map $u:\Sigma^{\infty} \bS^0 \to \bold{k}$. This unit map can be factored as $\Sigma^{\infty} \bS^0 \to K \bZ \to K \bC^{\delta} \to \bold{k}$. The first is the unit map into $K \bZ$, the spectrum $K \bC^{\delta}$ is the algebraic $K$-theory spectrum of the discrete ring $\bC$ and $K \bZ \to K \bC^{\delta}$ is induced from the map of rings $\bZ \to \bC$. The third map is the map induced by $\id_{\bC}$ on algebraic $K$-theory, viewed as a map from the discrete ring $\bC$ to the topological ring $\bC$.\\
Thus the index map $\ind \delbar_{1-\Lambda}= \trf^{*} (1)$ of the operator $\delbar_{1-\Lambda}$ factors through the sphere spectrum. In particular, all rational characteristic classes of the index bundle in degrees larger than one vanish, because $Q_0 (\bS^0)$ is rationally acyclic. This is not too surprising, since there is another description of $\ind(\delbar_{1-\Lambda})$.
Let $b \in B$ and let $C:=\pi^{-1}(b)$ be the compact Riemann surface over $b$. Then the fiber of the index bundle over $b$ is the virtual vector space $H^0(C, \cO) - H^1(C; \cO)- H^0(C, \Lambda_C) + H^1(C, \lambda_C)= \bC - H^0(C; \Lambda_C)^{\prime}- H^0(C, \Lambda_C) + \bC = 2-H^1(C, \bC)$. Here we used Serre-duality and the Hodge-decomposition. All isomorphisms are natural, and so we obtain a bundle isomorphism $\ind(\delbar_{1-\Lambda}) \cong 2- H^{1}_{v} (\pi)$. The latter is the flat bundle of first cohomology spaces. Thus the index bundle has a canonical flat structure. Moreover, the bundle $H^{1}_{v} (\pi)$ can be identified with the sum of the Hodge bundles $V_1 (\pi) \oplus \overline{V_1 (\pi)}$ and it is induced by the maps

$$B \to B \Gamma_g \to B \Sp_{2g}(\bZ) \to B \Gl_{2g}(\bC)^{\delta} \to BU(2g).$$

We have seen that the diagram of spectra

$$
\xymatrix{
B^{\infty}(\bZ \times B \Gamma_{\infty}) \ar[dd]^{\alpha} \ar[rrr]^{B\rho} & & & K \bZ \ar[r]^{u} & K \bC^{\delta} \ar[d]\\ 
 & & & &   \bold{k}\\
\bG_{-2}^{SO} \ar[r]^{\omega} & \Sigma^{\infty} \bC \bP^{\infty}_{+} \ar[r]^{\projection} & \Sigma^{\infty} \bS^0 \ar[r]^{u} & K \bZ \ar[r] &  K \bC^{\delta} \ar[u]\\
}
$$

is commutative (up to homotopy).\\
At this point, it starts to look plausible that already the two maps $B^{\infty}(\bZ \times B \Gamma_{\infty}^{+}) \to K \bZ$ are homotopic. This is indeed true, at least if we take the induced maps on infinite loop spaces. This is pointed out by Tillmann in \cite{Till2} and relies on the main theorem of \cite{DWW}. More precisely, we can say

\begin{prop}\label{ktheoryindex}(\cite{Till2}, \cite{DWW})
The maps $\bZ \times B \Gamma_{\infty}^{+} \to \bZ \times B \Gl_{\infty} (\bZ)^{+}$ given by the action on the first homology group and by the index map for the bundle $1- \Lambda$ are homotopic (as maps of spaces, not necessarily as infinite loop maps).
\end{prop}

The index bundle $\ind \delbar_V$ does not factor through $\Sigma^{\infty} \bC \bP^{\infty}_{+}$ if $V$ is not a multiple of $1- \Lambda$. An interesting example for this phenomenon is provided by the bundle $V=1$.\\
Recall the map of spectra $\lambda:\bG_{-2}^{SO} \to \Sigma^{-2}\bold{MU}$, which can be composed with the Conner-Floyd map $\xi: \bold{MU} \to \bold{k}$ to give a map $\bG_{-2}^{SO} \to \Sigma^{-2} \bold{k}$. The Bott periodicity map is a map $\Sigma^{-2} \bold{k} \to \bold{k}$.

\begin{prop}
Let $E \to B$ be a surface bundle on a compact space $B$. Then the composition $\Sigma^{\infty}B_+ \to \bG_{-2}^{\SO} \to \Sigma^{-2} \bold{k} \to \bold{k}$ represents $\ind(\delbar_1)$.
\end{prop}

This is an immediate consequence of Theorem \ref{indexsatz}, Definition \ref{defpusj}, and the fact that $\xi$ represents the Thom class of the universal complex vector bundle.\\
There is another, more naive description for $\ind(\delbar_1)$. We consider the map $B \rho: B \Gamma_g \to B\Sp_{2g}(\bZ) \to B \Sp_{2g}(\bR ) \simeq BU(g)$ from \ref{symplecto}. It classifies the universal Hodge bundle $\overline{V_1}$. Serre duality tells us that $\overline{V_1} \cong \coker(\delbar_1)$. The kernel bundle of $\delbar_1$ is the trivial complex line bundle and so we see that

$$\ind(\delbar_1) = 1- \overline{V_1}.$$

\subsection{From the Mumford conjecture to the Madsen-Tillmann map}\label{mumfordtomadsen}

Here we study in more detail the example $V= T_v E$. The index bundle $\ind (\delbar_{T (1- \Lambda)})= \ind (\delbar_{T - 1})$ is represented by the composition

$$
\xymatrix{
\Sigma^{\infty} B_+ \ar[r]^-{\alpha} &  \bG_{-2}^{SO} \ar[r]^-{\omega} &  \Sigma^{\infty} \bC \bP^{\infty}_+ \ar[r]^-{l} & \bold{k}.\\
}
$$

The whole composition pulls back the class $s_n$ to the MMM-class $\kappa_n$, because $l^{*} s_n = e^n$ and because of Proposition \ref{mummorid}. It is important to mention that $s_n $ is only a \emph{rational} spectrum cohomology class of $\bold{k}$, not an integral one. Thus the adjoint map $B \to \Omega^{\infty} \bold{k} = \bZ \times BU$ pulls back $s_n \in H^{2n} (\bZ \times BU; \bQ)$ to $\kappa_n \in H^{2n}(B; \bQ)$. The map $l$ represents the class of the canonical line bundle in $K^0 (\bC \bP^{\infty})$ and its adjoint $Q (\bC \bP^{\infty}_+) \to \Omega^{\infty} \bold{k} = \bZ \times BU$ is the rational equivalence from \cite{Segcp} mentioned above. Thus we see that rationally the integral Chern character classes of the index bundle $\ind (\delbar_{T - 1})$ agree with the MMM-classes.

\begin{rem}
A few remarks are to be done at this place. If $g \geq2$, and if $C$ is the fiber at a point $b \in B$, then the fiber of the index bundle under consideration is just the formal difference $-H^1(C, \cO(T^{(1,0)}C)) - \bC + H^1(C, \cO_C)$.
By the Kodaira-Spencer deformation theory of complex manifolds (of complex dimension 1), it follows that the vector space $H^1(C, \cO(T^{(1,0)}C))$ parameterizes all infinitesimal deformations of the complex mani\-fold $C$, that means, it should be the tangent space to the moduli space of Riemann surfaces of genus $g$. Teichm\"uller theory tells us that this is indeed true. By Serre duality $H^1(C, \cO(T^{(1,0)}C))=H^0(C, \Lambda^{\otimes 2}_{C})^{\prime}$, and this is the tangent space to $\cT_g$ (and not the cotangent space, because the definition of the Teichm\"uller map $H^0(C, \Lambda_{C}^{\otimes 2}) \to \cT_g$ involves the choice of a metric on the space of quadratic differentials.)\
There is still the problem that the moduli space is not a smooth manifold, but only an orbifold. Nevertheless, it has a complex tangent orbibundle, which has well-defined Chern classes. However, these Chern classes are only rational cohomology classes and not integral ones.\\
It should be more than a curiosity that the Chern character of the difference of these two natural bundles agrees with the tautological classes, which ultimately yield the isomorphism $\bQ[\kappa_1, \kappa_2, \ldots] \to H^{*} (B \Gamma_{\infty}; \bQ)$ predicted by the Mumford conjecture.
\end{rem}

I now show how one could \emph{find} the Madsen-Tillmann map when one starts from the Mumford conjecture \ref{mumconject}. This may be also the way Madsen and Tillmann found their map. The Mumford conjecture predicts that 

\begin{equation}\label{mumconj} \bQ[\kappa_1, \kappa_2, \ldots] \to H^*(B\Gamma_g, \bQ) 
\end{equation} 

is an isomorphism in the stable range of Harer's stability theorem.
The first observation is that the left-hand side of the conjectural isomorphism \ref{mumconj} is isomorphic to the cohomology of the infinite Grassmannian, $H^{*}(BU,\bQ) \cong \bQ[s_1, s_2,\ldots]$. We choose this system of generators instead of $c_1, c_2, \ldots$ because of their better algebraic properties (see below).\\
It is natural to look out for a map 

$$\Phi:B\Gamma_g \to BU$$

realizing the conjectural isomorphism \ref{mumconj} in rational cohomology.
By Tillmann's theorem (see \cite{Till}), $B\Gamma_{\infty}^{+}$ is an infinite loop space, in particular, $H^{*}(B\Gamma_{\infty};\bQ)$  is a Hopf algebra. The multiplication map $m: B\Gamma_{\infty}^{+} \times B\Gamma_{\infty}^{+} \to B\Gamma_{\infty}^{+}$ (homotopic to the loop sum) is homotopic to the colimit of the glueing maps

$$ a_{g,h}:B \Gamma_{g,2} \times B \Gamma_{h,2} \to B \Gamma_{g+h,2}.$$

What is the behavior of the classes $\kappa_n$ under the structural map? The answer is that they are primitive, i.e.

\begin{prop}

$$a_{g,h}^{*}\kappa_n = 1 \otimes \kappa_n + \kappa_n \otimes 1.$$

\end{prop}

\textbf{Proof:} Let $E$, $F$ be surface bundles over X of genus $g$, $h$ with two trivialized boundary components and $G$ be the bundle obtained by glueing at two boundary components. Then the vertical tangent bundles are trivialized along the boundary. Thus $e(T_v F) e(T_v E)=0$. It follows that $\kappa_n(G)=\pi_{*} e(T_v G)^{n+1} = \pi_{*} e(T_v E)^{n+1} + \pi_{*} e(T_v F)^{n+1}$.\qed\\

If we want that the conjectural map $\Phi$ is a map of $H$-spaces (this is natural, because of Tillmann's theorem), then we are almost forced to have $\Phi^{*} (n! \ch_n )= \kappa_n$, since $s_n = n!\ch_n$ generate the integral primitive classes. So, one is led to the question:

\begin{quest}
Is there a virtual vector bundle $V$ on $B\Gamma_g$ such that $n! \ch_n(V)=\kappa_n$?
\end{quest}

Index bundles give virtual vector bundles. We want an index bundle $\ind \delbar_W$ which is defined naturally for all bundles of Riemann surfaces. There are not too much possibilities, because the tensor powers of the cotangent bundle are the only natural line bundles on Riemann surfaces. The Grothendieck-Riemann-Roch theorem shows that $W= T_v -1$ is the correct choice: 

\begin{prop}\label{gleusch}
$s_n(\ind(\delbar_{T_v})-\ind(\delbar_{1}))= \kappa_n$ in rational cohomology, where $1$ denotes the trivial holomorphic line bundle.
\end{prop}

We can ask whether the equation in \ref{gleusch} is an equation of \emph{integral} cohomology classes. To consider this more closely, we define the analytical MMM-classes $\lambda_n \in H^{2n}(B \Gamma_g; \bZ)$.

\begin{defn}\label{analytical}
The \emph{analytical MMM-classes} of a surface bundle $\pi:E \to B$ are the classes $s_n(\ind(\delbar_{T_v})-\ind(\delbar_{1})) \in H^{2n}(B; \bZ)$.
\end{defn}

It is not true that the equation in Proposition \ref{gleusch} is an equation of \emph{integral} cohomology classes. The topological MMM-classes differ from the analytical ones. This can be seen as follows.
If one takes a surface $M$ and an action of a finite (cyclic groups are enough) group $G$, then one studies the bundle $E(G;M):= EG \times_G M \to BG$. There are explicit methods to compute both, the analytical and the topological MMM-classes of these examples. The topological classes $\kappa_n$ can be computed using the higher Riemann-Hurwitz formula of \cite{KU}, see also \cite{AKU}. The result only depends on the fixed-point data of the action.\\
The analytical classes can be computed using the Lefschetz-formula of Atiyah and Bott: There is an invariant complex structure on $M$; and the index bundles are just given by the appropriate combination of bundles of the form $EG \times_G H^i (M; \Lambda^{\otimes n}_M)$.\\
If one chooses a cyclic group action almost at random, then one sees that both classes differ. Also, we have shown in section \ref{spherebundles} that both classes differ for sphere bundles and in section \ref{torusbundles} that the same is true for tori.\\
The work \cite{Segcp} by Graeme Segal, which was already mentioned, contains a stronger statement than merely the rational homotopy equivalence $Q( \bC \bP^{\infty}_{+}) \simeq \bZ \times BU$. Segal shows that there exist a splitting (not an infinite loop map) $ \bZ \times BU \to Q(\bC \bP^{\infty}_{+})$, which shows that $Q(\bC \bP^{\infty}_{+}) \simeq \bZ \times BU \times F$, where $F$ is a quite mysterious space which is rationally acyclic. The difference classes $\kappa_n - \lambda_n$ seems to be related to the cohomology of the space $F$.
A systematic study of the difference classes would go beyond the scope of this work.

\subsection{The Madsen-Tillmann-diagram}

In the following diagram of spectra, we recollect the discussion of the present chapter.

\begin{equation}\label{mtddiagram}
\xymatrix{
 &  \bold{k} &   &   & \\
  & K (\bZ)\ar@{=}[r] \ar[u] & K (\bZ)  & \Sigma^{\infty} \bS^0 \ar[l]^{u} & \\
K \Sp(\bZ) \ar[ur]^{\tau} \ar[d]^{\zeta} & B^{\infty} (\bZ \times B \Gamma_{\infty}^{+}) \ar[l]^-{\rho} \ar[d]^{\beta \ind \delbar_1} \ar[u]_{\ind \delbar_{1-\lambda}} \ar[r]^-{\alpha} & \bG_{-2}^{\SO} \ar[d]^{\lambda} \ar[r]^-{\omega} & \Sigma^{\infty} \bC \bP^{\infty}_{+}  \ar[r]^-{l}  \ar[u]^{\projection} & \bold{k} \\        
\bold{k} & \Sigma^{-2} \bold{k} \ar[l]^{\beta^{-1}}  & \Sigma^{-2} \gM \gU  \ar[l]^-{\Sigma^{-2} \xi} &  & \\
}
\end{equation}

\pagebreak

\section{Divisibility of MMM-classes for spin surface bundles}\label{divisibility}

\subsection{The general case}

In this section, we study the divisibility properties of characteristic classes for surface bundles with spin structures. All our results are modulo torsion, i.e.:

\begin{conv}
If $B$ is a space and $x \in H^*(B ,\bZ)$ and $0 \neq a \in \bZ$, then the statement "$x$ is divisible by $a$" means that the image of $x$ in the free group $ H^{*}_{\mathrm{free}} (B; \bZ) := H^{*}(B, \bZ) / T$ is divisible by $a$ ($T$ denotes the torsion subgroup).
\end{conv}

The divisibility of characteristic classes of oriented surface bundles was studied by Galatius, Madsen and Tillmann in \cite{GMT}. Their main result is:

\begin{thm}(\cite{GMT})
Let $D_n$ be the maximal divisor of the class $\kappa_n$. Then $D_{2i}=2$ and $D_{2i-1}= \den (\frac{B_{i}}{2i})$.
\end{thm}

$B_i$ is the $i$th \emph{Bernoulli number}, which is defined by the expansion

\begin{equation}\label{bernoull}
\td(z):= \frac{z}{1-e^{-z}}= 1+\frac{z}{2} + \sum_{k=1}^{\infty} (-1)^{k+1} \frac{B_{k}}{(2k)!} z^{2k}.
\end{equation}

$B_k$ is a rational positive number; and von Staudt's theorem (see \cite {BorSha}, p. 410-411) gives the prime decomposition of its denominator:

$$\den (\frac{B_{k}}{2k}) = \prod_{(p-1)|2k} p^{1+ \nu_p (2k)}.$$

Note that this is always an even number. That $\den (\frac{B_{k}}{2k})$ divides $\kappa_{2n-1}$ was already noted by Morita (\cite{Mor1}) and follows from the Grothendieck-Riemann-Roch theorem in the following way. We compute $\ch (\ind (\delbar_{1}))$ in two ways. First, it is:

$$\ch (\ind (\delbar_{1})) = 1-g + \sum_{k=1}^{\infty} \frac{\gamma_k}{k!}$$

by the definition of $\gamma_k$. On the other hand, by the Grothendieck-Riemann-Roch theorem and by the definition of $\kappa_n$, we have
($c=c_1(T_v))$:

$$\ch (\ind (\delbar_{1})) = \pi_{!} (\td(c)) = 1-g + \sum_{k=1}^{\infty} (-1)^{k+1} \frac{B_{k}}{(2k)!} \kappa_{2k-1}.$$

By the way, this shows that $s_{2k}$ vanishes rationally. Comparison of coefficients leads to

\begin{equation}\label{toddd}
s_{2k-1} = \frac{B_{k}}{2k} \kappa_{2k-1}.
\end{equation}

Thus $\kappa_{2k-1}$ is divisible by $\den (\frac{B_{k}}{2k})$.\\

\subsection{The even classes}

\begin{assumption}\label{assai}
Let us now assume that the surface bundle $\pi:E \to B$ admits a spin structure. Then there exists a complex line bundle $S\to E$ and an
isomorphism $S^2 \cong \Lambda_v$. We denote $y:= -c_1 (S)$; thus $2y = c=c_1 (T_v)$.
\end{assumption}

\begin{prop}\label{obvious}
For a spin surface bundle, $\kappa_n$ is divisible by $2^{n+1}$. This holds even integrally.
\end{prop}

\textbf{Proof:}
$\kappa_n = \pi_{!}(c^{n+1}) = 2^{n+1} \pi_{!} (y^{n+1})$. \qed

\begin{thm}\label{geradteilbar}
For spin surface bundles, the class $\kappa_{2n}$ is not divisible by any nontrivial multiple of $2^{2n+1}$. This holds in the stable range for spin mapping class groups.
\end{thm}

Before we prove the result, we note that the statement is false if $g \leq 2$ and if $4n$ is larger than the homological dimension of the spin moduli space. However, our proof relies on an unstable computation for $g=0$ and on the Madsen-Weiss theorem.

\textbf{Proof:} First, we study the universal spin surface bundle of genus $0$. This is the projective bundle of the $2$-dimensional complex vector bundle

$$V:=ESU(2) \times_{SU(2)} \bC^2 \to B SU(2) \cong \bH \bP^{\infty}.$$

Its even MMM-classes are (compare \ref{hpunendlich}):

$$\kappa_{2n} = \pm 2^{2n+1} u^n,$$

where $u \in H^4 (\bH \bP^{\infty}; \bZ)$ is a polynomial generator of the ring $H^* (\bH \bP^{\infty}; \bZ) \cong \bZ[u]$. Thus it follows that the divisibility for $\kappa_{2n}$ is optimal for $g=0$.\\
We now enhance this result, where we use the Madsen-Weiss theorem in the spin case (\cite{GMTW}). Because the Arf-invariant of the spin
structure on the sphere is $0$, the Madsen-Tillmann map of the surface bundle studied above

$$\bH \bP^{\infty} \to \Omega^{\infty} \bG_{-2}^{Spin}$$

goes into the $(-1,0)$-component. The reasoning above, together with Proposition \ref{obvious} tells us that the maximal divisor of $\kappa_{2n}
\in H^{4n}(\Omega_{(-1,0)}^{\infty} \bG_{-2}^{Spin},\bZ)$ is $2^{2n+1}$. Now let $g \in \bN$ and $\epsilon \in \bZ/2$ be arbitrary. By the
Madsen-Weiss theorem in the spin case, the map

$$B\Gamma^{Spin,\epsilon}_{g}  \to \Omega_{(g-1,\epsilon)}^{\infty} \bG_{-2}^{Spin}$$

induces an isomorphism in integral homology in low degrees. This shows the claim.\qed

\subsection{The odd classes}

Now we discuss the case of the odd MMM-classes. Let $E \to B$ be a spin surface bundle. We have seen that $2^{2n} |\kappa_{2n-1}$ and that $\den
(\frac{B_{n}}{2n})|\kappa_{2n-1}$ and we have seen that in the presence of a spin structure, the divisibility of the even classes improves by a
power of $2$. It is natural to ask whether this also happens for the odd classes.

\begin{thm}\label{Clifford}
If the surface bundle $\pi:E \to B$ has a spin structure, then $\kappa_{2n-1}(\pi)$ is divisible by $2^{2n} \den (\frac{B_{n}}{2n})$ if $B$ is compact (this is probably not a serious restriction, compare the remarks after Theorem \ref{grorieroc}).
\end{thm}

We have already seen that both $2^{2n}$ and $\den (\frac{B_{n}}{2n})$ divide $\kappa_{2n-1}$. But since the denominators of the Bernoulli numbers are even, $2^{2n}$ and $ \den (\frac{B_{n}}{2n}) $ are \emph{not} coprime, and we need new arguments to prove Theorem \ref{Clifford}.\\
We keep the notations from \ref{assai} and prove Theorem \ref{Clifford} in two steps.

\begin{prop}\label{dirac}
For spin surface bundles, the class $\kappa_{2n-1}$ is divisible by $2^{2n-1} \den (\frac{B_{n}}{2n}) $
\end{prop}

\textbf{Proof:}
The proof relies on the Grothendieck-Riemann-Roch theorem. We study the index bundle of the operator\footnote{This choice looks accidental, but it is not. The reason is that $\Psi^2 (\ind \delbar_{1+S})= \ind \delbar_1$; $\Psi^2$ is the Adams operation. The choice of the bundle $1+S$ arises naturally from this requirement. Also, it should be more natural to consider the Dirac operator on spin surface bundles. This is true and the Dirac operator and the Cauchy-Riemann operator are closely related, but we prefer to use the Cauchy-Riemann operators which are more common.} $\delbar_{1+S}$. We keep the notation and compute

$$\ch (\ind \delbar_{1+S}) = \pi_{!} (\td(2y) (1+e^{-y})) = \pi_{!}(2 \td(y)).$$

Thus $v_{2n-1}:=s_{2n-1} (\ind \delbar_{1+S}) = 2 \frac{B_n}{2n} \pi_{!}(y^{2n})$ is an integral class.
By the definition of the MMM-class, it follows that

\begin{equation}\label{vau}
\kappa_{2n-1} = 2^{2n} \pi_{!}(y^{2n}) = 2^{2n-1} \frac{2n}{B_n} v_{2n-1}= 2^{2n-1} \frac{\den (\frac{B_{n}}{2n}) }{\num (\frac{B_{n}}{2n})} v_{2n-1}.
\end{equation}

Since $\num (\frac{B_{n}}{2n})$ and $2^{2n-1} \den (\frac{B_{n}}{2n}) $ are coprime by von Staudt's theorem, the proposition follows.\qed\\

Now we care about the last factor of $2$. In the proof of Proposition \ref{dirac}, we did not really use the index theorem, because only the interpretation of $\gamma_{2n-1}$ as the class of the Hodge bundle on the base and not its integrality required the use of analysis. This is not true for the proof of the next theorem, where the last factor of $2$ is detected.
In view of \ref{vau}, it suffices to show that $v_{2n-1}$ is divisible by $2$. In other words:

\begin{prop}
The class $s_{2n-1} (\ind \delbar_{1+S}) \in H^{2n}(B; \bZ)$ is divisible by $2$.
\end{prop}

\textbf{Proof:}
We compute

$$s_{2n-1} (\ind \delbar_{1+S}) = s_{2n-1}(\ind \delbar_1) + s_{2n-1}(\ind \delbar_S) = \gamma_{2n-1} + s_{2n-1}(\ind \delbar_S).$$

By \ref{vau} and \ref{toddd}, $\gamma_{2n-1}$ is divisible by $2$.\\
If the dimension of the kernel of $\delbar_S$ is constant (and thus the kernel is a vector bundle), then we see from Serre duality (\ref{serredual}) that

$$\ind \delbar_S = \ker \delbar_S - \overline{\ker \delbar_S}.$$

Here it is really important that $S^2= \Lambda_v$ as holomorphic line bundles and not merely as complex line bundles. It is a general fact that
for any complex vector bundle $V$, we have $s_k (\bar{V}) = (-1)^k s_k (V)$. Thus it follows that

$$s_{2n-1}(\ind \delbar_S) = 2 s_{2n-1}(\ker \delbar_S),$$

as required. If the dimension of the kernel is not a constant function, then the argument above also applies. To see this, one has to look at the construction of the index bundle in that case (see \cite{LM}).\qed\\

\pagebreak

\section{The icosahedral group and $\pi_{3}(B \Gamma_{\infty}^{+})$}\label{icosahedral}

One possible construction for elements in the homotopy groups of $B\Gamma_{g}^{+}$ arises in the follow\-ing way. Take a homology $n$-sphere $M$ and a representation $\rho:\pi_1(M) \to \Gamma_g$. We obtain a homotopy class of maps $B\rho: M \to B \pi_1(M) \to B\Gamma_g$, or an isomorphism class of surface bundles on $M$. Note, however, that the representation $\rho$ does not need to be induced by an action of $\pi_1 (M)$ on $\Sigma_g$ if $\pi_1(M)$ is not finite.\\
Since the fundamental group of a homology sphere $M$ is perfect, we can apply Quillen's plus construction to the whole fundamental group. The result is a $1$-connected space $M^+$ and a homology equivalence $M \to M^+$. The Hurewicz theorem tells us that $\pi_n (M^+)\cong \bZ$ and any generator of this group yields a homotopy equivalence $M^+ \simeq \bS^n$.\\
Thus, if we apply the plus construction to the map $B\rho$, we obtain a homotopy class of maps $\bS^n \to B \Gamma_g$, in other words, an element in $\pi_n(B \Gamma_{g}^{+})$. A similar construction, applied to the general linear group instead of the mapping class group, was considered by Jones and Westbury in their paper \cite{JonesWestb}, which we will use later and which was the stimulus for considering the icosahedral group in connection with the homotopy of the mapping class group.\\
The most famous case of a homology sphere, which was also the first example to be discovered is that of the Poincar\'e sphere, which I will
describe now.

\subsection{The icosahedral group and the Poincar\'e sphere}\label{poincaree}

Consider a regular icosahedron in Euclidean $3$-space, centered at $0$. It has $20$ faces (which are triangles), $12$ vertices (at every vertex,
exactly $5$ edges and $5$ faces meet) and $30$ edges. Let $G$ be the symmetry group of the icosahedron. It acts transitively on the vertices as
well as on the edges and faces of the icosahedron. We denote the isotropy subgroups of a chosen edge, face, vertex, respectively, by $G_1$,
$G_2$ and $G_3$, respectively. The orders of this subgroups are $2,3,5$, respectively.
Thus $G$ is a subgroup of $SO(3)$ of order $60$. It is well-known that $G$ is isomorphic to the groups $\bP \Sl_2(\bF_5)$ and $A_5$ (the alternating group). In particular, $G$ is perfect.\\
There is a universal central extension $\hat{G}$ of $G$ by $\bZ/2$ which is isomorphic to $\Sl_2(\bF_5)$. We call this the \emph{binary icosahedral group}. It can be obtained by taking the preimage of $G$ under the $2$-fold covering $SU(2) \to SO(3)$. The group $\hat{G}$ is also perfect and its center is the same as the kernel of the map $\phi:\hat{G} \to G$ and contains exactly one nontrivial element $h$. The quotient $ \hat{G} \backslash SU(2) $ is the Poincar\'e homology $3$-sphere.\\
The Poincar\'e sphere also can be described as a Seifert $3$-manifold (\cite{Saveliev}). I give a short description of Seifert homology spheres.\\
Let $n \in \bN$, and let $(a_i, b_i)$ be pairs of coprime integers, $i=1, \ldots ,n$. Take $n$ disjoint embedded discs in $\bS^2$ and let $F:=
\bS^2 \setminus (\bD^{2}_{1} \cup \ldots \cup \bD^{2}_{n})$. Consider $\bS^1 \times F$. Let $h$ be the curve represented by $\bS^1$ and let
$x_i$ be the boundary of $\bD^{2}_{i}$. Then $n$ solid tori are glued in such that the meridian of the $i$th solid torus is glued to the curve
$a_i x_i + b_i h$. The result is a closed manifold, the Seifert manifold $M((a_1,b_1), \ldots , (a_n,b_n))$. By the Seifert-van Kampen theorem,
it follows that its fundamental group has the presentation

$$\langle h,x_1, \ldots ,x_n|[h,x_i]=1, x_1 x_2 \ldots x_n =1,x_{i}^{a_i}=h^{-b_i} \rangle .$$

By the Hurewicz theorem, $M((a_1,b_1), \ldots , (a_n,b_n))$ is an integral homology sphere if and only if $a_1 \ldots a_n \sum_{i=1}^{n} \frac{b_i}{a_i}= \pm 1$. \\
If we take $a_1 =2, a_2=3, a_3=5, b_1 = -1, b_2 =1, b_3=1$, then we obtain the Poincar\'e sphere. This is not hard to see: Consider the quotient
map $\hat{G} \backslash \bS^3 \to \hat{G} \backslash \bS^3 / \bS^1 = \hat{G} \backslash \bC \bP^1 $. The latter is a compact Riemann surface
which has necessarily genus $0$ (If $\bC \bP^1 \to X$ is a nonconstant holomorphic map of Riemann surfaces, then $X \cong \bC \bP^1$). The
quotient map $M \to \bS^2 = G \backslash \bC \bP^1 $ is not a fibration, but it has three singular fibers which lie over the orbit of the
vertices, the orbit of the midpoints of the edges and the orbit of the midpoints of the faces of the icosahedron. With this construction, it is
easy to derive the presentation of the Poincar\'e sphere as a Seifert sphere. This also gives us a presentation of the group $\hat{G}$:

$$\hat{G} \cong \langle h,x_1,x_2,x_3 | [x_i, h]=1, x_1 x_2 x_3 =1, x_{1}^{2} = x_{2}^{-3} = x_{3}^{-5} =h\rangle.$$

The relation $h^2=1$ does not seem to follow immediately from this presentation. From this presentation we derive a presentation of $G$:

$$G \cong \langle y_1,y_2,y_3|y_1 y_2 y_3 =1;y_{1}^{2} = y_{2}^{-3} = y_{3}^{-5} =1\rangle.$$

The map $\hat{G} \to G$ is given by $h \mapsto 1, x_i \mapsto y_i$. $y_i$ generates an isotropy group $G_i \subset G$.\\
The Poincar\'e sphere seems to be the only example of a homology 3-sphere with \emph{finite} fundamental group (cf \cite{Milgro}).

\subsection{Surfaces with an action of the icosahedral group}\label{kleinsaction}

Now we construct certain actions of the binary icosahedral group on Riemann surfaces which in the end will give interesting elements in $\pi_3(B
\Gamma_{g}^{+})$. It is likely that these actions were already known to Felix Klein. The construction is based on an easy lemma. Let $\cO(k)$ be
the $k$th tensor power of the Hopf bundle on $\bC \bP^1$. Recall that the global holomorphic sections of $\cO(k)$ can be identified with the
vector space of homogeneous polynomials of degree $k$ on $\bC^2$. Further, $\cO(k)$ is an $\Sl_2(\bC)$-equivariant bundle over the
$\Sl_2(\bC)$-space $\bC \bP^1$. If $k$ is odd, then the central element $- 1$ acts as $- 1$ on $\cO(k)$ (i.e. on any fiber of the bundle). If
$k$ is even, then $-1$ acts trivially on $\cO(k)$, i.e. the action descends to an action of $\bP \Sl_2(\bC)$.

\begin{lem}
Let $G \subset \bP \Sl_2(\bC)$ be a finite subgroup and let $\hat{G} \subset \Sl_2(\bC)$ be the extension of $G$ by $\bZ/2$. Let $m \in \bN$ be positive. Let $s$ be a $G$-invariant holomorphic section of $\cO(2m)$ having only simple zeroes. Then there exists a hyperelliptic Riemann surface $f:X \to \bC \bP^1$ with a $\hat{G}$-action, such that $f$ is equivariant. The construction has the following properties.\\
If $m$ is odd, then the central element $h \in \hat{G}$ is the hyperelliptic involution, if $m$ is even, then $h$ acts trivially on $X$.\\
Further, the branch points of $f$ are precisely the zeroes of $s$.\footnote{Conversely, the automorphism group of any hyperelliptic surface of
genus $g \geq 2$ is a possibly trivial extension of a finite subgroup of $\bP \Gl_2(\bC)$ by $\bZ/2$.}
\end{lem}

\textbf{Proof:} Let $S \subset \cO(2m)$ be the graph of the section $s$. It is a surface of genus $0$ and it is stable under the $G$-action on
$\cO(2m)$. Let $q: \cO(m) \to \cO(2m)$ be the squaring map, let $X:=q^{-1}(S)$ and let $f:= q|_{X}$. Clearly, $X$ has a $\hat{G}$-action and $f$
is equivariant. Also, a generic point of $S$ has exactly $2$ preimages, being permuted by the involution $-1$ on the fibers on $\cO(m)$.
We need to show that $X$ is a smooth connected Riemann surface. The smoothness of $X$ is equivalent to the condition that the zeroes are simple (the locus $\{(x,y) \in \bC^2| y^2 = x^m\}$ is smooth if and only if $m=1$), and the connectivity is clear since $s$ has precisely $2m >0$ zeroes.\\
We have seen that the antipodal map of $\cO(m)$ induces the hyperelliptic involution. This implies the second statement.
The last sentence is clear from the construction. \qed\\

From now on, let $G$ again be the icosahedral group. I will give three examples of the construction above.

\begin{expl}\label{1}
Let $z_1, \ldots , z_{30} $ be the midpoints of the edges of the icosahedron, considered as points on $\bC \bP^1$ instead on $\bS^2$. The precise value of the points does not play a significant role in this discussion.\\
Now we take holomorphic sections $s_i$, $i=1, \ldots ,30$ of the Hopf bundle $\cO(1)$ (alias square root of the tangent bundle of $\bC \bP^1$) having a simple zero at $z_i$ and being nonzero elsewhere. Such sections exist and are unique up to multiplication with a complex constant. Set $s := s_1 \otimes \ldots \otimes s_{30} \in H^0(\bC \bP^1, \cO(30))$.\\
For $g \in \hat{G}$, there exists a $c(g) \in \bC^{\times}$ with $g s = c(g) s$, because $gs$ has the same zeroes as $s$. The map $c:g \mapsto c(g)$ is a homomorphism $\hat{G} \to \bC^{\times}$. Since $\hat{G}$ has no Abelian quotient, $c$ is constant. Thus, $s$ is an invariant section.\\
If we apply the construction of the lemma to $s$, we obtain a surface of genus $14$ (by the Riemann-Hurwitz formula) with a $\hat{G}$-action.
\end{expl}

\begin{expl}\label{2} Similarly, we can take the $12$ vertices as branch points for the hyperelliptic covering. The result is a surface of genus $5$. The central element $h$ acts trivially on the surface. It is the surface of genus $5$ which will give the best result (i.e., the order of the resulting element in $\pi_3 ((B \Gamma_g)^+)$ is the largest).
\end{expl}

\begin{expl}\label{3}  Similarly, we can take the $20$ midpoints of the faces of the icosahedron as branch points. The result is a surface of genus $9$ with an action of $\hat{G}$. The central element acts as the identity.\\
\end{expl}

\begin{prop}
\begin{enumerate}
\item The number of fixed points of the elements $x_1,x_2,x_3 \in \hat{G}$ in Example \ref{1} is $2, 0,0$, respectively.\\
\item The number of fixed points of the elements $x_1,x_2,x_3 \in \hat{G}$ in Example \ref{2} is $4, 2,4$, respectively.\\
\item The number of fixed points of the elements $x_1,x_2,x_3 \in \hat{G}$ in Example \ref{3} is $4, 4,2$, respectively.\\
\end{enumerate}
\end{prop}

\textbf{Proof:}\\
\textbf{Example \ref{1}:}
Since $h$ is the hyperelliptic involution, it has precisely $2g+2 =30$ fixed points, namely the branch points which lie over the midpoints of the edges of the icosahed\-ron.\\
There are two fixed points of $y_1$ on $\bP^1$, and they are two opposite branch points (i.e. midpoints of edges). Thus, $x_1$ has exactly 2 fixed points on $X$.\\
$y_2$ fixes precisely two opposite midpoints of faces. Hence all fixed points of $x_2$ lie over these two points. But if $x_2$ would have at least one fixed point, then this must also be a fixed point of $x_{2}^{3} =h$. This is impossible, since $h$ has already the $30$ fixed points mentioned above and since no nontrivial automorphism of a Riemann surface of genus $g \leq 2$ can have more than $2g+2$ fixed points (\cite{FaKr}, p.257). Thus $x_2$ is fixed-point-free.\\
The same argument shows that $x_3$ is also fixed-point free.\\
\textbf{Example \ref{2}:}
Recall that $h$ acts as the identity and that the branch points are precisely the midpoints of the $20$ faces of the icosahedron, two of which are fixed by $y_2$.\\
Down on $\bP^1$, $y_2$ has exactly $2$ fixed points. Since both are branch points, it follows that the number of fixed points of $x_2$ is $2$.\\
Down on $\bP^1$, $y_3$ has exactly $2$ fixed points, which are not branch points and thus have $2$ preimages. The two points in the preimage are permuted by $x_3$. Because $x_3$ acting on the surface has \emph{odd order}, the permutation of the two preimages must be the identity. Thus $x_3$ has $4$ fixed points.\\
The consideration for $x_1$ is slightly more complicated. $y_1$ has exactly $2$ fixed points on $\bP^1$ (they are midpoints of two opposite edges of the icosahedron). A more careful look at the action of $h$ and $x_1$ shows that the number of fixed points must be $4$.\\
\textbf{Example \ref{3}:} Now the vertices of the icosahedron are the branch points and the central element $h$ acts as the identity.\\
The generator $y_3$ has two fixed points, and since they are branch points, $x_3$ also has $2$ fixed points.\\
$y_2$ has $2$ fixed points on $\bP^1$ and so $x_2$ has $4$ fixed points (since its order is odd and the fixed points lie over the two fixed points of $x_2$).\\
$x_1$ has $4$ fixed points. To see this, let $q \in \bC \bP^1$ be the midpoint of an edge fixed by $x_1$ and let $p_1, p_2$ be the adjacent
vertices. The preimage of the edge containing $q$ under the hyperelliptic covering is a simple closed curve $S$ which is fixed by $x_1$ (as a
set, not pointwise). $y_1$ defines a diffeomorphism $f$ of the curve. Since $x_1$ changes the adjacent faces of the edge, the resulting
diffeomorphism is \emph{orientation-reversing} and interchanges two opposite points (the preimages of $p_i$). Further, we can choose a
parameterization of $S$ by $\bS^1$ and of the arc in $\bP^1$ by $[-1,1]$, such that the branch points correspond to $\pm$ and the projection to
the arc in $\bP^1$ to the real part and such that $f$ is linear. Then it follows that $f(z) = - \bar{z}$. This shows that $x_2$ has two fixed
points on the curve and thus the total number of fixed points is $4$.
\qed\\

\subsection{Review of the work by Jones-Westbury}\label{joneswestbury}

Now I will describe the work \cite{JonesWestb} by Jones and Westbury. They consider the following situation. Assume that $M$ is a homology
$n$-sphere and assume that $\rho: \pi_1 (M) \to \Sl_{N}(\bC)$ is a representation (in fact, any representation $\pi_1 (M) \to \Gl_N (\bC)$ takes
values in $\Sl_N (\bC)$, because $\pi_1 (M)$ is perfect). The representation defines a map $B \rho: M \to B\pi_1(M) \to B \Gl_{N}(
\bC^{\delta})$ (the latter is the classifying space of $\Gl_N(\bC)$, considered as a discrete group). After stabilization ($\Gl_N \to
\Gl_{\infty} = \colim \Gl_N $) and application of the Quillen plus construction to, we obtain a homotopy class of maps $\bS^n \simeq M^{+} \to
(B\Gl_{\infty}(\bC)^{\delta})^{+}$, or, in other words, an element in $\pi_{n}((B\Gl_{\infty}(\bC)^{\delta})^{+}) = K_n(\bC^{\delta})$.
More precisely, $\Sl_{\infty}(\bC^{\delta}) \subset \Gl_{\infty}(\bC^{\delta})$ is the commutator subgroup and it is perfect. Thus we can take the Quillen plus construction of $B \Gl_{\infty}(\bC^{\delta})$ with respect to the group $\Sl_{\infty}(\bC^{\delta})$. The result is the unit component $(K \bC^{\delta})_0$ of the K-theory space of the discrete ring $\bC^{\delta}$. Because the representation takes values in $\Sl_{\infty}(\bC)$, the image of the map induced by $B\rho$ on fundamental groups is in the group $\Sl_{\infty}(\bC)$. Hence the Quillen plus construction applies to $B \rho$.\\
Further, for odd values of $n$, there exists a homomorphism $e: K_n (\bC^{\delta}) \to \bC/\bZ$. The name for this homomorphism results from the
fact that it extends the classical $e$-invariant $e: \pi_{2n+1}^{st} \to \bQ/\bZ$ defined by Adams. Let us explain the meaning of the last
phrase. By the Barratt-Priddy-Quillen theorem (see, for example \cite{Seg}), there is a homotopy equivalence $\bZ \times B \Sigma_{\infty} ^{+} \simeq Q(\bS^0)$, where the plus construction on the left hand side is taken with respect to the alternating group, which is perfect. The obvious representation $\Sigma_{\infty} \to \Gl_{\infty}(\bZ)$, together with the inclusion $\bZ \subset \bC^{\delta}$, gives us a map $q: Q(\bS^0) \to K\bZ \to K \bC^{\delta}$ and thus a homomorphism $ \pi_{2n+1}^{st} \to K_{2n+1}(\bC^{\delta})$ (\cite{Quil}) \footnote{An alternative way to see the map $q$: Recall that the K-theory space of a ring $R$ is an infinite loop space. A map $\bS^0 \to K R$ is given by sending the non-basepoint in $\bS^0$ to a point in the 1-component of $KR$. By the universal property of the functor $Q$ from spaces to infinite loop spaces, this map extends uniquely (up to homotopy) to $Q (\bS^0)$.}.\\
Let us also recall the classical fact that $e: \pi^{st}_{3} \to \frac{1}{24}\bZ/\bZ$ is an isomorphism.\\
Jones and Westbury considered the problem of computing the $e$-invariant for those classes in $K_{3}(\bC^{\delta})$ defined by Seifert-spheres and representations of their fundamental groups, at least for special classes of homology spheres and special classes of representations.\\
Let $M=M((a_1,b_1), \ldots , (a_n,b_n))$ be a Seifert homology sphere and let $\rho: \pi_1 (M) \to \Sl_N(\bC) $ be a representation (alias flat
vector bundle on $M$), such that $\rho(h) $ acts as a scalar $\lambda_h$. Define

$$\alpha_i:= \rho(x_i)$$

and let

$$\lambda_{1}(i), \ldots , \lambda_{N}(i)$$

be the eigenvalues of $\alpha_i$. We have

$$\lambda_h =: \zeta_{N}^{r_h},\zeta_N := \exp(\frac{2\pi i}{N}).$$

We define rational numbers $s_k(j)$ by the relation

$$\lambda_{k}(j)=: \zeta_{N a_j}^{Ns_k(j) - b_j r_h}.$$

Let $e$ be the $e$-invariant of the element in $K_3(\bC^{\delta})$ defined by the map $M \to B\Sl_N(\bC^{\delta}) \to B\Gl(\bC^{\delta})$ given
by the representation $\rho$. Theorem C in \cite{JonesWestb} states that

$$2 \Re Ne=- a \sum_{j=1}^{n} \sum_{k=1}^{N} \sum_{l=1}^{N} \frac{(s_k(j)-s_l(j))^2}{2a_{j}^{2}} \in \bC / \bZ,$$

where $a:= a_1 \ldots a_n$. In the case where $h$ acts as the identity, there is a simpler formula to compute the $e$-invariant
(\cite{JonesWestb}, Lemma 5.3 and Theorem A), namely

$$e=-\sum_{j=1}^{n} \sum_{k=1}^{N}\frac{a s_{k}(j)^{2}}{2a_{j}^{2}}.$$

The proof of these formulae depend strongly on the Atiyah-Patodi-Singer index theorem for manifolds with boundary.

\subsection{Application of the Jones-Westbury formula to the mapping class group}\label{applijonwest}

In this section, we compute the $e$-invariants of the representations of the icosahedral group on the first homology spaces of the surfaces from
the examples \ref{1}, \ref{2}, \ref{3} constructed above. We keep the notation from the last subsections. Recall that we have

$$(a_1,b_1)= (2,-1); (a_2,b_2)=(3,1 ); (a_2,b_2)=(5,1 ).$$

\begin{prop}
The $e$-invariant of the element in $K_3(\bZ)$ given by example \ref{1} has order $6$, $12$ or $24$ in $\bQ / \bZ$.
\end{prop}

We cannot obtain a sharper result with these methods.

\textbf{Proof:} Because $h$ is the hyperelliptic involution, it acts as $- 1$ on the first homology group.\\
Since $N=28$, we have $r_h=14$.\\
Since $x_1$ has $2$ fixed points, it follows that $\tr(\alpha_1)=0$ by the Lefschetz fixed point formula. Because $x_1^2=h$, the eigenvalues are $\pm i$, the multiplicities being 14 for both signs.\\
Thus $s_k(1)=0$ for $k=1, \ldots ,14$; $s_k(1)=1$ for $k=15, \ldots ,28$. The first summand in the Jones-Westbury formula amounts to

$$56 e_1 =-30  \sum_{k=1}^{28} \sum_{l=1}^{28} \frac{(s_k(1)-s_l(1))^2}{2a_1^2} =- \cdot 30 \cdot 49 \equiv 0 \pmod{ \bZ}.$$

Because $x_2$ is fixed-point-free, it follows that $\tr(\alpha_2)= 2$. Because $x_2^3=h$, it follows that the eigenvalues of $\alpha_2$ are of
the form $\zeta_{6}^{j}$ with $j=1,3$ or $5$. Let $\mu_i$ be the number of eigenvalues with power $i$. We must determine $\mu_i$. The relevant
equations are

\begin{enumerate}
\item $\mu_1 +\mu_3 +\mu_5 = 28$;
\item $\zeta_6 \mu_1 + \zeta_{6}^{3} \mu_2 +\zeta_{6}^{5} \mu_3 = 2$ (trace formula);
\item $\mu_1 = \mu_5$ ($\alpha_1$ is real).
\end{enumerate}

These equations have a unique integral nonnegative solution, $(\mu_1 ,\mu_3, \mu_5) =(10,8,10)$. Thus $
s_1 (2) = \ldots = s_{10}(2) = 1$; $s_{11}(2) = \ldots = s_{20}(2) =0$; $s_{21}(2) = \ldots = s_{28}(2) =2$. \\
The contribution to the $e$-invariant is

$$56 e_2 = -\frac{30}{18} (1000) = -\frac{5}{3} (1000)\equiv   -\frac{5}{3} \pmod{ \bZ} \equiv \frac{1}{3} \pmod{ \bZ}.$$

The computation of the third summand follows the same pattern. We end up with $s_1 (3)= \ldots s_6 (3)= 1$, $s_7 (3)= \ldots s_{12}(3) = 2$, $s_{13}(3) = \ldots s_{18} (3)= 4$, $s_{19}(3) = \ldots s_{24}(3) = 0$, $s_{25} (3)= \ldots s_{28} (3)= 3$. It follows that $56 e_3 = - \frac{30}{50}(36 + 6) 80 \equiv 0 \pmod{ \bZ}$.\\
In the end, we obtain $56 e \equiv \frac{1}{3} $.\\
Choose a representative of $e$ in $\bQ$, also denoted $e$. Because $e$ must be annihilated by $24$ in $\bQ/\bZ$ (see section \ref{conclusion}), we can write $e= \frac{s}{24}$ for $s \in \bZ$. We obtain $\frac{56 s}{24}= r + \frac{1}{3}$ which implies the proposition. \qed\\

\begin{prop}\label{einvaria}
For Example \ref{2}, the $e$-invariant is $\frac{1}{2}$.\\
For Example \ref{3}, the $e$-invariant is $-\frac{1}{12}$.\\
\end{prop}

\textbf{Proof:}
We keep the notation as before and treat first Example \ref{2}. Then $N=18$. By the count of fixed points and the Lefschetz formula, we see that $\tr(\alpha_1)=-2$, $\tr(\alpha_2)=0$, $\tr(\alpha_3)=-2$. \\
The eigenvalues of $\alpha_1$ are of the form $\pm 1$. It follows that $s_1(1)= \ldots s_8(1)=0$ and $s_9(1) = \ldots s_{18}(1)=1$. The contribution to the $e$-invariant is $e_1= -30 \sum_{k=1}^{18}\frac{s_k(1)^2}{8} = -\frac{75}{2}\equiv \frac{1}{2} \pmod{ \bZ}$ (we use the simplification of the Jones-Westbury formula).\\
For $y_2$, it follows that $s_1(2)=\ldots =s_6(2)=0$, $s_7(2)=\ldots = s_{12}(2) = 1$, $s_{13} (2)= \ldots =s_{18} (2)=2$ and that $e_2 = \frac{-900}{18} \equiv 0 \pmod{ \bZ}$. \\
For the last summand, we obtain $s_1(3)=s_2(3)=0$; $s_3 (3)= \ldots = s_6(3) =1 $, $s_7 (3)= \ldots = s_{10}(3) =2 $; $s_{11}(3) = \ldots = s_{14}(3) = 3$, $s_{15}(3) = \ldots = s_{18}(3) = 4$. The summand for the $e$-invariant is $-180 \equiv 0 \pmod{ \bZ}$.\\
Summing everything up, gives us $e= \frac{1}{2}$.\\
For Example \ref{3}, we only give the results, because the computation follows the same pattern.\\
$s_{1}(1) = \ldots = s_{4}(1) = 0$, $s_{5}(1) = \ldots = s_{10}(1) = 1$. Thus $e_1 = \frac{-75}{4}$.\\
$s_{1}(2) = \ldots = s_{2}(2) = 0$, $s_{3}(2) = \ldots = s_{6}(2) = 1$; $s_{7}(2) = \ldots = s_{10}(2) = 2$. Thus $e_2= \frac{-100}{3}$.\\
$s_{1}(3) =  s_{2}(3) = 0$, $s_{3}(3) = s_{4}(3) = 1$, $s_{5}(3)  = s_{6}(3) = 2$; $s_{7}(3) = s_{8}(3) = 3$, $s_{9}(3) =  s_{10}(3) =4 $. Hence $e_3 \equiv 0 \pmod{ \bZ}$.\\
In the end, we obtain $e\equiv  -\frac{1}{12}$. \qed\\

\subsection{Increasing the genus}\label{increasing}

Example \ref{3} produces an element $\pi_3(B \Gamma_{5}^{+})$ whose $e$-invariant is $\frac{1}{12}$. The question arises whether we can increase the genus of the surface, since $\pi_3(B \Gamma_5)^+$ is certainly not in the stable range for Harer stability.\\
Geometrically, the stabilization procedure is as follows. Consider our surface $X$ of genus $5$ with an action of the icosahedral group $G$,
which is \emph{faithful} (the central element $h \in \hat{G}$ acts as the identity on $X$). Then we choose a point $p \in X$ which is not fixed
by any nontrivial $g \in G$ and an embedding $j:\bD^2 \to X$, $p \in D= j(\bD^2)$, such that $gD \cap hD = \emptyset$ whenever $g,h \in G$, $g
\neq h $. Let $F_n$ be a surface of genus $n$ with one parameterized boundary component. Define $X_n$ as the pushout of the diagram

$$\xymatrix{  \coprod_{g \in G} \bS^1  \ar[r] \ar[d] & \coprod_{g \in G} F_n\\
X \setminus (\cup_{g \in G} gD) & \\
}$$

Then $X_n$ is a closed smooth surface of genus $5 + 60 n$ with an faithful $G$-action. The representation of $\hat{G}$ on $H_1(X_n)$ is the
direct sum $H_1(X) \oplus 2n \cdot \phi^{*} \bZ G$, where $\phi: \hat{G} \to G$ was the extension homomorphism and $\bZ G$ denotes the regular
representation of $G$. Since the $e$-invariant is additive as a map $R \hat{G} \to \bQ / \bZ$ (see \cite{JonesWestb}), we need to compute the
$e$-invariant of the regular representation of $G$, pulled back to $\hat{G}$. We do this with the Jones-Westbury formula again.
Because $\hat{G}$ is perfect, the regular representation takes values in the special linear group.\\
Let $\alpha: \hat{G} \to \Sl_{120} (  \bC)$ represent $2 \cdot \phi^* \bZ G$. Certainly, the central element acts as the identity. If we
restrict the representation $2 \cdot \bZ G$ to the subgroup $\langle y_1\rangle \subset G$ of order $2$, it is the sum of $60$ copies of the
regular representation of $\bZ/2$. Thus, the contribution to the $e$-invariant is $e_1 = - \frac{30}{8} \cdot 60 \equiv  0 \pmod{ \bZ}$.
Similarly, the second summand in the Jones-Westbury formula becomes $- \frac{1000 }{3} \equiv -\frac{1}{3} \pmod{ \bZ}$ (the representation
restricted to $\langle y_2 \rangle$ is $40$ times the regular representation of $\bZ/3$. The third summand is $e_3  \equiv 0 \pmod{ \bZ}$. Thus
we obtain

\begin{prop}
The $e$-invariant of the action $\hat{G} \actson X_n$ is $-\frac{1}{12} - \frac{n}{3}$.
\end{prop}

In particular, we do not obtain element with $e$-invariants of higher order.\\
It seems reasonable that this result is the best possible, i.e. that there does not exist an action of the binary icosahedral group on some
surface such that the resulting element in $\pi_3(B\Gamma_{\infty}^{+})$ has a smaller $e$-invariant. To prove this, it would suffice to show
that the $e$-invariant of any \emph{integral} representation of $\hat{G}$ belongs to $(\frac{1}{12} \bZ) / \bZ \subset \bQ / \bZ$. How can one
show this? The $e$-invariant defines an additive homomorphisms $R\hat{G} \to \bC / \bZ$. The representation ring of the finite group $\hat{G}$
is a free finitely generated abelian group and one has a rather small explicit basis. By direct checking it is possible to find a $\bZ$-basis
for the subgroup of $R\hat{G}$ consisting of representations which are defined over the integers. However, the amount of time and paper one
needs for this may be quite large and so we refrain from doing anything in this direction.

\subsection{Concluding remarks}\label{conclusion}

Proposition \ref{drittehomotopie} and the Madsen-Weiss theorem show that $\pi_3(B \Gamma_{\infty}^{+}) \cong \bZ/24$ and an isomorphism is
induced by $B \Gamma_{\infty}^{+} \to Q( \bC \bP^{\infty}_{+}) \to Q (\bS^0)$. Theorem \ref{ktheoryindex} tells us that the maps

$$\bZ \times B \Gamma_{\infty}^{+} \to \bZ \times B \Gl_{\infty}(\bZ)^{+}$$

given by the action on the first homology and

$$\bZ \times B \Gamma_{\infty}^{+} \to \Omega^{\infty } \bG_{-2}^{SO} \to Q (\bC \bP^{\infty}_{+} )\to
Q (\bS^0 )\to \bZ \times B \Gl_{\infty}(\bZ)^{+}$$

are homotopic.\\
Since the $e$-invariant $\pi_3 (Q (\bS^0)) \to \bQ / \bZ$ is injective, the composition

$$ \pi_3 (B \Gamma_{\infty}^{+}) \to \pi_3 (Q (\bS^0)) \to \bQ / \bZ$$

is also injective and it agrees with the composition

$$ \pi_3 (B \Gamma_{\infty}^{+}) \to K_3 (\bZ) \to \bQ / \bZ$$

of the representation map with the $e$-invariant. This shows that the order of the element in $\pi_3 (B \Gamma_{\infty}^{+})$ with $e$-invariant
$1/12$ constructed in this chapter has order $12$ (and not $24$).

\pagebreak

\appendix

\section{Appendix on groupoids}\label{appendiy}

In this appendix, we give the proof of two abstract statements about classifying spaces of topological groupoids.

\begin{lem}\label{groupoid}
%\begin{special}{Lemma \ref{groupoid}}
Let $G$ be a topological groupoid with a discrete object set. Let $R$ be a system of representatives for the set $\Ob(G) / \cong$. Then the functor

$$F:\coprod_{x \in R}  \Aut (x) \to G$$

induces a homotopy equivalence of classifying spaces:

$$\coprod_{x \in R}  B\Aut (x) \to BG.$$
%\end{special}
\end{lem}

\textbf{Proof:} We show that $F$ is an equivalence of topological categories by defining an inverse functor $H$.\\
For any $y \in \Ob (G)$, there is a unique object $H(y) \in R$ and an isomorphism $a_{y} :H( y) \to y$. We define $H$ by sending

$$y \mapsto x; (g:y \to z) \mapsto a_{z}^{-1} \circ g \circ a_y.$$

We write composition in a category contravariantly, pretending that the morphisms are maps between sets.
Since the object set is discrete, $H$ is a continuous functor. It is trivial to check that $H$ is right adjoint (it is even an equivalence of categories). Thus $F$ and $H$ induce mutually inverse homotopy equivalences. \qed\\

Let $H$ be a topological groupoid (with discrete object set); $G$ a topological group with unit $1$ and $\phi: H \to G$ a continuous morphism of groupoids. The \emph{kernel} of $\phi $ is the groupoid whose objects are the objects of $H$ and whose morphisms are the elements $h \in \Mor (H)$ with $\phi(h)=1$.
We say that a morphism of groupoids (alias functor) $\Phi: H \to G$ is \emph{surjective} if for any object $x \in \Ob (H)$ and any $g \in \Mor(G) =G$, there exists $y \in \Ob (H)$ and $h \in \Mor_H (x;y)$ with $\Phi(h)=g$.\\

\begin{prop}\label{groupoid2}
Let $G$, $H$ and $\phi$ as before and let $K= \ker \phi$. Assume that $\phi$ is surjective in the sense that for any $g \in G$, there exists $h \in \Mor H$ with $\phi(h)=g$. Assume further that the map $\Mor (H) \to G$ is a Serre fibration.\\
Then the canonical map from $BK$ to the homotopy fiber of $B \phi: BH \to BG$ is a homotopy equivalence.
\end{prop}

The proof is a bit longer. First, we prove the proposition under the assumption that $G$ and $H$ are \emph{discrete}. We can use Quillen's Theorem B (\cite{QuiKT}, p. 97) to do the job. We state it in greater generality than actually needed. Let $F: \cA \to \cB$ be a functor between small categories. For an object $Y \in \Ob (\cB)$, we define the \emph{right fiber} $Y /F$ to be the category whose objects are the pairs $(X,b)$; $X$ an object of $\cA$ and $b:Y \to F(X)$ a morphism in $\cB$. A morphism $(X;b) \to (X^{\prime}, b^{\prime}) $ in $Y/F$ is a morphism $a : X \to X^{\prime}$ in $\cA$ such that $b^{\prime} = F(a) \circ b$.\\
Any morphism $\beta: Y \to Y^{\prime}$ in $\cB$ determines a functor $\beta^{*} : Y^{\prime} /F \to Y/F$ by sending $(X^{\prime},b^{\prime})$ to $(X^{\prime},b^{\prime} \circ \beta)$. This functor is called the \emph{transition map}. There is an obvious forgetful functor $v:Y/F \to \cA$. There is also an obvious functor $u:Y/F \to Y / \cB$, where the latter is the category of all objects over $Y $ in $\cB$. The classifying space of $Y/B$ is contractible, since it has a terminal object. Quillen's Theorem B states that if all transition maps induce homotopy equivalences on the level of classifying spaces, then the following commutative diagram is homotopy cartesian:

$$\xymatrix{
B(Y/F) \ar[r]^{v} \ar[d]^{u} & B \cA \ar[d]^{F}\\
B(Y/\cB) \ar[r]  & B \cB.\\
}
$$

In other words, $B(Y/F)$ is the homotopy fiber of $BF : B \cA \to B \cB$. If $\cB=G$ is a group, then every morphism in $\cB$ is an isomorphism and thus every transition map is a homotopy equivalence after realization. Thus the assumption of Quillen's Theorem B are satisfied. It remains to identify $B(Y/F)$ with $BK$. Let $F: H \to G$ be a morphism of groupoids, let $G$ be a group and let $\phi$ be surjective. Let $K$ be the kernel of $\phi$.\\
Let $*$ be the object of $\cB$. We define a functor $\cF : K \to * / \phi$:\\
The objects of $K$ are the objects of $H$, and $\cF$ sends an object $a$ to the object $(a,1) \in \Ob (* / \phi)$. A morphism $f: a \to b$ in $H$ is in $K$ if and only if $\phi(f)=1$; and $  \cF$ sends it to the morphism $f: (a,1) \to (b,1)$. It is readily checked that $\cF$ is essentially surjective (every object in $*/\phi$ is isomorphic to one of the form $\cF(a)$) and that $\cF$ is bijective on morphism sets. Under these circumstances, $\cF$ induces a homotopy equivalence on classifying spaces. Note that the construction of $\cF$ is natural in $G,H$ and $\phi$.\\
This finishes the proof of the proposition in the case that $H$ and $G$ are discrete.\\
In the general case, we want to apply a generalized version of Quillen's Theorem B, due to Waldhausen (\cite{Wald}). It is stated for simplicial categories. Recall that a \emph{simplicial category} is a simplicial object in the category of small categories and functors. The standard example is the singular simplicial category $S_{\bullet} \cC$ of a topological category $\cC$. The objects of the category $S_n \cC$ in the $n$th degree are the continuous maps from the standard $n$-simplex $\Delta^n$ to the object space of $\cC$, the morphisms of $S_n \cC$ are the continuous maps $\Delta^n \to \Mor (\cC)$, with pointwise composition.\\
If $\cA_{\bullet}$ is a simplicial category, we can form its \emph{nerve} $N_{\bullet} \cA_{\bullet}$, which is the bisimplicial set $([n], [m]) \mapsto N_n \cA_m:= Funct ([n] ; \cA_m)$, where $[n]$ is the partially ordered set $ 0 < 1 < \ldots n$, viewed as a category.\\
For any bisimplicial set $([n], [m]) \mapsto X_{n,m}$, there is a canonical homeomorphism between the iterated realizations:

$$| [n] \mapsto |[m] \mapsto X_{n,m}|| \cong | [m] \mapsto |[n] \mapsto X_{n,m}||,$$

see \cite{QuiKT}, p. 94. We denote both spaces with the symbol $B|\cA_{\bullet}|$. In the case that $X_{n,m} = N_n S_m \cC$ for a topological
category $\cC$, this shows that

$$|[m] \mapsto | [n] \mapsto N_n S_m \cC|| \cong |[m] \mapsto B(S_m \cC)|$$

and

$$|[n] \mapsto | [m] \mapsto N_n S_m \cC|| \cong |[n] \mapsto |S_{\bullet} N_n \cC|$$

are homeomorphic. Because for any space $X$; $X$ and $|S_{\bullet} X|$ are weakly homotopy equivalent, we see that

$$|[n] \mapsto |S_{\bullet} N_n \cC| \simeq |[n] \mapsto N_n \cC| = B \cC.$$

All homeomorphisms and homotopy equivalences occurring in this discussion are natural.\\
It follows that we can replace $H$, $G$ and $\phi$ in the proposition by their singular simplicial counterpart. More precisely, we study
$S_{\bullet} \phi: S_{\bullet} H \to S_{\bullet} G$. It is a simplicial functor. Let us state Waldhausens improvement of Quillen's Theorem B.
Assume that $F_{\bullet}: \cA_{\bullet} \to \cB_{\bullet}$ is a simplicial functor. We make the additional assumption that any object in
$\cB_{\bullet}$ is isomorphic to a zero-dimensional one. That means the following: For any $y \in \Ob(\cB_q)$, there exist $z \in \Ob (\cB_0)$
and an isomorphism $f: y \to \epsilon^{*} z$; where $\epsilon: [q] \to[0]$ is the unique map. Then the right fiber over $y \in \cB_0$ is the
simplicial category

$$q \mapsto F_n / \epsilon^{*} y.$$

If all transition maps of the right fibers which are induced by morphisms in $\cB_0$ are homotopy equivalences, then the analogous statement to Quillen's Theorem B is true. There is a version without the assumption about the zero-dimensional objects which we do not need here.\\
In our case, $\cB = S_{\bullet} G$ is a simplicial group and thus it has only one object in any degree. Thus the additional assumption is satisfied. The transition maps are all homotopy equivalences because all morphisms in $\cB_0$ are isomorphisms.\\
The assumption that our original $\phi$ was a surjective Serre fibration implies that $S_n \phi: S_n H \to S_n G$ is a surjective homomorphism of discrete groupoids with kernel $S_n K $. We conclude that we have a homotopy equivalence

$$BK = B| [n] \mapsto S_n K| \to B| [n] \mapsto F_n / \epsilon^{*} * | = B(F_{\bullet} / *)$$

by the consideration for the discrete case. By the improved Theorem B, it follows that

$$B(F/ *) \to \hofib B|S_{\bullet} \phi|$$

is a homotopy equivalence, which finishes the proof of the proposition.\qed\\

\pagebreak

\section{List of notations}\label{listof}

$\At(M, \sigma)$ --- the Atiyah invariant of the spin structure $\sigma$ on the surface $M$ (\ref{atiyahinv})

$B_n$ --- $n$th Bernoulli number, \ref{bernoull}

$BG$ ---  classifying space of the group $G$

$C(X;Y)$ --- space of continuous maps $X \to Y$

$\chi_A (V)$ --- Euler class in the generalized cohomology theory $A$ of the $A$-oriented vector bundle $V$

$\delbar_V$ --- Cauchy-Riemann operator on a surface bundle twisted with a vector bundle $V$

$\den(q)$ --- denominator of the rational number $q$

$\Diff(M)$ --- usually the group of \emph{orientation-preserving} diffeomorphisms of the manifold $M$

$\Diff_0 (M)$ --- the group of diffeomorphisms of $M$ homotopic to the identity

$\Diff(M, \sigma)$ --- subgroup of $\Diff(M)$, consisting of diffeomorphisms fixing the spin structure $\sigma$ (\ref{spindiffff})

$e$ --- the $e$-invariant (only in chapter \ref{icosahedral})

$e(V)$ --- Euler class in ordinary cohomology of the oriented vector bundle $V$

$\chi_A (V)$ --- Euler class in the generalized cohomology theory $A$

$EG \to BG$ --- universal $G$-principal bundle; $E(G;X) \to BG$ --- the bundle $EG \times_G X$; where $X$ is a $G$-space

$G \actson X$ --- the group $G$ acts on the space $X$

$G^{\delta}$ ($\bC^{\delta}$) --- the Lie group $G$ (the ring $\bC$) considered as discrete group (ring)

$\bG_{-d}^{\gF}$ --- the generalized Madsen-Tillmann spectrum (\ref{defmadtil})

$\Gl_{d}^{+} (\bR)$ --- the group of all real matrices with positive determinant

$\spingroup$ --- the twofold connected covering group of $\Gl_{d}^{+} (\bR)$ (\ref{spinvector})

$\Gamma_g, \Gamma_{g,n}, \Gamma_{\infty}$ --- mapping class group of closed surfaces of genus $g$, mapping class group of surfaces of genus $g$ with $n$ boundary components, infinite mapping class group

$\Gamma_{g}^{\epsilon}, \hat{\Gamma_{g}^{\epsilon}}$ --- subgroup of the mapping class group fixing a spin structure of Atiyah invariant $\epsilon$, spin mapping class group of genus $g$ and Atiyah invariant $\epsilon$

$\gamma_n$, $\gamma_n (E)$, $\gamma_n(\pi)$ --- the symplectic cohomology classes (\ref{defnsymp})

$\bold{k}$ --- connective complex topological $K$-theory spectrum

$\kappa_n$, $\kappa_n (\pi)$, $\kappa_n (E)$ --- the Morita-Miller-Mumford classes or MMM-classes (\ref{defnmum})

$KR$ ---  algebraic $K$-theory spectrum of the discrete ring $R$

$\Lambda_v E$ --- vertical cotangent bundle of the manifold bundle $ \pi:E \to B$

$\lambda_n$ --- the analytical Morita-Mumford classes (\ref{analytical})

$\fM_g$ --- moduli space of Riemann surfaces of genus $g$

$\Mor (C)$ --- morphism space of the topological category $C$

$\cO_C$ --- sheaf of local holomorphic functions on the Riemann surface $C$

$\Ob (C)$  --- object space of the topological category $C$

$\Omega^{\gF}_{d}$ ---  the normal $\gF$-bordism group

$\Omega^{\infty} \gE$ --- infinite loop space of the spectrum $\gE$

$\Omega_a \gE$ --- the component of $\Omega^{\infty} \gE$ which corresponds to $a \in \pi_0 (\gE) = \pi_0 (\Omega^{\infty} \gE)$

$\pi:E \to B$ --- default notation for a manifold bundle

$Q(X)$, $Q_0(X)$ --- the infinite loop space of the spectrum $\Sigma^{\infty} X$, the unit component of $Q(X)$

$\cS(M)$ --- the space of all almost complex structures on the surface $M$

$s_n$ --- integral Chern character class (\ref{chclosb})

$\Sp_{2n} (R)$ --- the symplectic group with entries in $R$ (\ref{sympgroup})

$\Sspin(M)$ --- the groupoid of all spin structures on the manifold $M$ (\ref{spiccc})

$\Spiff(M,\sigma)$ --- the group of spin diffeomorphisms of $(M, \sigma)$ (\ref{spinnnnnn})

$\Spiff(M)$ --- the groupoid of all spin structures and spin diffeomorphisms of $M$ (\ref{spinnnnnn})

$\Spiff(M)^{\epsilon}$ --- the groupoid of all spin structures with Atiyah invariant $\epsilon$ and all spin diffeomorphisms (if $M$ is a surface)

$\Sigma_n, \Sigma_{\infty}$ --- the symmetric groups

$\Sigma^{\infty} X$ --- suspension spectrum of a pointed space $X$

$\cT_g$ --- Teichm\"uller space of Riemann surfaces of genus $g$

$T_v E$ --- vertical tangent bundle of the manifold bundle $\pi:E \to B$

$\Th(V)$ --- Thom space of the vector bundle $V$

$\bTh (V)$ --- Thom spectrum of the stable vector bundle $V$ (\ref{defthomspec})

$V_n(\pi)$ --- $n$th Hodge bundle of the surface bundle $\pi:E \to B$.

$X_+$ --- the space $X$ with an additional base point

$X^+$ --- the Quillen plus construction applied to $X$ 

\pagebreak

\pagebreak

\section*{Curriculum vitae}

Name: Johannes Felix Ebert\\
\\
Date and place of birth: 1976, June 16, in Bonn-Bad Godesberg\\
\\
June 1995: Abitur at the Friedrich-Ebert-Gymnasium, Bonn\\
August 1995 until August 1996: Civil service\\
\\
October 1996 until September 1998: Undergraduate studies of mathematics at the Rhein\-ische Friedrich-Wilhelm-Universit\"at Bonn\\
September 1998: Vordiplom\\
\\
October 1998 until May 2003: Graduate studies in mathematics, University of Bonn\\
October 1998 until September 2002: Tutor at the Mathematical Insitute, University of Bonn\\
May 2003: Diplom\\
Title of diploma thesis: "\"Uber den Modulraum mehrfach gerichteter und punktierter Kleinscher Fl\"achen". Supervisor: C.-F. B\"odigheimer\\
\\
July 2003 until June 2006: PhD studies at the University of Bonn, financally supported by the Max-Planck-Institut f\"ur Mathematik Bonn; member of the International Max-Planck Research School "Moduli spaces", member of the Bonner Internationale Graduiertenschule, associate member of the Graduiertenkolleg "Homotopie und Kohomologie."\\

\thispagestyle{empty}

\end{document}